\documentstyle[12pt]{amsart}
\let\text=\mbox
\newcommand{\skipthistext}[1]{}
\newcommand{\C}{{\mathbb C}}
\newcommand{\R}{{\mathbb R}}
\newcommand{\Z}{{\mathbb Z}}
\newcommand{\Q}{{\mathbb Q}}
\renewcommand{\P}{{\mathbb P}}
\newcommand{\s}{{\mathbb S}}
\newcommand{\B}{{\mathbb B}}
\newcommand{\I}{{\mathbb I}}
\newcommand{\h}{{\mathbb H}}

\renewcommand{\pf}{\noindent {\bf Proof} \hspace{2mm}}
              
              \newcommand{\M}{{\cal M}}
              
              \newcommand{\U}{{\cal U}}

              \newcommand{\E}{{\cal E}}
              \newcommand{\F}{{\cal F}}

              \newcommand{\G}{{\cal G}}

\newcommand{\abs}[1]{\centerline{\sc Abstract}
\vspace{3mm}
\centerline{\parbox{150mm}{\small #1}}}
\newcounter{thm}

\def\thmnumber{
\addtocounter{thm}{1}\if\arabic{section}0\relax
\else\arabic{section}.\fi%
\if\arabic{subsection}0\relax
\else\arabic{subsection}.\fi%
\arabic{thm}
\hspace{2mm}}

\def\setthmnumber{\setcounter{thm}{0}}

\newenvironment{thm}{
\vspace{3mm}\par\noindent
{\bf Theorem \thmnumber} \it}
{\vspace{3mm}\par\rm}

\newenvironment{conj}{
\vspace{3mm}\par\noindent
{\bf Conjecture \thmnumber} \it}
{\vspace{3mm}\par\rm}

\newenvironment{prob}{
\vspace{3mm}\par\noindent
{\bf Problem \thmnumber} \it}
{\vspace{3mm}\par\rm}

\newenvironment{lem}{
\vspace{3mm}\par\noindent
{\bf Lemma \thmnumber} \it}
{\vspace{3mm}\par\rm}

\newenvironment{prop}{
\vspace{3mm}\par\noindent
{\bf Proposition \thmnumber}\it}
{\vspace{3mm}\par\rm}

\newenvironment{re}{
\vspace{3mm}\par\noindent
{\bf Remark \thmnumber}\it}
{\vspace{3mm}\par\rm}

\newcommand{\sectioni}[1]{\setthmnumber\setcounter{subsection}{0}
\section{#1}}

\renewcommand{\subsection}[1]{%
\vspace{7mm}\par\noindent%
\addtocounter{subsection}{1}%
{\sc\arabic{section}.\arabic{subsection} \ {#1}}
\vspace{4mm}\par\noindent\setthmnumber}

\begin{document}
\title{Smooth $s$-cobordisms of elliptic $3$-manifolds}
\author{Weimin Chen}
\date{\today}
\maketitle

\abs{The main result of this paper states that a symplectic
$s$-cobordism of elliptic $3$-manifolds is diffeomorphic to a product
(assuming a canonical contact structure on the boundary). Based on this 
theorem, we conjecture that a smooth $s$-cobordism of elliptic $3$-manifolds 
is smoothly a product if its universal cover is smoothly a product. We 
explain how the conjecture fits naturally into the program of Taubes of 
constructing symplectic structures on an oriented smooth $4$-manifold with 
$b_2^{+}\geq 1$ from generic self-dual harmonic forms. The paper also contains 
an auxiliary result of independent interest, which generalizes Taubes' theorem 
``$SW\Rightarrow Gr$'' to the case of symplectic $4$-orbifolds. 
}

\sectioni{Introduction: conjecture and main result}

One of the fundamental results in topology is the so-called
$s$-cobordism theorem, which allows one to convert
topological problems into questions of algebra and homotopy theory.
This theorem says that if $W$ is a compact $(n+1)$-dimensional manifold
with boundary the disjoint union of manifolds $Y_1$ and $Y_2$, then when
$n\geq 5$, $W$ is diffeomorphic, piecewise linearly homeomorphic, or
homeomorphic, depending on the category, to the product $Y_1\times [0,1]$,
provided that the inclusion of each boundary component into $W$ is a 
homotopy equivalence and a certain algebraic invariant $\tau(W;Y_1)\in
Wh(\pi_1(W))$, the Whitehead torsion, vanishes. (Such a $W$ is called an
$s$-cobordism from $Y_1$ to $Y_2$; when $\pi_1(W)=\{1\}$, the theorem is 
called the $h$-cobordism theorem, first proved by Smale.) 
Note that the $s$-cobordism theorem is trivial for the dimensions where 
$n\leq 1$. However, great effort has been made to understand the remaining 
cases where $n=2,3$ or $4$, and the status of the $s$-cobordism theorem in 
these dimensions, for each different category, reflects the fundamental
distinction between topology of low-dimensional manifolds and that of
the higher dimensional ones.

For the case of $n=2$, the $s$-cobordism theorem is equivalent to the
original Poincar\'{e} conjecture, which asserts that a closed, simply
connected $3$-manifold is homeomorphic to the $3$-sphere
(cf. e.g. \cite{He,Li,Ru}). For $n=4$, work of Freedman \cite{F} yielded
a topological $s$-cobordism theorem for $W$ with a relatively small
fundamental group, e.g. finite or polycyclic. On the other hand, Donaldson
\cite{D} showed that the $h$-cobordism theorem fails in this dimension for
the smooth (and equivalently, the piecewise linear) category.

This paper is concerned with $4$-dimensional $s$-cobordisms with boundary
components homeomorphic to elliptic $3$-manifolds. (An elliptic $3$-manifold
is one which is homeomorphic to $\s^3/G$ for some finite subgroup
$G\subset SO(4)$ acting freely on $\s^3$.) Building on the aforementioned
work of Freedman, the classification of topological $s$-cobordisms of
elliptic $3$-manifolds up to orientation-preserving homeomorphisms was
completed in a series of papers by Cappell and Shaneson \cite{CS1, CS2},
and Kwasik and Schultz \cite{KS1, KS2}. Their results showed that for each
elliptic $3$-manifold, the set of distinct topological $s$-cobordisms is
finite, and is readily determined from the fundamental group of the
$3$-manifold. In particular, there are topologically nontrivial (i.e.
non-product), orientable $s$-cobordisms in dimension four\footnote
{Existence of nontrivial, non-orientable $4$-dimensional $s$-cobordisms, 
which is of a different nature, had been fairly understood, cf. \cite{MS}, 
also \cite{Kw} for a concrete example.}, and the nontriviality of these 
$s$-cobordisms is evidently related to the fundamental group of the 
$3$-manifold. On the other hand, not much is known in the smooth category. 
Note that the construction of the aforementioned nontrivial $s$-cobordisms 
involves surgery on some topologically embedded $2$-spheres, and it is 
generally a difficult problem to determine whether these $2$-spheres are 
smoothly embedded. In particular, it is not known whether these nontrivial 
$s$-cobordisms are smoothable or not. As for smooth $s$-cobordisms obtained 
from constructions other than taking a product, examples can be found in 
Cappell and Shaneson \cite{CS3, CS4} (compare also \cite{AR}), where the
authors exhibited a family of smooth $s$-cobordisms $W_r$ from $\s^3/Q_r$ 
to itself, with
$$
Q_r=\{x,y\mid x^2=y^{2^r}=(xy)^2=-1\}
$$
being the group of generalized quaternions of order $2^{r+2}$ (note that
$Q_r$ with $r=1$ is the group of order $8$ generated by the quaternions
$i,j$). It has been an open question, only until recently, as whether any
of these $s$-cobordisms or their finite covers is smoothly nontrivial.
In \cite{Ak}, Akbulut showed that the universal cover of $W_r$ with $r=1$
is smoothly a product. Despite this result, however, the following general
questions have remained untouched.

\begin{itemize}
\item [{(1)}] Are there any exotic smooth structures on a 
              trivial 4-dimensional $s$-cobordism?
\item [{(2)}] Is any of the nontrivial topological 4-dimensional
              $s$-cobordisms smoothable?
\end{itemize}

In this paper, we propose a program for understanding smooth $s$-cobordisms
of elliptic $3$-manifolds. At the heart of this program is the following
conjecture, which particularly suggests that in the smooth category,
any nontrivial $s$-cobordism (should there exist any) will have nothing
to do with the fundamental group of the $3$-manifold.

\begin{conj}
A smooth $s$-cobordism of an elliptic $3$-manifold to itself is smoothly
a product if and only if its universal cover is smoothly a product.
\end{conj}

We propose two steps toward Conjecture 1.1, and undertake the first one
in this paper.

In order to describe the first step, we recall some relevant definitions
from symplectic and contact topology. Let $Y$ be a $3$-manifold. A contact
structure on $Y$ is a distribution of tangent planes $\xi\subset TY$ where
$\xi=\ker\alpha$ for a $1$-form $\alpha$ such that $\alpha\wedge d\alpha$
is a volume form on $Y$. Note that the contact manifold $(Y,\xi)$ has a
canonical orientation defined by the volume form $\alpha\wedge d\alpha$.
Let $(Y_i,\xi_i)$, $i=1,2$, be two contact $3$-manifolds given with the
canonical orientation. A symplectic cobordism from $(Y_1,\xi_1)$ to
$(Y_2,\xi_2)$ is a symplectic $4$-manifold with boundary $(W,\omega)$
such that $\partial W=Y_2-Y_1$ (here $W$ is given with the canonical 
orientation defined by $\omega\wedge\omega$), and that there exists a normal 
vectorfield $v$ in a neighborhood of $\partial W$ where $L_v\omega=\omega$, 
and $\xi_i=\ker (i_v\omega|_{Y_i})$ for $i=1,2$. 

Notice that every elliptic $3$-manifold is diffeomorphic to $\s^3/G$,
where $\s^3\subset\C^2$, for a finite subgroup $G\subset U(2)$ acting
freely on $\s^3$ (we shall explain later in this section). The $3$-manifold 
$\s^3/G$ has a canonical 
contact structure, i.e., the descendant of the distribution of complex lines 
on $\s^3$ under the map $\s^3\rightarrow \s^3/G$. Furthermore, the canonical 
orientation from the contact structure coincides with the one induced from 
the canonical orientation on $\s^3$. With the preceding understood, the
following theorem is the main result of this paper.

\begin{thm}
A symplectic $s$-cobordism from an elliptic $3$-manifold $\s^3/G$ to 
itself is diffeomorphic to a product. Here $G$ is a subgroup of $U(2)$ and
$\s^3/G$ is given with the canonical contact structure.
\end{thm}

Thus, in order to prove Conjecture 1.1, it suffices to show, which
is the second step, that a smooth $s$-cobordism of an elliptic
$3$-manifold to itself is symplectic if its universal cover is smoothly
a product. We shall explain next how this step fits naturally into
Taubes' program of constructing symplectic structures on an oriented
smooth $4$-manifold with $b_2^{+}\geq 1$ from generic self-dual harmonic
forms on the $4$-manifold, cf. \cite{T4}.

The starting point of Taubes' program is the observation that on
an oriented smooth $4$-manifold with $b_2^{+}\geq 1$, a self-dual harmonic
form for a generic Riemannian metric has only regular zeroes, which
consist of a disjoint union of embedded circles in the $4$-manifold. In the
complement of the zero set, the $2$-form defines a symplectic structure,
and furthermore, given the almost complex structure in the complement
which is canonically defined by the metric and the self-dual $2$-form,
Taubes showed that nontriviality of the Seiberg-Witten invariant of
the $4$-manifold implies existence of pseudoholomorphic subvarieties
in the complement which homologically bound the zero set. Having
said this, the basic idea of the program is to cancel the zeroes of
the self-dual $2$-form to obtain a symplectic form on the $4$-manifold,
by modifying it in a neighborhood of the pseudoholomorphic subvarieties.

As illustrated in \cite{T4}, it is instructive to look at the
case where the $4$-manifold is $S^1\times M^3$, the
product of a circle with a closed, oriented $3$-manifold. Let $\alpha$
be a harmonic $1$-form on $M^3$ with integral periods, which has
only regular zeroes for a generic metric. In that case, $\alpha=
df$ for some circle-valued harmonic Morse function $f$ on $M^3$.
Given with such a $1$-form $\alpha$, one can define a self-dual harmonic 
form $\omega$ on $S^1\times M^3$ for the product metric 
by 
$$
\omega=dt\wedge\alpha+\ast_3\alpha,
$$
where $t$ is the coordinate function on the $S^1$ factor, and $\ast_3$
is the Hodge-star operator on $M^3$. The zero set of $\omega$ is
regular, and can be easily identified with $\sqcup_{\{p|df(p)=0\}}
S^1\times\{p\}$. Moreover, the pseudoholomorphic subvarieties in
this case are nothing but the embedded tori or cylinders in
$S^1\times M^3$ which are of the form $S^1\times\gamma$, where
$\gamma$ is an orbit of the gradient flow of the Morse function $f$,
either closed or connecting two critical points of $f$. With
these understood, Taubes' program for the $4$-manifold $S^1\times M^3$,
if done $S^1$-equivariantly, is nothing but to cancel all critical
points of a circle-valued Morse function on $M^3$ to make a fibration
$M^3\rightarrow S^1$. It is well-known that there are substantial 
difficulties in canceling critical points in dimension $3$. This seems
to suggest that in general one may expect similar difficulties in
implementing Taubes' program in dimension $4$ as well.

With the preceding understood, our philosophy is to consider Taubes' program
in a more restricted context where the $4$-manifold is already symplectic,
so that one may use the existing symplectic structure as a reference point
to guide the cancellation of the zeroes of a self-dual harmonic form.
For a model of this consideration, we look at the case of
$S^1\times M^3$ where $M^3$ is fibered over $S^1$ with fibration $f_0:
M^3\rightarrow S^1$. Suppose the circle-valued Morse function $f$ is homotopic
to $f_0$. Then a generic path of functions from $f$ to $f_0$ will provide 
a guide to cancel the critical points of $f$ through a sequence of 
birth/death of critical points of Morse functions
on $M^3$. Note that Taubes' program in this restricted sense will not
help to solve the existence problem of symplectic structures in general,
but it may be used to construct symplectic structures with certain special
features, e.g., equivariant symplectic structures in the presence of
symmetry. (Note that this last point, when applied to the case of
$S^1\times M^3$, is related to the following conjecture which still remains
open: If $S^1\times M^3$ is symplectic, $M^3$ must be fibered over $S^1$.
Under some stronger conditions, the conjecture was verified in \cite{CM}
through a different approach.) Now we consider

\begin{prob}
Let $G$ be a finite group acting smoothly on $\C\P^2$ which has an
isolated fixed point $p$ and an invariant embedded $2$-sphere $S$
disjoint from $p$, such that $S$ is symplectic with respect to the
K\"{a}hler structure $\omega_0$ and generates the second homology.
Suppose $\omega$ is a $G$-equivariant, self-dual harmonic form which
vanishes transversely in the complement of $S$ and $p$. Modify $\omega$
in the sense of Taubes \cite{T4}, away from $S$ and $p$, to construct a
$G$-equivariant symplectic form on $\C\P^2$.
\end{prob}

A positive solution to Problem 1.3 will confirm Conjecture 1.1,
as we shall explain next.

\vspace{2mm}

Let $W$ be a $s$-cobordism with boundary the disjoint union of elliptic
$3$-manifolds $Y_1,Y_2$. Clearly $W$ is orientable. We note first that
for any orientations on $Y_1,Y_2$ induced from an orientation on $W$,
there exists an orientation-preserving homeomorphism from $Y_1$ to $Y_2$.
Such a homeomorphism may be obtained as follows. Let $h:Y_1\rightarrow Y_2$
be the simple homotopy equivalence induced by the $s$-cobordism $W$.
Then $h$ is easily seen to be orientation-preserving for any induced
orientations on $Y_1,Y_2$. On the other hand, as a simple homotopy
equivalence between geometric $3$-manifolds, $h$ is homotopic to a
homeomorphism $\hat{h}:Y_1\rightarrow Y_2$ (cf. \cite{Tu}, and for a
proof, cf. \cite{KS3}), which is clearly also orientation-preserving.

Next we recall the fact that every finite subgroup $G\subset SO(4)$
which acts freely on $\s^3$ is conjugate in $O(4)$ to a subgroup
of $U(2)$. In order to understand this, we fix an identification
$\R^4=\C^2=\h$, where $\C^2$ is identified with the space of quaternions
$\h$ as follows
$$
(z_1,z_2)\mapsto z_1+z_2 j.
$$
Consequently, the space of unit quaternions $\s^3$ is canonically identified
with $SU(2)$. Consider the homomorphism $\phi:\s^3\times\s^3\rightarrow
SO(4)$ which is defined such that for any $(q_1,q_2)\in \s^3\times\s^3$,
$\phi(q_1,q_2)$ is the element of $SO(4)$ that sends $x\in\R^4=\h$ to
$q_1xq_2^{-1}\in\h=\R^4$. It is easily seen that $\phi$ is surjective with
$\ker \phi=\{(1,1), (-1,-1)\}$, where we note that the center of $\s^3$
consists of $\{\pm 1\}$. Let $\s^1\subset\s^3$ be the subset consisting
of elements of the form $(z,0)\in\C^2=\h$. Then it is easily seen that a
subgroup of $SO(4)$ acts complex linearly on $\C^2=\h$ if it lies in the image
$\phi(\s^1\times\s^3)$. Note on the other hand that one can switch the
two factors of $\s^3$ in $\s^3\times\s^3$ by an element of $O(4)$ which
sends $x\in\R^4=\h$ to its conjugate $\bar{x}\in\h=\R^4$. With these
understood, it suffices to note that every finite subgroup of $SO(4)$
which acts freely on $\s^3$ is conjugate in $SO(4)$ to a subgroup of either
$\phi(\s^1\times\s^3)$ or $\phi(\s^3\times\s^1)$, cf. Theorem 4.10
in \cite{Sct}.

Now suppose $W$ is a smooth $s$-cobordism of elliptic $3$-manifolds
$Y_1,Y_2$. By combining the aforementioned two facts, it is easily seen
that for any fixed orientation on $W$, one can choose a normal
orientation near $\partial W$ such that with respect to the induced
orientations on $Y_1,Y_2$, there exist orientation-preserving
diffeomorphisms $f_1:Y_1\rightarrow \s^3/G$, $f_2:Y_2\rightarrow \s^3/G$,
where $\s^3\subset\C^2$ and $G$ is a finite subgroup of $U(2)$ acting freely 
on $\s^3$, and $\s^3/G$ is given with the canonical orientation. Call the
regular neighborhood of a component of $\partial W$ the positive end
(resp. negative end) of $W$ if it is identified by an orientation-preserving
map with $(-1,0]\times (\s^3/G)$ (resp. $[0,1)\times (\s^3/G)$).

\begin{lem}
By further applying a conjugation in $SO(4)$ to the $G$-action on the
negative end if necessary, one can fix an identification $\R^4=\C^2=\h$
and regard $G$ canonically as a subgroup of $U(2)$, such that there exists 
a $2$-form $\omega$ on $W$, which is self-dual and harmonic with
respect to some Riemannian metric, and has the following properties.
\begin{itemize}
\item [{(1)}] There are constants $\lambda_{+}>\lambda_{-}>0$ for which
$\omega=\lambda_{+}\omega_0$ on the positive end and $\omega
=\lambda_{-}\omega_0$ on the negative end. Here $\omega_0$ is the
standard symplectic form on $\C^2/G$ {\em (}i.e. the descendant of
$\frac{\sqrt{-1}}{2}\sum_{i=1}^2 dz_i\wedge d\bar{z}_i$ on $\C^2$ {\em )},
and the two ends of $W$ are identified via $f_1,f_2$ to the corresponding
neighborhoods of $\s^3/G$ in $\C^2/G$.
\item [{(2)}] The $2$-form $\omega$ has only regular zeroes.
\end{itemize}
\end{lem}

With this understood, let $\widetilde{W}$ be the universal
cover of $W$ and $\tilde{\omega}$ be the pull-back of $\omega$ to
$\widetilde{W}$. Then $\partial\widetilde{W}=\s^3\sqcup\s^3$ and
$\tilde{\omega}$ equals a constant multiple of the standard symplectic
form on $\C^2$ near $\partial\widetilde{W}$. In particular, both ends
of $\widetilde{W}$ are of contact type with respect to $\tilde{\omega}$,
with one end convex and one end concave. As $\tilde{\omega}$ is invariant
under the Hopf fibration $\s^3\rightarrow\s^2$, we can close up
$\widetilde{W}$ by collapsing each fiber of the Hopf fibration on the
convex end and capping off the concave end with the standard symplectic
$4$-ball $\B^4$. The resulting smooth $4$-manifold $X$ is a homotopy 
$\C\P^2$, with a smoothly embedded $2$-sphere $S$ representing a generator 
of $H_2(X;\Z)$ and a self-dual harmonic form $\omega^\dag$, which
is $G$-equivariant with respect to the obvious $G$-action on $X$, 
has only regular zeroes, and obeys $\omega^\dag |_{S}>0$.

If furthermore, the universal cover $\widetilde{W}$ is smoothly a product,
then the $4$-manifold $X$ is diffeomorphic to $\C\P^2$, with a symplectic
form $\omega_0^\dag$ such that $\omega_0^\dag |_{S}>0$.
It is clear now that a positive solution to Problem 1.3, when applied to
$X$, would yield a symplectic structure on $W$, making it into a symplectic
$s$-cobordism. By Theorem 1.2, $W$ is smoothly a product.

\begin{re}{\em
We add a remark here about the smooth $s$-cobordism $W_r$ with $r=1$ in 
the examples of Cappell and Shaneson \cite{CS3, CS4}. According to
Akbulut's theorem in \cite{Ak}, the universal cover of $W_r$ with $r=1$
is smoothly a product. Thus Conjecture 1.1 suggests that the $s$-cobordism
itself is smoothly a product. It would be interesting to find out by
direct means such as in Akbulut \cite{Ak} as whether $W_r$ with $r=1$
has an exotic smooth structure (it is known to be topologically a product
by the classification of Cappell and Shaneson \cite{CS1, CS2}). 
}
\end{re}

We now turn to the technical aspect of this paper. 

\vspace{1.5mm}

Despite the tremendous progress over the last two decades, topology of
smooth $4$-manifolds is still largely obscure as far as classification
is concerned. In particular, there is lack of effective methods for 
determining the diffeomorphism type of a $4$-manifold in a given homotopy
class. However, in some rare cases and under an additional assumption that the 
$4$-manifold is symplectic, Gromov in \cite{Gr} showed us how to recover
the diffeomorphism type using certain moduli space of pseudoholomorphic 
curves (if it is nonempty). Later, Taubes showed in \cite{T3} that in 
Gromov's argument, the existence of pseudoholomorphic curves may be replaced 
by a condition on the Seiberg-Witten invariant of the $4$-manifold, 
which is something more manageable. As a typical example one obtains the 
following theorem.

\begin{thm}{\em (Gromov-Taubes)} \hspace{2mm}
Let $X$ be a symplectic $4$-manifold with the rational homology of $\C\P^2$.
Then $X$ is diffeomorphic to $\C\P^2$ if the Seiberg-Witten invariant of
$X$ at the $0$-chamber vanishes, e.g., if $X$ has a metric of positive
scalar curvature.  
\end{thm}

Our proof of Theorem 1.2, in a nutshell, is based on an orbifold analog
of the above theorem.

More precisely, in order to prove Theorem 1.2 we extend in this paper 
(along with the earlier one \cite{C2}) Gromov's pseudoholomorphic curve 
techniques and Taubes' work on
the Seiberg-Witten invariants of symplectic $4$-manifolds to the
case of $4$-orbifolds. (See \cite{C3} for an exposition.) In particular, 
we prove the following theorem (see Theorem 2.2 for more details).

\begin{thm}{\em (Orbifold Version of Taubes' Theorem ``$SW\Rightarrow Gr$'')} 
\hspace{2mm}
Let $(X,\omega)$ be a symplectic $4$-orbifold. Suppose $E$ is an orbifold
complex line bundle such that the corresponding Seiberg-Witten invariant 
(in Taubes chamber when $b_2^{+}(X)=1$) is nonzero. Then for any 
$\omega$-compatible almost complex structure $J$, the Poincar\'{e} dual
of $c_1(E)$ is represented by $J$-holomorphic curves in $X$. 
\end{thm}

\begin{re}{\em
(1) The proof of Theorem 1.7 follows largely the proof of Taubes in 
\cite{T3}. However, we would like to point out that Taubes' proof involves
in a few places Green's function for the Laplacian $\Delta=d^\ast d$ and
a covering argument by geodesic balls of uniform size. This part of the proof
requires the assumption that the injective radius is uniformly bounded 
from below, which does not generalize to the case of orbifolds
straightforwardly. Some modification or reformulation is needed here. 

(2) The situation of the full version of Taubes' theorem ``$SW=Gr$'' is 
more complicated for $4$-orbifolds. In fact, the proof of ``$SW=Gr$'' relies 
on a regularity result (i.e. embeddedness) of the $J$-holomorphic curves in
Taubes' theorem ``$SW\Rightarrow Gr$'' for a generic almost complex
structure. While this is generally no longer true for $4$-orbifolds, 
how ``regular'' the $J$-holomorphic curves in Theorem 1.7 could be depends,
in a very interesting way, on what types of singularities the 
$4$-orbifold has. We plan to explore this issue on a future occasion.

(3) There are only a few examples of $4$-manifolds which are symplectic 
and have a metric of positive scalar curvature. Hence the $4$-manifold
in Theorem 1.6 rarely occurs. On the other hand, there are numerous 
examples of symplectic $4$-orbifolds which admit positive scalar curvature
metrics. In fact, there is a class of normal complex surfaces, called
log Del Pezzo surfaces, which are K\"{a}hler orbifolds with positive first
Chern class. (By Yau's theorem, these surfaces admit positive Ricci curvature 
metrics.) Unlike their smooth counterpart, log Del Pezzo surfaces occur
in bewildering abundance and complexity (cf. e.g. \cite{KMc}). Recently,
these singular surfaces appeared in the construction of Sasakian-Einstein
metrics on certain $5$-manifolds (including $\s^5$). In particular, the 
following question arose naturally in this context: What are the log Del 
Pezzo surfaces that appear as the quotient space of a fixed point free 
$\s^1$-action on $\s^5$? (See Koll\'{a}r \cite{Kol}.) We believe that the 
techniques developed in this paper would be useful in answering this 
question.
}
\end{re}

We end this section with an outline for the proof of Theorem 1.2. 
First of all, note that the case where the elliptic $3$-manifold is a lens 
space was settled in \cite{C2} using a different method. Hence in this paper, 
we shall only consider the remaining cases, where the elliptic $3$-manifold 
is diffeomorphic to $\s^3/G$ with $G$ being a non-abelian subgroup of $U(2)$.

Let $W$ be a symplectic $s$-cobordism as in Theorem 1.2. Note that
near the boundary the symplectic form on $W$ is standard, and is invariant
under the obvious Seifert fibration on the boundary. We close up $W$
by collapsing each fiber of the Seifert fibration on the convex end
of $W$ and capping off the concave end with a standard symplectic cone
--- a regular neighborhood of $\{0\in\C^2\}/G$ in the orbifold $\C^2/G$
which is given with the standard symplectic structure. The diffeomorphism
type of $W$ can be easily recovered from that of the resulting symplectic
$4$-orbifold $X$. In order to determine the diffeomorphism type of $X$, 
we compare it with the ``standard'' $4$-orbifold $X_0$, which is $\B^4/G$ 
with boundary $\s^3/G$ collapsed along the fibers of the Seifert fibration. 
More concretely, we consider the space $\M$ of pseudoholomorphic 
maps into $X$, which corresponds, under the obvious homotopy equivalence 
$X\rightarrow X_0$, to the family of complex lines in $\B^4/G$ with boundary 
collapsed. Using the pseudoholomorphic curve theory of $4$-orbifolds 
developed in \cite{C2}, one can easily show that $X$ is diffeomorphic to
$X_0$ provided that $\M\neq\emptyset$, from which Theorem 1.2 follows.

Thus the bulk of the argument is devoted to proving that $\M\neq\emptyset$.
We follow the usual strategy of applying Taubes' theorem ``$SW\Rightarrow
Gr$''. More concretely, the proof of $\M\neq\emptyset$ consists of the 
following three steps. 

\begin{itemize}
\item [{(1)}] Construct an orbifold complex line bundle $E$ such that the
homology class of a member of $\M$ is Poincar\'{e} dual to $c_1(E)$. Note 
that when $X$ is smooth, a complex line bundle is determined by its 
Chern class in $H^2(X;\Z)$. This is no longer true for orbifolds.
In particular, we have to construct $E$ by hand, which is given in 
Lemma 3.6. The explicit construction of $E$ is also needed in order to
calculate the contribution of singular points of $X$ to the dimension $d(E)$
of the Seiberg-Witten moduli space corresponding to $E$, which is a
crucial factor in the proof. (See Lemma 3.8.) 
\item [{(2)}] Show that the Seiberg-Witten invariant corresponding to $E$
is zero in the $0$-chamber. This follows from the fact that the $4$-orbifold
$X$ contains a $2$-suborbifold $C_0$ which has a metric of positive 
curvature and generates $H_2(X;\Q)$. Here $C_0$ is the image of the
convex boundary component of $W$ in $X$. (See Lemma 3.7.) 
\item [{(3)}] By a standard wall-crossing argument, with the fact that
$d(E)\geq 0$, the Poincar\'{e} dual of $c_1(E)$ is represented 
by $J$-holomorphic curves by Theorem 1.7. The main issue here is to show 
that there is a component of the $J$-holomorphic curves which is the image 
of a member of $\M$, so that $\M$ is not empty. When $c_1(E)\cdot c_1(E)$ is 
relatively small, one can show that this is indeed the case by using the 
adjunction formula in \cite{C2}. The key observation is that, 
when $c_1(E)\cdot c_1(E)$ is not small, the dimension $d(E)$ of the 
Seiberg-Witten moduli space is also considerably large, so that one 
may break the $J$-holomorphic curves from Theorem 1.7 into smaller components 
by requiring them to pass through a certain number (equaling half of the 
dimension $d(E)$) of specified points. It turns out that one of the resulting 
smaller components is the image of a member of $\M$, so that 
$\M$ is also nonempty in this case. This part of the proof is the content 
of Lemma 3.9, which is the most delicate one, often involving a case by case 
analysis according to the type of the group $G$ in $\s^3/G$. 
\end{itemize}

The organization of this paper is as follows. In \S 2 we briefly go
over the Seiberg-Witten-Taubes theory for $4$-orbifolds, ending with
a statement of the orbifold version of Taubes' theorem,
whose proof is postponed to \S 4. The proof of the main result, Theorem 1.2,
is given in \S 3. There are three appendices. Appendix A contains a brief
review of the index theorem over orbifolds in Kawasaki \cite{Ka}, and a 
calculation
for the dimension of the relevant Seiberg-Witten moduli space. Appendix B
is concerned with some specific form of Green's function for the Laplacian
$\Delta=d^\ast d$ on orbifolds, which is involved in the proof of Taubes'
theorem for $4$-orbifolds. In Appendix C, we give a proof of Lemma 1.4.

\vspace{7mm}

\centerline{\bf Acknowledgments}

\vspace{3mm}

It is a great pleasure to acknowledge the generous help of S\l awomir
Kwasik on $3$-manifold topology, the valuable conversations with him
about $4$-dimensional $s$-cobordism theory, and the comments and suggestions
after reading the draft version of this paper. I am also very grateful to
Cliff Taubes for the useful communications regarding his work
\cite{T1, T2, T3}, and to Dagang Yang for helpful discussions on some
relevant aspects of Riemannian geometry. This work was partially
supported by NSF Grant DMS-0304956.

\vspace{3mm}

\sectioni{The Seiberg-Witten-Taubes theory for $4$-orbifolds}

In this section, we first go over the Seiberg-Witten theory for smooth
$4$-orbifolds, and then we extend Taubes' work \cite{T1,T2,T3} on symplectic
$4$-manifolds to the orbifold setting. The discussion will be brief since the
theory is parallel to the one for smooth 4-manifolds.

Let $X$ be an oriented smooth 4-orbifold. Given any Riemannian metric on
$X$, a $Spin^\C$ structure is an orbifold principal $Spin^\C(4)$ bundle over
$X$ which descends to the orbifold principal $SO(4)$ bundle of oriented
orthonormal frames under the canonical homomorphism $Spin^\C(4)\rightarrow
SO(4)$. There are two associated orbifold $U(2)$ vector bundles (of rank $2$)
$S_{+},S_{-}$ with $\det(S_{+})=\det(S_{-})$, and a Clifford multiplication
which maps $T^\ast X$ into the skew adjoint endomorphisms of
$S_{+}\oplus S_{-}$.

The Seiberg-Witten equations associated to the $Spin^\C$ structure (if
there is one) are equations for a pair $(A,\psi)$, where $A$ is a connection
on $\det(S_{+})$ and $\psi$ is a section of $S_{+}$. Recall that the
Levi-Civita connection together with $A$ defines a covariant derivative
$\nabla_A$ on $S_{+}$. On the other hand, there are two maps $\sigma:
S_{+}\otimes T^\ast X\rightarrow S_{-}$ and $\tau:\mbox{End}(S_{+})
\rightarrow\Lambda_{+}\otimes\C$ induced by the Clifford multiplication,
with the latter being the adjoint of $c_{+}:\Lambda_{+}\rightarrow
\mbox{End}(S_{+})$, where $\Lambda_{+}$ is the orbifold bundle of
self-dual $2$-forms. With these understood, the Seiberg-Witten equations
read
$$
D_A\psi=0 \mbox{  and  } P_{+}F_A=\frac{1}{4}\tau(\psi\otimes\psi^\ast)
+\mu,
$$
where $D_A\equiv\sigma\circ\nabla_A$ is the Dirac operator, $P_{+}:
\Lambda^2 T^\ast X\rightarrow \Lambda_{+}$ is the orthogonal projection,
and $\mu$ is a fixed, imaginary valued, self-dual $2$-form  which
is added in as a perturbation term.

The Seiberg-Witten equations are invariant under the gauge transformations
$(A,\psi)\mapsto (A-2\varphi^{-1}d\varphi,\varphi\psi)$, where
$\varphi\in C^\infty(X;S^1)$ are circle valued smooth functions on $X$.
The space of solutions modulo gauge equivalence, denoted by $M$, is compact,
and when $b^{+}_2(X)\geq 1$ and when it is nonempty, $M$ is a smooth
orientable manifold for a generic choice of $(g,\mu)$, where $g$ is the
Riemannian metric and $\mu$ is the self-dual $2$-form of perturbations.
Furthermore, $M$ contains no classes of reducible solutions (i.e. those
with $\psi\equiv 0$), and if let $M^0$ be the space of solutions modulo
the based gauge group, i.e., those $\varphi\in C^\infty(X;S^1)$ such that
$\varphi(p_0)=1$ for a fixed base point $p_0\in X$, then $M^0\rightarrow M$
defines a principal $S^1$-bundle. Let $c$ be the first Chern class of
$M^0\rightarrow M$, $d=\dim M$, and fix an orientation of $M$. Then
the Seiberg-Witten invariant associated to the $Spin^\C$ structure is
defined as follows.
\begin{itemize}
\item When $d<0$ or $d=2n+1$, the Seiberg-Witten invariant is zero.
\item When $d=0$, the Seiberg-Witten invariant is a signed sum
of the points in $M$.
\item When $d=2n>0$, the Seiberg-Witten invariant equals $c^n[M]$.
\end{itemize}
As in the case of smooth 4-manifolds, the Seiberg-Witten invariant of $X$
is well-defined when $b^{+}_2(X)\geq 2$, depending only on the diffeomorphism
class of $X$ (as orbifolds). Moreover, there is an involution on the set of
$Spin^\C$ structures which preserves the Seiberg-Witten invariant up to
a change of sign. When $b^{+}_2(X)=1$, there is a chamber structure and the
Seiberg-Witten invariant also depends on the chamber where the pair $(g,\mu)$
is in. Moreover, the change of the Seiberg-Witten invariant when crossing
a wall of the chambers can be similarly analyzed as in the smooth 4-manifold
case.

For the purpose of this paper, we need the following wall-crossing formula.
Its proof is identical to the manifold case, hence is omitted, cf. e.g.
\cite{KM}.

\begin{lem}
Suppose $b_1(X)=0$, $b_2^{+}(X)=b_2(X)=1$, and $c_1(S_{+})\neq 0$. Then
there are two chambers for the Seiberg-Witten invariant associated to
the $Spin^\C$ structure $S_{+}\oplus S_{-}$: the $0$-chamber where
$\int_X\sqrt{-1}\mu\wedge\omega_g$ is sufficiently close to $0$, and the
$\infty$-chamber where $\int_X\sqrt{-1}\mu\wedge\omega_g$ is sufficiently
close to $+\infty$. Here $\omega_g$ is a fixed harmonic $2$-form with respect
to the Riemannian metric $g$ such that $c_1(S_{+})\cdot [\omega_g]>0$.
Moreover, if the dimension of the Seiberg-Witten moduli space $M$
{\em (}which is always an even number in this case{\em )} is non-negative,
the Seiberg-Witten invariant changes by $\pm 1$ when considered in the other
chamber.
\end{lem}

Now we focus on the case where $X$ is a symplectic 4-orbifold. Let $\omega$
be a symplectic form on $X$. We orient $X$ by $\omega\wedge\omega$, and fix
an $\omega$-compatible almost complex structure $J$. Then with respect to
the associated Riemannian metric $g=\omega(\cdot,J\cdot)$, $\omega$ is
self-dual with $|\omega|=\sqrt{2}$. The set of $Spin^\C$ structures
on $X$ is nonempty. In fact, the almost complex structure $J$ gives rise to
a canonical $Spin^\C$ structure where the associated orbifold $U(2)$ bundles
are $S_{+}^0=\I\oplus K_X^{-1}$, $S_{-}^0=T^{0,1}X$. Here $\I$ is the trivial
orbifold complex line bundle and $K_X$ is the canonical bundle
$\det(T^{1,0}X)$. Moreover, the set of $Spin^\C$ structures is canonically
identified with the set of orbifold complex line bundles where each
orbifold complex line bundle $E$ corresponds to a $Spin^\C$ structure whose
associated orbifold $U(2)$ bundles are $S_{+}^E=E\oplus (K_X^{-1}\otimes E)$
and $S_{-}^E=T^{0,1}X\otimes E$. The involution on the set of $Spin^\C$
structures which preserves the Seiberg-Witten invariant up to a change of
sign sends $E$ to $K_X\otimes E^{-1}$.

As in the manifold case, there is a canonical (up to gauge equivalence)
connection $A_0$ on $K_X^{-1}=\det(S_{+}^0)$ such that the fact $d\omega=0$
implies that $D_{A_0}u_0=0$ for the section $u_0\equiv 1$ of $\I$ which is
considered as the section $(u_0,0)$ in $S_{+}^0=\I\oplus K_X^{-1}$.
Furthermore, by fixing such an $A_0$, any connection $A$ on
$\det(S_{+}^E)=K_X^{-1}\otimes E^2$ is canonically determined by a
connection $a$ on $E$. With these
understood, there is a distinguished family of the Seiberg-Witten equations
on $X$, which is parametrized by a real number $r>0$ and is for a triple
$(a,\alpha,\beta)$, where in the equtions, the section $\psi$ of $S_{+}^E$ is
written as $\psi=\sqrt{r}(\alpha,\beta)$ and the perturbation term $\mu$
is taken to be $-\sqrt{-1}(4^{-1}r\omega)+P_{+}F_{A_0}$. (Here $\alpha$
is a section of $E$ and $\beta$ a section of $K_X^{-1}\otimes E$.)
Note that when $b_2^{+}(X)=1$, this distinguished family of Seiberg-Witten
equations (with $r\gg 0$) belongs to a specific chamber for the 
Seiberg-Witten invariant.
This particular chamber will be referred to as Taubes chamber.

\vspace{2mm}

The following is the analog of the relevant theorems of Taubes in the
orbifold setting. (Its proof is postponed to \S 4.)

\begin{thm}
Let $(X,\omega)$ be a symplectic $4$-orbifold. Then the following are true.
\begin{itemize}
\item [{(1)}] The Seiberg-Witten invariant associated to the canonical
$Spin^\C$ structure equals $\pm 1$. {\em (}When $b_2^{+}(X)=1$, the
Seiberg-Witten invariant is in Taubes chamber.{\em)}
Moreover, when $b_2^{+}(X)\geq 2$, the Seiberg-Witten invariant corresponding
to the canonical bundle $K_X$ equals $\pm 1$, and for any orbifold complex
line bundle $E$, if the Seiberg-Witten invariant corresponding to $E$ is
nonzero, then $E$ must satisfy
$$
0\leq c_1(E)\cdot [\omega]\leq c_1(K_X)\cdot [\omega],
$$
where $E=\I$ or $E=K_X$ when either equality holds.
\item [{(2)}] Let $E$ be an orbifold complex line bundle. Suppose there
is an unbounded sequence of values for the parameter $r$ such that the
corresponding Seiberg-Witten equations have a solution $(a,\alpha,\beta)$.
Then for any $\omega$-compatible almost complex structure $J$, there are
$J$-holomorphic curves $C_1,C_2,\cdots,C_k$ in $X$ and positive integers
$n_1,n_2,\cdots,n_k$ such that $c_1(E)=\sum_{i=1}^k n_i PD(C_i)$. Moreover,
if a subset $\Omega\subset X$ is contained in $\alpha^{-1}(0)$
throughout, then $\Omega\subset\cup_{i=1}^k C_i$ also.
\end{itemize}
{\em (}Here $PD(C)$ is the Poincar\'{e} dual of the $J$-holomorphic curve $C$.
See \S 3 of \cite{C2} for the definition of $J$-holomorphic curves in
an almost complex $4$-orbifold and the definition of Poincar\'{e}
dual of a $J$-holomorphic curve in the $4$-orbifold.{\em )}
\end{thm}

\begin{re}
{\em
There are two typical sources for the subset $\Omega$ in the theorem.
For the first one, suppose $p\in X$ is an orbifold point such that
the isotropy group at $p$ acts nontrivially on the fiber of $E$ at $p$.
Then $p\in\alpha^{-1}(0)$ for any solution $(a,\alpha,\beta)$, and
consequently $p\in\cup_{i=1}^k C_i$. For the second one, suppose the
Seiberg-Witten invariant corresponding to $E$ is nonzero and the dimension
of the moduli space $M$ is $d=2n>0$. Then for any subset of distinct
$n$ points $p_1,p_2,\cdots,p_n\in X$, and for any value of parameter $r$,
there is a solution $(a,\alpha,\beta)$ such that $\{p_1,p_2,\cdots,p_n\}
\subset \alpha^{-1}(0)$. Consequently, we may require the $J$-holomorphic
curves $C_1,C_2,\cdots,C_k$ in the theorem to contain any given subset of
less than or equal to $n$ points in this circumstance. (The proof goes
as follows. Observe that the map $(a,\alpha,\beta)\mapsto \alpha(p)$
descends to a section $s_p$ of the complex line bundle associated
to the principal $S^1$ bundle $M^0\rightarrow M$, where $M^0$ is the
moduli space of solutions modulo the based gauge group with base point $p$.
Moreover, there are submanifolds $\Sigma_1,\Sigma_2,\cdots,\Sigma_n$ 
of codimension $2$ in $M$ such that each $\Sigma_i$ is Poincar\'{e} dual 
to the first Chern class $c$ of $M^0\rightarrow M$ and is arbitrarily close 
to $s_{p_i}^{-1}(0)$. Now if there were no solution $(a,\alpha,\beta)$ such 
that $\{p_1,p_2,\cdots,p_n\}\subset \alpha^{-1}(0)$, which means that 
$s_{p_1}^{-1}(0)\cap s_{p_2}^{-1}(0)\cap\cdots\cap s_{p_n}^{-1}(0)=\emptyset$, 
then one would have $\Sigma_1\cap \Sigma_2\cap\cdots\cap \Sigma_n=\emptyset$. 
But this contradicts the assumption that $c^n[M]$, the Seiberg-Witten 
invariant, is nonzero.)

\hfill $\Box$
}
\end{re}

\sectioni{Proof of the main result}

We begin by recalling the classification of finite subgroups of
$\mbox{GL}(2,\C)$ without quasi-reflections, which is due to Brieskorn
\cite{Br}. The following is a list of the non-abelian ones up to
conjugations in $\mbox{GL}(2,\C)$.
\begin{itemize}
\item $\langle Z_{2m},Z_{2m};\widetilde{D}_n,\widetilde{D}_n\rangle$,
where $m$ is odd, $n\geq 2$, and $m,n$ are relatively prime.
\item $\langle Z_{4m},Z_{2m};\widetilde{D}_n,C_{2n}\rangle$,
where $m$ is even, $n\geq 2$, and $m,n$ are relatively prime.
\item $\langle Z_{2m},Z_{2m};\widetilde{T},\widetilde{T}\rangle$,
where $m$ and $6$ are relatively prime.
\item $\langle Z_{6m},Z_{2m};\widetilde{T},\widetilde{D}_2\rangle$,
where $m$ is odd and is divisible by $3$.
\item $\langle Z_{2m},Z_{2m};\widetilde{O},\widetilde{O}\rangle$,
where $m$ and $6$ are relatively prime.
\item $\langle Z_{2m},Z_{2m};\widetilde{I},\widetilde{I}\rangle$,
where $m$ and $30$ are relatively prime.
\end{itemize}
Here $Z_k\subset\mbox{ZL}(2,\C)$ is the cyclic subgroup of order $k$ in
the center of $\mbox{GL}(2,\C)$, $C_k\subset SU(2)$ is
the cyclic subgroup of order $k$, and $\widetilde{D}_n,\widetilde{T},
\widetilde{O}$, $\widetilde{I}\subset SU(2)$ are the binary dihedral,
tetrahedral, octahedral and icosahedral groups of order $4n$, $24$,
$48$ and $120$ respectively, which are the double cover of the corresponding
subgroups of $SO(3)$ under the canonical homomorphism $SU(2)\rightarrow
SO(3)$. As for the notation $\langle H_1,N_1;H_2,N_2\rangle$, it stands
for the image under $(h_1,h_2)\mapsto h_1h_2$ of the subgroup of
$H_1\times H_2$, which consists of pairs $(h_1,h_2)$ such that the
classes of $h_1$ and $h_2$ in $H_1/N_1$ and $H_2/N_2$ are equal under
some fixed isomorphism $H_1/N_1\cong H_2/N_2$. (In the present case,
the group does not depend on the isomorphism $H_1/N_1\cong H_2/N_2$,
at least up to conjugations in $\mbox{GL}(2,\C)$.)

We shall assume throughout that the elliptic $3$-manifolds under
consideration are diffeomorphic to $\s^3/G$ for some finite subgroup
$G\subset \mbox{GL}(2,\C)$ listed above. (Note that a finite subgroup
$G\subset U(2)$ acts freely on $\s^3$ if and only if $G$ is a subgroup
of $\mbox{GL}(2,\C)$ containing no quasi-reflections.)

Next we begin by collecting some preliminary but relevant information
about the elliptic $3$-manifold $\s^3/G$. First of all, note that $G$
contains a cyclic subgroup of order $2m$, $Z_{2m}$, which is the subgroup
that preserves each fiber of the Hopf fibration on $\s^3$. Evidently,
the Hopf fibration induces a canonical Seifert fibration on $\s^3/G$,
which can be obtained in two steps as follows. First, quotient $\s^3$ and
the Hopf fibration by the subgroup $Z_{2m}$ to obtain the lens space $L(2m,1)$
and the $\s^1$-fibration on it. Second, quotient $L(2m,1)$ and the
$\s^1$-fibration by $G/Z_{2m}$ to obtain $\s^3/G$ and the Seifert fibration
on $\s^3/G$. It follows immediately that the Euler number of the Seifert
fibration is
$$
e=\frac{|Z_{2m}|}{|G/Z_{2m}|}=\frac{4m^2}{|G|}.
$$
The Seifert fibration has three singular fibers, and the normalized
Seifert invariant
$$
(b,(a_1,b_1),(a_2,b_2),(a_3,b_3)), \mbox{ where }
0<b_i<a_i,\; a_i,b_i \mbox{ relatively prime},
$$
can be determined from the Euler number $e$ and the induced action
of $G/Z_{2m}$ on the base of the $\s^1$-fibration on $L(2m,1)$. We collect
these data in the following list.
\begin{itemize}
\item $\langle Z_{2m},Z_{2m};\widetilde{D}_n,\widetilde{D}_n\rangle$,
$\langle Z_{4m},Z_{2m};\widetilde{D}_n,C_{2n}\rangle$:
$(a_1,a_2,a_3)=(2,2,n)$, and $b,b_1,b_2$ and $b_3$ are given by
$b_1=b_2=1$, $m=(b+1)n+b_3$.
\item $\langle Z_{2m},Z_{2m};\widetilde{T},\widetilde{T}\rangle$,
$\langle Z_{6m},Z_{2m};\widetilde{T},\widetilde{D}_2\rangle$:
$(a_1,a_2,a_3)=(2,3,3)$, and $b,b_1,b_2$ and $b_3$ are
given by $b_1=1$, $m=6b+3+2(b_2+b_3)$.
\item $\langle Z_{2m},Z_{2m};\widetilde{O},\widetilde{O}\rangle$:
$(a_1,a_2,a_3)=(2,3,4)$, and $b,b_1,b_2$ and $b_3$ are
given by $b_1=1$, $m=12b+6+4b_2+3b_3$.
\item $\langle Z_{2m},Z_{2m};\widetilde{I},\widetilde{I}\rangle$:
$(a_1,a_2,a_3)=(2,3,5)$, and $b,b_1,b_2$ and $b_3$ are
given by $b_1=1$, $m=30b+15+10b_2+6b_3$.
\end{itemize}

Now let $\omega_0=\frac{\sqrt{-1}}{2}\sum_{i=1}^2 dz_i\wedge d\bar{z}_i$.
We consider the Hamiltonian $\s^1$-action on $(\C^2,\omega_0)$ given by
the complex multiplication, with the Hamiltonian function given by
$\mu(z_1,z_2)=\frac{1}{2}(|z_1|^2+|z_2|^2)$. It commutes with the action
of $G$ on $\C^2$, hence an Hamiltonian $\s^1$-action on the symplectic
orbifold $(\C^2,\omega_0)/G$ is resulted, with the Hamiltonian function
$\mu^\prime$ equaling $\frac{1}{2m}$ times the descendant of $\mu$ to
$\C^2/G$. Given any $r>0$, consider the subset $(\mu^\prime)^{-1}([0,r])
\subset\C^2/G$. According to \cite{Le}, we can collapse each fiber of
the $\s^1$-action in $(\mu^\prime)^{-1}(r)\subset (\mu^\prime)^{-1}([0,r])$
to obtain a closed symplectic $4$-dimensional orbifold, which we denote
by $X_r$. The symplectic $4$-orbifold $X_r$ contains
$(\mu^\prime)^{-1}([0,r))$
as an open symplectic suborbifold, and also contains a $2$-dimensional
symplectic suborbifold $C_0\equiv (\mu^\prime)^{-1}(r)/\s^1$. Note that
$(\mu^\prime)^{-1}(r)\rightarrow C_0$ is the canonical Seifert fibration
on $\s^3/G$ we mentioned earlier. Moreover, the Euler number of the
normal bundle of $C_0$ in $X_r$ equals the Euler number of the Seifert
fibration on $\s^3/G$.

Suppose we are given with a symplectic $s$-cobordism $W$ of the
elliptic $3$-manifold $\s^3/G$ to itself. Fix a sufficiently large
$r>0$, there are $0<r_1,r_2<r$ such that a neighborhood of the convex
end of $W$ is symplectomorphic to a neighborhood of the boundary
of $(\mu^\prime)^{-1}([0,r_1])$ and a neighborhood of the concave
end is symplectomorphic to a neighborhood of the boundary
of $X_r\setminus (\mu^\prime)^{-1}([0,r_2))$. We close up $W$ by
gluing $X_r\setminus (\mu^\prime)^{-1}([0,r_1))$ to the convex end
and gluing $(\mu^\prime)^{-1}([0,r_2])$ to the concave end. We
denote the resulting symplectic $4$-orbifold by $(X,\omega)$. Note
that $X$ inherits a $2$-dimensional symplectic suborbifold $C_0$
from $X_r\setminus (\mu^\prime)^{-1}([0,r_1))$, whose normal
bundle in $X$ has Euler number equaling that of the Seifert
fibration on $\s^3/G$. We fix an $\omega$-compatible almost
complex structure $J$ on $X$ such that $C_0$ is $J$-holomorphic
and $J$ is integrable in a neighborhood of each singular point of
$X$. (Note that the latter is possible because of the equivariant
Darboux' theorem.)

The $4$-orbifold $X$
has $4$ singular points. One of them, denoted by $p_0$, is inherited
from $(\mu^\prime)^{-1}([0,r_2])$ and has a neighborhod modeled by that
of $\{0\in\C^2\}/G$. The other three, denoted by $p_1,p_2$ and $p_3$, are
all contained in $C_0$, and are of type $(a_1,b_1)$, $(a_2,b_2)$ and
$(a_3,b_3)$ respectively, where $\{(a_i,b_i)\mid i=1,2,3\}$ is part of
the normalized Seifert invariant of the Seifert fibration on $\s^3/G$.
(Here a singular point is said of type $(a,b)$ if the isotropy
group is cyclic of order $a$ and the action on a local
uniformizing system is of weight $(1,b)$.) The Betti numbers of $X$
are $b_1(X)=b_3(X)=0$ and $b_2(X)=b_2^{+}(X)=1$. In fact, we have
$H_2(X;\Q)=\Q\cdot [C_0]$, and using the intersection product $C_0\cdot
C_0$, we may identify $H^2(X;\Q)$ with $H_2(X;\Q)=\Q\cdot [C_0]$
canonically. Finally, using the normalized Seifert invariant and
the adjunction formula (cf. \cite{C2}, Theorem 3.1), we obtain
$$
C_0\cdot C_0=\frac{4m^2}{|G|} \mbox{ and }
c_1(K_X)\cdot C_0=-\frac{4m(m+1)}{|G|},
$$
where $K_X$ is the canonical bundle of $(X,J)$, and $m$ is one half of 
the order of the subgroup $Z_{2m}$ of $G$ which preserves each fiber of 
the Hopf fibration.

With the preceding understood, we now introduce the relevant moduli
space of pseudoholomorphic curves. To this end, let $\Sigma$ be the
orbifold Riemann sphere of one orbifold point $z_\infty=\infty$ of
order $2m$. (Recall that $2m$ is the order of the cyclic subgroup
$Z_{2m}\subset G$ which preserves each fiber of the Hopf fibration.)
Note that $\Sigma$ has a unique complex structure. The group of automorphisms
of $\Sigma$, denoted by $\G$, is easily identified with the group
of linear translations on $\C$.

We shall consider the space $\M$ of $J$-holomorphic maps
$f:\Sigma\rightarrow X$ such that
\begin{itemize}
\item [{(1)}] The homology class $[f(\Sigma)]\in H_2(X;\Z)$ obeys
$[f(\Sigma)]\cdot C_0=1$.
\item [{(2)}] $f(z_\infty)=p_0$, and in a local representative
$(f_\infty,\rho_\infty)$ of $f$ at $z_\infty$,
$$
\rho_\infty(\mu_{2m})=\mu_{2m} I\equiv\left (\begin{array}{ll}
\mu_{2m} & 0\\
0 & \mu_{2m}\\
\end{array} \right )\in Z_{2m}, \mbox{ where }
\mu_k\equiv\exp(\sqrt{-1}\frac{2\pi}{k}).
$$
\end{itemize}

Here the notion of maps between orbifolds is as defined in
\cite{C1}. For the terminology used in this paper in connection
with $J$-holomorphic maps or curves, the reader is specially referred
to the earlier paper \cite{C2}. With these understood, we remark that
$\G$ acts on $\M$ by reparametrization.

The following proposition is the central technical result of this section.

\begin{prop}
The space $\M$ is nonempty, and is a smooth
$6$-dimensional manifold. Moreover, the quotient space
$\M/\G$ is compact.
\end{prop}

The proof of Proposition 3.1 will be given through a sequence of lemmas.
We begin with the Fredholm theory for pseudoholomorphic curves in
a symplectic $4$-orbifold $(X,\omega)$.

Let $(\Sigma,j)$ be an orbifold Riemann surface with a fixed complex
structure $j$. Consider the space of $C^k$ maps $[\Sigma;X]$ from
$\Sigma$ to $X$ for some fixed, sufficiently large integer $k>0$.
By Theorem 1.4 in Part I of \cite{C1}, the space $[\Sigma;X]$ is a smooth
Banach orbifold, which we may simply assume to be a smooth Banach
manifold for the sake of technical simplicity, because the relevant subset
$\M$ in the present case is actually contained in the
smooth part of $[\Sigma;X]$. There is a Banach bundle $\E\rightarrow
[\Sigma;X]$, with a Fredholm section $\underline{L}$ defined by
$$
\underline{L}(f)=df+J\circ df\circ j, \;\forall f\in [\Sigma;X].
$$
The zero locus $\underline{L}^{-1}(0)$ is the set of $J$-holomorphic
maps from $(\Sigma,j)$ into $X$. In the present case, $\M$
is contained in $\underline{L}^{-1}(0)$ as an open subset with respect
to the induced topology.

The index of the linearization $D\underline{L}$ at
$f\in\underline{L}^{-1}(0)$ can be computed using the index
theorem of Kawasaki \cite{Ka} for elliptic operators on orbifolds,
see Appendix A for a relevant review. To state the general index
formula for $D\underline{L}$ (cf. Lemma 3.2.4 of \cite{CR}), let
$(\Sigma,j)$ be an orbifold Riemann surface with orbifold points
$z_i$ of order $m_i$, where $i=1,2,\cdots, l$, and let
$f:\Sigma\rightarrow X$ be a $J$-holomorphic map from $(\Sigma,j)$
into an almost complex $4$-orbifold $(X,J)$. If a local representative
of $f$ at each $z_i$ is given by $(f_i,\rho_i)$ where $\rho_i(\mu_{m_i})$
acts on a local uniformizing system at $f(z_i)$ by $\rho_i(\mu_{m_i})
\cdot (z_1,z_2)=(\mu_{m_i}^{m_{i,1}}z_1,\mu_{m_i}^{m_{i,2}}z_2)$, with
$0\leq m_{i,1}, m_{i,2}<m_i$ (here $\mu_k\equiv\exp(\sqrt{-1}
\frac{2\pi}{k})$), then the index of $D\underline{L}$ at $f$ is $2d$ with
$$
d=c_1(TX)\cdot [f(\Sigma)]+2-2g_{|\Sigma|}-\sum_{i=1}^l
\frac{m_{i,1}+m_{i,2}}{m_i},
$$
where $g_{|\Sigma|}$ is the genus of the underlying Riemann surface
of $\Sigma$. In the present case, for each $f\in \M$,
one half of the index of $D\underline{L}$ at $f$ is
$$
d=\frac{m+1}{m}+2-\frac{1+1}{2m}=3.
$$
Thus $\M$ is a $6$-dimensional smooth manifold provided that it is
nomempty and $\underline{L}$ is transversal to the zero section at $\M$.

Transversality of the Fredholm section $\underline{L}$ at its
zero locus can be addressed in a similar fashion as in the case when
$X$ is a manifold. For the purpose of this paper, we shall use the
following regularity criterion, which is the orbifold analog of
Lemma 3.3.3 in \cite{McDS}.

\begin{lem}
Let $f:\Sigma\rightarrow X$ be a $J$-holomorphic map from an orbifold
Riemann surface into an almost complex $4$-orbifold. Suppose at each
$z\in\Sigma$, the map $f_z$ in a local representative $(f_z,\rho_z)$ of $f$
at $z$ is embedded. Then $f$ is a smooth point in the space of
$J$-holomorphic maps from $\Sigma$ into $X$ provided that $c_1(T\Sigma)
(\Sigma)>0$ and $c_1(TX)\cdot [f(\Sigma)]>0$.
\end{lem}

\pf
Let $E\rightarrow\Sigma$ be the pull back of $TX$ via $f$. Since $f_z$
is embedded for each $z\in \Sigma$, $T\Sigma$ is a subbundle of $E$,
and one has the decomposition $E=T\Sigma\oplus (E/T\Sigma)$. Then
a similar argument as in Lemma 3.3.3 of \cite{McDS} shows that
$D\underline{L}$ is surjective at $f$ if both
$(-c_1(T\Sigma)+c_1(K_\Sigma))(\Sigma)$ and
$(-c_1(E/T\Sigma)+c_1(K_\Sigma))(\Sigma)$ are negative. (Here $K_\Sigma$
is the canonical bundle of $\Sigma$.) The lemma follows easily.

\hfill $\Box$

Note that the conditions $c_1(T\Sigma)(\Sigma)>0$ and
$c_1(TX)\cdot [f(\Sigma)]>0$ in the previous lemma are met by
each $f\in\M$: $c_1(T\Sigma)(\Sigma)=1+\frac{1}{2m}>0$,
and $c_1(TX)\cdot [f(\Sigma)]=\frac{m+1}{m}>0$. Thus for the smoothness
of $\M$, it suffices to verify that for each $f\in\M$,
$f_z$ is embedded, $\forall z\in \Sigma$. This condition is verified
in the next lemma. But in order to state the lemma, it proves convenient
to introduce the following

\vspace{3mm}

\noindent{\bf Definition:}\hspace{3mm}
{\em
Let $C$ be a $J$-holomorphic curve in $X$ which contains the singular point
$p_0$, and is parametrized by $f:\Sigma\rightarrow X$. We call $C$ a
{\em quasi-suborbifold} if the following are met.
\begin{itemize}
\item $f$ induces a homeomorphism between the underlying Riemann surface
and $C$,
\item $f$ is embedded in the complement of the singular points in $X$,
\item a local representative $(f_z,\rho_z)$ of $f$ at each $z\in\Sigma$
where $f(z)$ is a singular point obeys {\em (i)} $f_z$ is embedded,
{\em (ii)} $\rho_z$ is isomorphic if $f(z)\neq p_0$, and if $f(z)=p_0$,
$\rho_z$ {\em (}which is injective by definition {\em)} maps onto the maximal
subgroup of $G$ that fixes the tangent space of $\mbox{Im }f_z$ at the
inverse image of $p_0$ in the local uniformizing system at $p_0$.
\end{itemize}
}

We remark that in terms of the adjunction formula (cf. Theorem 3.1 of
\cite{C2})
$$
g(C)=g_{\Sigma}+\sum_{\{[z,z^\prime]|z\neq z^\prime, f(z)=f(z^\prime)\}}
k_{[z,z^\prime]}+\sum_{z\in\Sigma} k_z,
$$
a $J$-holomorphic curve $C$ is a quasi-suborbifold if and only if
$k_{[z,z^\prime]}=0$ for all $[z,z^\prime]$, $k_z=0$ for any $z$
such that $f(z)\neq p_0$, and $k_{z_0}=\frac{1}{2m_0}(\frac{|G|}{m_0}-1)$
where $f(z_0)=p_0$. Here $m_0$ is the order of $z_0\in\Sigma$.
(Compare Corollary 3.3 in \cite{C2}, and note that
$\frac{1}{2m_0}(\frac{|G|}{m_0}-1)$ is the least of the possible values
of $k_{z_0}$.)

Now in the lemma below, we describe how the members of $\M$ look like.

\begin{lem}
Each $f\in\M$ is either a (multiplicity-one) parametrization
of a $J$-holomorphic quasi-suborbifold intersecting $C_0$ transversely
at a smooth point, or a multiply covered map onto a $J$-holomorphic
quasi-suborbifold intersecting $C_0$ at a singular point, such that
the order of the singular point equals the multiplicity of $f$.
Moreover, even in the latter case, the map $f_z$ in a local representative
$(f_z,\rho_z)$ of $f$ at $z$ is embedded for all $z\in\Sigma$.
\end{lem}

\pf
Set $C\equiv \mbox{Im }f$. We first consider the case where $f$ is not
multiply covered. Under this assumption, we have
$$
C\cdot C=\frac{|G|}{4m^2} \mbox{ and } c_1(K_X)\cdot C=-\frac{m+1}{m},
$$
which implies that the virtual genus
$$
g(C)=\frac{1}{2}(\frac{|G|}{4m^2}-\frac{m+1}{m})+1
=\frac{|G|}{8m^2}-\frac{m+1}{2m}+1.
$$
On the other hand, the orbifold genus $g_\Sigma=\frac{1}{2}(1-\frac{1}{2m})$
and $k_{z_\infty}\geq \frac{1}{4m}(\frac{|G|}{2m}-1)$. It follows easily
from the adjunction formula that $C$ is a quasi-suborbifold and
$C$ intersects $C_0$ transversely at a smooth point.

Now consider the case where $f$ is multiply covered. Let $s>1$ be the
multiplicity of $f$. Clearly $C$, $C_0$ are distinct, hence by the
intersection formula (cf. \cite{C2}, Theorem 3.2),
$$
\frac{1}{s}=C\cdot C_0=\sum_{i=1}^3 \frac{k_i}{a_i},
$$
where $a_i$ is the order of the singular point $p_i$, and $k_i\geq 0$
is an integer which is nonzero if and only if $p_i\in C\cap C_0$.
It follows immediately that $s\leq\frac{a_i}{k_i}$ if $k_i\neq 0$.
From the possible values of $(a_1,a_2,a_3)$, one can easily see that
$C$ intersects $C_0$ at exactly one singular point.

Suppose $C,C_0$ intersect at $p_i$ for some $i=1,2$ or $3$. Let
$\hat{f}:\hat{\Sigma}\rightarrow X$ be a (multiplicity-one)
parametrization of $C$ by a $J$-holomorphic map such that $f$ factors
through a map $\varphi:\Sigma\rightarrow\hat{\Sigma}$ to $\hat{f}$,
and let $\hat{f}(\hat{z}_0)=p_i$ for some $\hat{z}_0\in\hat{\Sigma}$
whose order is denoted by $\hat{m}_0$. Set $\hat{z}_\infty\equiv
\varphi(z_\infty)$. First, by the intersection formula, we get
$\frac{1}{s}\geq\frac{a_i/\hat{m}_0}{a_i}=\frac{1}{\hat{m}_0}$.
Hence $s\leq\hat{m}_0$. On the other hand, let $z_0\in\Sigma$ be
an inverse image of $\hat{z}_0$ under $\varphi$. Then $\hat{m}_0$ is
no greater than the degree of the branched covering $\varphi$ at $z_0$,
which is no greater than the total multiplicity $s$. This implies that
$s=\hat{m}_0$. Now we look at the point $\hat{z}_\infty$. Let $m_\infty$
be the degree of the branched covering $\varphi$ at $z_\infty$. Then (1)
the order of $\hat{z}_\infty$, denoted by $\hat{m}_\infty$, is no greater
than $2mm_\infty$, and (2) $m_\infty\leq s=\hat{m}_0$. In particular,
$\hat{m}_\infty\leq 2m\hat{m}_0$.

Now in the adjunction formula for $C$, the virtual genus
$$
g(C)=\frac{|G|}{8m^2\hat{m}_0^2}-\frac{m+1}{2m\hat{m}_0}+1,
$$
and on the right hand side,
$$
g_{\hat{\Sigma}}\geq\frac{1}{2}(1-\frac{1}{\hat{m}_0})+
\frac{1}{2}(1-\frac{1}{\hat{m}_\infty}), \;k_{\hat{z}_0}\geq
\frac{1}{2\hat{m}_0}(\frac{a_i}{\hat{m}_0}-1), 
k_{\hat{z}_\infty}
\geq \frac{1}{2\hat{m}_\infty}(\frac{|G|}{\hat{m}_\infty}-1).
$$
If $2m\hat{m}_0>\hat{m}_\infty$, then the adjunction
formula for $C$ gives rise to
$$
\frac{|G|}{4m\hat{m}_0}<1,
$$
which is impossible because
$$
1\leq\frac{|G|}{4ma_i}<\frac{|G|}{4m\hat{m}_0}.
$$
Hence $2m\hat{m}_0=\hat{m}_\infty$. With this in hand, the adjunction
formula further implies that $a_i=\hat{m}_0$ and $C$
is a quasi-suborbifold, and the multiplicity of $f$ equals the order
of the singular point where $C,C_0$ intersect.

It remains to check that $f_z$ is embedded for all $z\in\Sigma$.
But this follows readily from (1) $C$ is a quasi-suborbifold, and
(2) $\varphi:\Sigma\rightarrow\hat{\Sigma}$ is a cyclic branched
covering of degree $s$, branched at $z_0,z_\infty$, and $\hat{m}_0=s$,
$\hat{m}_\infty=2ms$.

\hfill $\Box$

Up to this point, we see that $\M$ is a $6$-dimensional
smooth manifold provided that it is nonempty. Next we show

\begin{lem}
The quotient space $\M/\G$ is compact.
\end{lem}

\pf
According to the orbifold version of the Gromov compactness theorem
proved in \cite{CR}, for any sequence of maps $f_n\in\M$,
there exists a subsequence which converges to
a cusp-curve after suitable reparametrization. More concretely,
after reparametrization if necessary, there is a subsequence of $f_n$,
which is still denoted by $f_n$ for simplicity, and there are at most
finitely many simple closed loops $\gamma_1,\cdots,\gamma_l\subset\Sigma$
containing no orbifold points, and a nodal orbifold Riemann surface
$\Sigma^\prime=\cup_\omega\Sigma_\omega$ obtained by collapsing
$\gamma_1,\cdots,\gamma_l$, and a $J$-holomorphic map
$f:\Sigma^\prime\rightarrow X$, such that (1) $f_n$ converges
in $C^\infty$ to $f$ on any given compact subset in the complement of
$\gamma_1,\cdots,\gamma_l$, (2) $[f(\Sigma^\prime)]=[f_n(\Sigma)]
\in H_2(X;\Q)$, (3) if $z_\omega\in\Sigma_\omega, z_\nu\in
\Sigma_\nu$ are two distinct points (here $\Sigma_\nu=\Sigma_\omega$
is allowed) with orders $m_\omega, m_\nu$ respectively,
such that $z_\omega,z_\nu$ are the image of the same simple closed loop
collapsed under $\Sigma\rightarrow\Sigma^\prime$, then $m_\omega=m_\nu$,
and there exist local representatives
$(f_\omega,\rho_\omega), (f_\nu,\rho_\nu)$
of $f$ at $z_\omega,z_\nu$, which obey $\rho_{\omega}(\mu_{m_\omega})=
\rho_\nu(\mu_{m_\nu})^{-1}$, and (4) if $f$ is constant over a component
$\Sigma_\nu$ of $\Sigma^\prime$, then either the underlying surface of
$\Sigma_\nu$ has nonzero genus, or $\Sigma_\nu$ contains at least $3$
special points, where a special point is either an orbifold point inherited
from $\Sigma$ or any point resulted from collapsing a simple closed loop
in $\Sigma$. Regarding the last point about constant components, since
in the present case $\Sigma$ is an orbifold Riemann sphere with only one
orbifold point $z_\infty$, any constant component in the limiting cusp-curve
must be obtained by collapsing at least $2$ simple closed loops, and if
$z_\infty$ is not contained, by collapsing at least $3$ simple closed loops.

With the preceding understood, note that there are two possibilities:
(1) none of the simple closed loops $\{\gamma_i\}$ is null-homotopic
in the complement of $z_0,z_{\infty}$ in $\Sigma$ where $f_n(z_0)\in C_0$
and $f_n(z_\infty)=p_0$, or (2) there is a simple closed loop
$\gamma\in\{\gamma_i\}$ such that $\gamma$ bounds a disc $D$ in the
complement of $z_0,z_{\infty}$ in $\Sigma$, such that $D$ contains none
of the simple closed loops $\gamma_1,\cdots,\gamma_l$.

Case (1): under this assumption, it is easily seen that there are no
constant components in the limiting cusp-curve. Moreover, there is a component
$\Sigma_\omega$ such that $f_\omega\equiv f|_{\Sigma_\omega}:
\Sigma_\omega\rightarrow X$ obeys the following conditions:
\begin{itemize}
\item There exists an orbifold point $w_\infty\in\Sigma_\omega$ of order
$2m$ inherited from $\Sigma$ (i.e. $w_\infty=z_\infty$) such that
$f_\omega(w_\infty)=p_0$, and a local representative of $f_\omega$ at
$w_\infty$ obeys the second condition in the definition of $\M$.
\item $f_\omega^{-1}(C_0)$ consists of only one point $w_0$, which is
necessarily obtained from collapsing one of the simple closed loops
$\gamma_1,\cdots,\gamma_l$.
\item All points in $\Sigma_\omega\setminus\{w_0,w_\infty\}$ are
regular, i.e. of order $1$ in $\Sigma_\omega$.
\end{itemize}
First of all, we show that $w_0$ is actually a regular point of
$\Sigma_\omega$. In order to see this, we only need to consider
the case where $f_\omega(w_0)$ is a singular point, say $p_i$, for
some $i=1,2$ or $3$. (Note that $w_0$ is automatically a regular
point if $f_\omega(w_0)$ is a smooth point.) Let $(w_1,w_2)$
be holomorphic coordinates on a local uniformizing system at $p_i$,
where $C_0$ is locally given by $w_2=0$, and the singular fiber of
the Seifert fibration at $p_i$ is defined by $w_1=0$, $|w_2|\equiv
\mbox{constant}$. The local $\Z_{a_i}$-action is given by $\mu_{a_i}
\cdot (w_1,w_2)=(\mu_{a_i}w_1,\mu_{a_i}^{b_i}w_2)$. (Here $(a_i,b_i)$
is the normalized Seifert invariant at $p_i$.) Let $m_0\geq 1$
be the order of $w_0$, and let $(f_0,\rho_0)$ be a local representative
of $f_\omega$ at $w_0$, where $\rho_0(\mu_{m_0})=\mu_{m_0}^r$ with
$0\leq r<m_0$, $r,m_0$ relatively prime, and $f_0(w)=(c(w^{l_1}+\cdots),
w^{l_2}+\cdots)$ (note that $\mbox{Im }f_\omega\neq C_0$). By a
$\Z_{m_0}$-equivariant change of coordinates $w^\prime\equiv
w(1+\cdots)^{1/l_2}$ near $w=0$, we may simply assume $f_0(w)
=(c(w^{l_1}+\cdots),w^{l_2})$. Furthermore, $l_2\equiv b_i r\pmod{m_0}$,
so that $l_2,m_0$ are relatively prime. With these understood, note that
the image of the link of $w_0$ in $\Sigma_\omega$ under $f_\omega$
is parametrized in the local uniformizing system by $f_0(\epsilon\exp
(\sqrt{-1}\frac{2\pi}{m_0}\theta))$, $0\leq\theta\leq 1$. Through
$f_t(w)\equiv (c(1-t)(w^{l_1}+\cdots),w^{l_2})$, $0\leq t\leq 1$,
it is homotopic to $(0,\epsilon^{l_2}\exp(\sqrt{-1}\frac{2\pi l_2}{m_0}
\theta))$ in the complement of $C_0$. It follows easily that the link
of $w_0$ in $\Sigma_\omega$ under $f_\omega$ is homotopic in $\s^3/G$ to
$\frac{l_2a_i}{m_0}$ times of the singular fiber of the Seifert fibration
at $p_i$, whose homotopy class in $\pi_1(\s^3/G)=G$ has order $2ma_i$.
On the other hand, the link of $w_0$ in $\Sigma_\omega$ under $f_\omega$
is homotopic in $W$ to the image of the inverse of the link of $w_\infty$
in $\Sigma_\omega$ under $f_\omega$. The latter's homotopy class in
$G$ is $\mu_{2m}^{-1} I\in Z_{2m}$, which implies that the former's
homotopy class is an element of order $2m$ in $Z_{2m}$ (in fact,
it is $\mu_{2m}^{-1} I\in Z_{2m}$, cf. Lemma 3.5 below). This gives
$\frac{l_2a_i}{m_0}\cdot 2m=2m a_i l$ for some $l>0$, which contradicts
the fact that $l_2$, $m_0$ are relatiely prime unless $m_0=1$.
Therefore $w_0$ is a regular point of $\Sigma_\omega$.

Now note that $f_\omega:\Sigma_\omega\rightarrow X$ satisfies
all the conditions in the definition of $\M$ except for
the first one. i.e. $[f_\omega(\Sigma_\omega)]\cdot C_0=1$, which
we prove next. To see this, let $\hat{f}_\omega:\hat{\Sigma}_\omega
\rightarrow X$ be the multiplicity-one parametrization of $C_\omega
\equiv \mbox{Im }f_\omega$ obtained by factoring $f_\omega$ through
a branched covering map
$\varphi:\Sigma_\omega\rightarrow\hat{\Sigma}_\omega$ of degree $s$.
(If $f_\omega$ is not multiply covered, we simply
let $\hat{\Sigma}_\omega\equiv\Sigma_\omega$, $\hat{f}_\omega\equiv
f_\omega$, and $s=1$.) Set $\hat{w}_0\equiv\varphi(w_0)\in
\hat{\Sigma}_\omega$, and let $\hat{m}_0\geq 1$ be the order of $\hat{w}_0$
in $\hat{\Sigma}_\omega$. Then by the intersection formula, we get
$$
\frac{1}{s}\geq\frac{1}{s}\cdot [f_\omega(\Sigma_\omega)]\cdot C_0=
C_\omega\cdot C_0\geq \frac{a_i/\hat{m}_0}{a_i}=\frac{1}{\hat{m}_0}.
$$
Hence $s\leq\hat{m}_0$. On the other hand, $\hat{m}_0$ is
no greater than the degree of the branched covering $\varphi$ at $w_0$,
which is no greater than the total multiplicity $s$. This implies that
$s=\hat{m}_0$ and $C_\omega\cdot C_0=\frac{1}{s}$. Hence
$[f_\omega(\Sigma_\omega)]\cdot C_0=s\cdot C_\omega\cdot C_0=1$.
It is clear that $f=f_\omega\in\M$ in this case.

Case (2): Let $\Sigma_\omega$ be the component obtained from the
disc $D$ that $\gamma$ bounds, and let $z_0\in\Sigma_\omega$ be
the point which is the image of $\gamma$ under
$D\rightarrow\Sigma_\omega$. Note that $f$ is nonconstant over
$\Sigma_\omega$. Set $f_\omega\equiv f|_{\Sigma_\omega}$ and
$C_\omega\equiv\mbox{Im }f_\omega$. Since $f_n(D)$ is disjoint from $C_0$,
either $f_\omega^{-1}(C_0)$ consists of only one point $z_0$, or
$C_\omega=C_0$. However, the latter case can be ruled out for the
following reason. Note that $\Sigma_\omega$ contains at most one
orbifold point, hence is simply connected as an orbifold. Consequently,
the degree of the map $f_\omega:\Sigma_\omega\rightarrow C_\omega=C_0$
is at least $\frac{|G|}{2m}$, which is the order of $G/Z_{2m}$, the
orbifold fundamental group of $C_0$. It follows that
$[f_\omega(\Sigma_\omega)]\cdot C_0\geq\frac{|G|}{2m}\cdot
\frac{4m^2}{|G|}=2m>1$, which is a contradiction.

Let $\hat{f}_\omega:\hat{\Sigma}_\omega\rightarrow X$ be the
multiplicity-one parametrization of $C_\omega$ obtained by factoring
$f_\omega$ through a map $\varphi:\Sigma_\omega\rightarrow
\hat{\Sigma}_\omega$ of degree $s$.
(If $f_\omega$ is not multiply covered, we simply
let $\hat{\Sigma}_\omega\equiv\Sigma_\omega$, $\hat{f}_\omega\equiv
f_\omega$, and $s=1$.) Set $\hat{z}_0\equiv\varphi(z_0)\in
\hat{\Sigma}_\omega$, and let $m_0$ be the order of $\hat{z}_0$ in
$\hat{\Sigma}_\omega$.

Note that in this case $\Sigma_\omega$ is necessarily not the only
component of $\Sigma^\prime$ over which $f$ is nonconstant.
Consequently $C_\omega\cdot C_0<1$, and $f_\omega(z_0)$ is a
singular point of $X$ on $C_0$, say $p_i$ for some $i=1,2$ or $3$.
Let $z_1,z_2$ be the holomorphic coordinates on a local
uniformizing system at $p_i$, with local group action given by
$\mu_{a_i}\cdot (z_1,z_2)=(\mu_{a_i}z_1,\mu_{a_i}^{b_i}z_2)$, such
that $C_0$ is locally given by $z_2=0$ and the singular fiber of
the Seifert fibration on $\s^3/G$ at $p_i$ is given by $z_1=0$,
$|z_2|\equiv\mbox{constant}$. (Here $(a_i,b_i)$ is the normalized Seifert
invariant at $p_i$.) Let $(f_0,\rho_0)$ be a local representative of
$\hat{f}_\omega$ at $\hat{z}_0$, where we write $f_0(z)=(c(z^{l_1}+\cdots),
z^{l_2})$ (note that $C_\omega\neq C_0$).
Then by the intersection formula, we have
$$
\frac{1}{s}\geq C_\omega\cdot C_0=\frac{(a_i/m_0)\cdot
l_2}{a_i}=\frac{l_2}{m_0}.
$$
On the other hand, by a similar argument, we see that the link of $p_i$
in $C_\omega$ is homotopic to $\frac{l_2a_i}{m_0}$ times of the singular
fiber of the Seifert fibration at $p_i$, whose homotopy class in $G$
is of order $2ma_i$. Since $\s^3/G\rightarrow W$ is a homotopy equivalence,
and $f_n(\gamma)$, which bounds a disc $f_n(D)\subset W$, is homotopic
to $s$ times of the link of $p_i$ in $C_\omega$, we have $s\cdot
\frac{l_2a_i}{m_0}=2ma_il$ for some $l>0$. But this contradicts
the inequality $\frac{1}{s}\geq \frac{l_2}{m_0}$ we obtained earlier.

Hence case (2) is impossible, and therefore the quotient space
$\M/\G$ is compact.

\hfill $\Box$

It remains to show, in the proof of Proposition 3.1, that $\M$
is not empty. This will be achieved in the following three steps:
\begin{itemize}
\item [{(1)}] Construct an orbifold complex line bundle $E\rightarrow X$
such that $c_1(E)\cdot C_0=1$.
\item [{(2)}] Show that the associated Seiberg-Witten invariant is
nonzero in Taubes chamber.
\item [{(3)}] Apply Theorem 2.2 (2) to produce a $J$-holomorphic curve
$C$ such that $C=\mbox{Im }f$ for some $f\in\M$.
\end{itemize}

For step (1), we derive a preliminary lemma first. To state the lemma,
let $h:\s^3/G\rightarrow \s^3/G$ be the simple homotopy equivalence induced
by the $s$-cobordism $W$. Then there is a pair $(\hat{h},\hat{\rho}):
(\s^3,G)\rightarrow (\s^3,G)$ where $\hat{h}:\s^3\rightarrow\s^3$ is
$\hat{\rho}$-equivariant and descends to $h:\s^3/G\rightarrow \s^3/G$.
The pair $(\hat{h},\hat{\rho})$ is unique up to conjugation by an element
of $G$.

\begin{lem}
The restriction of $\hat{\rho}$ to $Z_{2m}\subset G$ is the identity map.
\end{lem}

\pf
Recall the double cover $\phi:\s^3\times\s^3\rightarrow SO(4)$,
which is defined by sending $(q_1,q_2)\in \s^3\times\s^3$ to the
matrix in $SO(4)$ that sends $x\in\R^4=\h$ to $q_1xq_2^{-1}\in\h=\R^4$.
Regard $G$ as a subgroup of $\phi(\s^1\times\s^3)$.

Note that as a simple homotopy equivalence, $h:\s^3/G\rightarrow \s^3/G$
is homotopic to a diffeomorphism (cf. \cite{Tu}, and for a proof,
\cite{KS3}). On the other hand, any diffeomorphism between elliptic
$3$-manifolds is homotopic to an isometry, cf. e.g. \cite{McC}, hence
$h$ is homotopic to an isometry. It follows easily that
$\hat{h}:\s^3\rightarrow\s^3$ is $\hat{\rho}$-equivariantly homotopic
to an isometry $\xi\in SO(4)$. In particular, $\hat{\rho}(g)=\xi g\xi^{-1}$.

Now let $\xi=\phi(q,q^\prime)$ and $g=\phi(x,y)$. Then $\hat{\rho}(g)=
\phi(qxq^{-1},q^\prime y (q^\prime)^{-1})$. Note that for any $g\in
Z_{2m}\subset G$, $g=\phi(x,1)$ with $x=(\mu_{2m}^l,0)$,
$0\leq l\leq 2m-1$. If we let $q=(w_1,w_2)$, then
$$
qxq^{-1}=(|w_1|^2\mu_{2m}^l+|w_2|^2\mu_{2m}^{-l},
w_1w_2(\mu_{2m}^{-l}-\mu_{2m}^l)).
$$
Note that when $m=1$, $qxq^{-1}=x$ so that the lemma holds trivially.
For the case where $m\neq 1$, the fact that $qxq^{-1}\in\s^1$ implies
that either $w_1$ or $w_2$ must be zero. Clearly, for any $g\in Z_{2m}$,
$\hat{\rho}(g)=g$ iff $w_2=0$ and $\hat{\rho}(g)=g^{-1}$ iff $w_1=0$.

It remains to show that $\hat{\rho}(g)=g^{-1}$, $\forall g\in Z_{2m}$,
is impossible. Here we need to use the assumption that $W$ is
symplectic. Let $\widetilde{W}$ be the universal cover of $W$.
Note that the canonical bundle $K_{\widetilde{W}}$ is trivial.
This gives rise to a representation $\theta:G=\pi_1(W)\rightarrow
\s^1$, which obeys $\theta=\theta\circ\hat{\rho}$. Let $g\in Z_{2m}$
be the matrix $\mu_{2m}I$. Then $\theta(g)=\mu_{2m}^2$, which
implies $\mu_m=\mu_m^{-1}$ if $\hat{\rho}(g)=g^{-1}$. But this is
impossible unless $m=2$, which occurs only when $G=\langle
Z_{4m},Z_{2m};\widetilde{D}_n,C_{2n}\rangle$. But even in this
case, $\hat{\rho}(g)\neq g^{-1}$ because otherwise, we would have
$q=(0,w_2)$, which implies that $\hat{\rho}(\mu_{4m}\alpha)=
\mu_{4m}^{-1}q^\prime\alpha(q^\prime)^{-1}$ for any
$\alpha\in\widetilde{D}_n$ whose class is nonzero in
$\widetilde{D}_n/C_{2n}$. But $\theta(\mu_{4m}\alpha)=\mu_{2m}$
and $\theta(\mu_{4m}^{-1}q^\prime\alpha (q^\prime)^{-1})=\mu_{2m}^{-1}$,
which contradicts $\theta=\theta\circ\hat{\rho}$ and $m=2$.
Hence the lemma.

\hfill $\Box$

Now back to step (1) of the proof. In the following lemma, we give an 
explicit construction of the orbifold complex line bundle $E$. 

\begin{lem}
There exists a canonically defined orbifold complex line bundle
$E\rightarrow X$ such that $c_1(E)\cdot C_0=1$.
\end{lem}

\pf
Note that $X$ is decomposed as $N\bigcup W\bigcup N_0$, where $N$
is a regular neighborhood of $C_0$, which is diffeomorphic to the
unit disc bundle associated to the Seifert fibration on $\s^3/G$,
and $N_0=\B^4/G$ is a regular neighborhood of the singular point
$p_0$. The orbifold complex line bundle $E$ will be defined by
patching together an orbifold complex line bundle on each of
$N, W$ and $N_0$, which agree on the intersections.

The bundle on $N$ is defined as follows. Take the complex line
bundle on the complement of the singular points $p_1,p_2$ and $p_3$
in $N$, which is Poincar\'{e} dual to a regular fiber of $N$. (The
regular fibers of $N$ are so oriented that the intersection with
$C_0$ has a $+$ sign.) This bundle is trivial on the link of each $p_i$
in $N$, so we can simply extend it over to the whole $N$ trivially to
obtain the orbifold complex line bundle on $N$.

The restriction of the bundle on $N$ to $\partial N=\s^3/G$
is Poincar\'{e} dual to a regular fiber of the Seifert fibration.
By Lemma 3.5, there exists a map $\psi:\s^1\times [0,1]\rightarrow W$
such that $\psi(\s^1\times\{0\})$ is a regular fiber of the Seifert
fibration on $\partial N=\s^3/G$, and $\psi(\s^1\times\{1\})$ is
the image of the boundary of a generic unit complex linear disc in
$\B^4$ under the quotient map $\partial \B^4=\s^3\rightarrow\s^3/G$.
We let the bundle on $W$ be the Poincar\'{e} dual of
$\psi(\s^1\times [0,1])$.

It remains to construct an orbifold complex line bundle $E_0$ on
$N_0=\B^4/G$ such that the restriction of $E_0$ on $\partial\B^4/G$ is
Poincar\'{e} dual to $\psi(\s^1\times\{1\})$. The resulting
orbifold complex line bundle $E\rightarrow X$ clearly obeys
$c_1(E)\cdot C_0=1$.

To this end, note that given any representation $\rho:G\rightarrow\s^1$,
there exists an orbifold complex line bundle on $N_0$, which is given
by the projection $(\B^4\times\C,G)\rightarrow (\B^4,G)$ on the 
uniformizing system,
where the action of $G$ on $\B^4\times\C$ is given by
$g\cdot (z,w)=(gz,\rho(g)w)$, $\forall (z,w)\in\B^4\times\C, g\in G$.
With this understood, the definition of $E_0\rightarrow N_0$ for the
various cases of $G$ is given below.
\begin{itemize}
\item $\langle
Z_{2m},Z_{2m};\widetilde{D}_n,\widetilde{D}_n\rangle$:
$\rho(h)=\mu_{2m}^{2n}$, $\rho(x)=(-1)^n$, and $\rho(y)=1$, where
$h=\mu_{2m}I\in Z_{2m}$, and $x,y$ are the generators of
$\widetilde{D}_n$ with relations $x^2=y^n=(xy)^2=-1$.
\item $\langle Z_{4m},Z_{2m};\widetilde{D}_n,C_{2n}\rangle$:
$\rho(h^2)=\mu_{2m}^{2n}$, $\rho(hx)=(-\mu_{2m})^n$, 
$\rho(y)=1$, where $h=\mu_{4m}I\in Z_{4m}$, and $x,y$ are the
generators of $\widetilde{D}_n$ with relations $x^2=y^n=(xy)^2=-1$.
\item $\langle Z_{2m},Z_{2m};\widetilde{T},\widetilde{T}\rangle$:
$\rho(h)=\mu_{2m}^{12}$, $\rho(x)=1$, and $\rho(y)=1$, where
$h=\mu_{2m} I\in Z_{2m}$, and $x,y$ are the generators of
$\widetilde{T}$ with relations $x^2=y^3=(xy)^3=-1$.
\item $\langle Z_{6m},Z_{2m};\widetilde{T},\widetilde{D}_2\rangle$:
$\rho(h^3)=\mu_{2m}^{12}$, $\rho(x)=1$, and $\rho(hy)=\mu_{2m}^4$,
where $h=\mu_{6m} I\in Z_{6m}$, and $x,y$ are the generators of
$\widetilde{T}$ with relations $x^2=y^3=(xy)^3=-1$.
\item $\langle Z_{2m},Z_{2m};\widetilde{O},\widetilde{O}\rangle$:
$\rho(h)=\mu_{2m}^{24}$, $\rho(x)=1$, and $\rho(y)=1$, where
$h=\mu_{2m} I\in Z_{2m}$, and $x,y$ are the generators of
$\widetilde{O}$ with relations $x^2=y^4=(xy)^3=-1$.
\item $\langle Z_{2m},Z_{2m};\widetilde{I},\widetilde{I}\rangle$:
$\rho(h)=\mu_{2m}^{60}$, $\rho(x)=1$, and $\rho(y)=1$, where
$h=\mu_{2m} I\in Z_{2m}$, and $x,y$ are the generators of
$\widetilde{I}$ with relations $x^2=y^5=(xy)^3=-1$.
\end{itemize}

The verification that the restriction of $E_0\rightarrow N_0$
to $\partial N_0$ is Poincar\'{e} dual to $\psi(\s^1\times\{1\})$
goes as follows. Fix a generic vector $u=(u_1,u_2)\in\C^2$, we let
$f_u$ be the linear function on $\C^2$ defined by
$$
f_u(z_1,z_2)\equiv u_1z_1+u_2z_2.
$$
The action of $g\in G$ as a $2\times 2$ complex valued matrix on
$f_u$ is given by $g^\ast f_u=f_{ug}$, where $ug=(u_1,u_2)g$ is
the row vector obtained from multiplying by $g$ on the right.
With this understood, consider the epimorphism $\pi:G\rightarrow\Gamma
\equiv G/Z_{2m}$, where $\Gamma$ is isomorphic to the corresponding
subgroup (dihedral, tetrahedral, octahedral, or icosahedral) in $SO(3)$.
For any $\gamma\in\Gamma$, we fix a $\hat{\gamma}\in G$ such that
$\pi(\hat{\gamma})=\gamma$. Then consider the product
$$
f(z)\equiv\prod_{\gamma\in\Gamma} f_{u\hat{\gamma}}(z), \;\forall z\in\C^2.
$$
The claim is that for any $g\in G, z\in\C^2$, $f(gz)=\rho(g)f(z)$, so that
$z\mapsto (z,f(z))$ is an equivariant section of the $G$-bundle
$\B^4\times\C\rightarrow\B^4$, which descends to a section $s$ of
the orbifold complex line bundle $E_0\rightarrow N_0$. The zero
locus of $s$ in $\partial N_0$ is the image of $f^{-1}_u(0)\cap\s^3$
under $\s^3\rightarrow\s^3/G=\partial N_0$, which can be so arranged
that it is actually $\psi(\s^1\times\{1\})$.

So it remains to verify the claim that for any $g\in G, z\in\C^2$,
$f(gz)=\rho(g)f(z)$. This is elementary but tedious, so we shall only
illustrate it by a simple example but also with some general remarks. The
details for all other cases are left out to the reader.

Consider the case $G=\langle Z_{2m},Z_{2m};\widetilde{D}_3,\widetilde{D}_3
\rangle$. The dihedral group $D_3$ is generated by $\alpha,\beta$
with relations $\alpha^2=\beta^3=(\alpha\beta)^2=1$, while the
binary dihedral group $\widetilde{D}_3$ is generated by $x,y$ with
relations $x^2=y^3=(xy)^2=-1$. Clearly $x\mapsto \alpha$, $y\mapsto\beta$
under $\widetilde{D}_3\rightarrow D_3$. Set $h\equiv\mu_{2m} I\in Z_{2m}$.
In this case, we may take
$$
f\equiv f_{u}f_{uy}f_{uy^2}f_{ux}f_{uyx}f_{uy^2x}.
$$
One can easily check that $f(hz)=\mu_{2m}^6 f(z)$,
$$
f(xz)=f_{ux}(z)f_{uyx}(z)f_{uy^2x}(z)f_{ux^2}(z)f_{uyx^2}(z)f_{uy^2x^2}(z)
=(-1)^3 f(z),
$$
and similarly $f(yz)=(-1)^2 f(z)=f(z)$.

As for the general remarks, the dihedral case can be similarly
handled as in the above example. For the tetrahedral case, the
order of the group $\Gamma=G/Z_{2m}$ is $12$, so it is not terribly
complicated. For the octahedral case, the trick is to fix an explicit
identification between the octahedral group $O$ and the symmetric group
$S_4$, e.g. $\alpha\mapsto (12), \beta\mapsto (1234)$ where
$\alpha,\beta$ are generators of $O$ with relations $\alpha^2=\beta^4=
(\alpha\beta)^3=1$, and use the identification between $O$ and $S_4$
to guide the manipulation of the words generated by $\hat{\alpha}$
and $\hat{\beta}$, where $\hat{\alpha},\hat{\beta}\in\widetilde{O}$
are some fixed choice of elements which obey
$\hat{\alpha}\mapsto\alpha$, $\hat{\beta}\mapsto
\beta$ under $\widetilde{O}\rightarrow O$. The case of icosahedral group
is actually quite simple. The observation is that
$H_1(\s^3/\widetilde{I};\Z)$ is trivial, so that any representation
${\rho}^\prime: G\rightarrow\s^1$ obtained from $f(gz)={\rho}^\prime(g)f(z)$
has to satisfy ${\rho}^\prime(g)=1$, $\forall g\in\widetilde{I}$, because
${\rho}^\prime$ factors through $H_1(\s^3/G;\Z)$.

\hfill $\Box$

Next for step (2), we show that the Seiberg-Witten invariant corresponding 
to $E$ is nonzero in Taubes chamber. First of all, we observe the following
lemma. 

\begin{lem}
The Seiberg-Witten invariant corresponding to $E$ is zero in the
$0$-chamber.
\end{lem}

\pf
Decompose $X$ as $X_1\bigcup X_2$ where $X_1$ is a regular
neighborhood of $C_0$. Note that $X_1$ is diffeomorphic to the unit disc
bundle associated to the Seifert fibration on $\s^3/G$.

The lemma follows readily from the fact that $X_1$ has a
Riemannian metric of positive scalar curvature which is
a product metric near $\partial X_1$. Accept this fact
momentarily, and suppose that the Seiberg-Witten invariant is nonzero
in the $0$-chamber. Then one can stretch the neck along
$\partial X_1=\partial X_2$, such that any solution of
the Seiberg-Witten equations on $X$ will yield a solution
$(A,\psi)$ on $\hat{X}_1\equiv X_1\cup [0,-\infty)\times\partial X_1$,
where $|\psi|$ converges to zero exponentially fast along the cylindrical
end of $\hat{X}_1$. Since the natural metric on $\hat{X}_1$ is
of positive scalar curvature, we must have $\psi\equiv 0$ by
the Weitzenb\"{o}ck formula. But this implies that
$P_{+}F_A=\frac{1}{4}\tau(\psi\otimes\psi^\ast)\equiv 0$, which
contradicts the fact that $\frac{\sqrt{-1}}{2\pi}
\int_{C_0}F_A=c_1(E^2\otimes K_X^{-1})\cdot C_0\neq 0$.

As for the fact that $X_1$ has a metric of positive scalar
curvature, here is a proof. Note that $X_1=((\s^3/G)\times D^2)/\s^1$,
where the $\s^1$-action on $\s^3/G$ defines the Seifert fibration, and
where $D^2$ is the unit $2$-disc with the $\s^1$-action given by
complex multiplication. Give $(\s^3/G)\times D^2$ a product
metric such that on the factor $\s^3/G$, it is the metric of
constant curvature which is clearly invariant under the
$\s^1$-action, and on the factor $D^2$, it is an $\s^1$-invariant
metric with nonnegative curvature which is a product metric near the
boundary. Now observe that the orthogonal complement of the vector
field generated by the $\s^1$-action on $(\s^3/G)\times D^2$ is an
$\s^1$-equivariant subbundle of the tangent bundle of $(\s^3/G)\times D^2$,
which canonically defines a Riemannian metric on $X_1$ through the projection
$(\s^3/G)\times D^2\rightarrow X_1$, making it into a Riemannian
submersion in the sense of O'Neill \cite{O'N}. It follows easily from
the calculation therein that the metric on $X_1$ has positive scalar
curvature.

\hfill $\Box$

Observe that $c_1(K_X)\cdot C_0<0$, so that $c_1(S_{+}^E)\cdot [\omega]
=c_1(K_X^{-1}\times E^2)\cdot [\omega]>0$. By the wall-crossing
formula in Lemma 2.1, the Seiberg-Witten invariant corresponding to $E$
is nonzero in Taubes chamber provided that the dimension of the
corresponding moduli space of the Seiberg-Witten equations is
nonnegative, which is shown in the next lemma.

\begin{lem}
The dimension of the Seiberg-Witten moduli space corresponding to
$E$, denoted by $d(E)$, is given for the various cases of $G$
in the following list.
\begin{itemize}
\item $\langle Z_{2m},Z_{2m};\widetilde{D}_n,\widetilde{D}_n\rangle$
or $\langle Z_{4m},Z_{2m};\widetilde{D}_n,C_{2n}\rangle$:
$d(E)=\delta+2+\frac{1}{2}((-1)^\delta-1)$
if $m<n$, where $n=\delta m+r$ with $0\leq r\leq m-1$,
and $d(E)=2$ if $m>n$.
\item $\langle Z_{2m},Z_{2m};\widetilde{T},\widetilde{T}\rangle$
or $\langle Z_{6m},Z_{2m};\widetilde{T},\widetilde{D}_2\rangle$:
$d(E)=2$ if $m\neq 1$, and $d(E)=8$ if $m=1$.
\item $\langle Z_{2m},Z_{2m};\widetilde{O},\widetilde{O}\rangle$:
$d(E)=2$ if $m\neq 1$, and $d(E)=14$ if $m=1$.
\item $\langle Z_{2m},Z_{2m};\widetilde{I},\widetilde{I}\rangle$:
$d(E)=2$ if $m\neq 1,7$, $d(E)=4$ if $m=7$, and $d(E)=32$ if
$m=1$.
\end{itemize}
\end{lem}

The proof of Lemma 3.8 is given in Appendix A.

\vspace{2mm}

Now the last step, where we apply Theorem 2.2 (2) to produce a 
$J$-holomorphic curve $C$ such that $C=\mbox{Im }f$ for some
$f\in\M$.
.
\begin{lem}
The space $\M$ is nonempty.
\end{lem}

The proof of Lemma 3.9 is given at the end of this section.
Accepting it for now, and hence Proposition 3.1, we shall prove
Theorem 1.2 next.

\vspace{3mm}

\noindent{\bf Proof of Theorem 1.2}

\vspace{2mm}

The group $\G$ acts on $\M$ smoothly (see the
general discussion at the end of \S 3.3, Part I of \cite{C1}).
Moreover, the action is free at any $f\in\M$ which is not
multiply covered. At a multiply covered $f\in\M$ with
multiplicity $a>1$, the isotropy subgroup is
$\{(\mu_a^l,0)\mid l=0,1,\cdots,a-1\}\subset\G$ up to
conjugations in $\G$. Here $\G$
is canonically identified with the group of linear translations on $\C$,
$\{(s,t)\mid s\in\C^\ast,t\in\C\}$. Clearly, $\M\rightarrow \M^\dag
\equiv \M/\G$ is a smooth orbifold principle
$\G$-bundle over a compact $2$-dimensional orbifold.

Recall that the domain of each $f\in\M$ is the orbifold
Riemann sphere $\Sigma$ of one orbifold point $z_\infty\equiv \infty$
of order $2m$. We identify $\Sigma\setminus\{z_\infty\}$ canonically
with $\C$ such that the action of $\G$ on $\Sigma\setminus\{z_\infty\}$ 
is given by linear translations on $\C$. We introduce the associated 
orbifold fiber bundle $Z\equiv\M\times_{\G} (\Sigma\setminus\{z_\infty\})
\rightarrow\M^\dag$. Then as shown in our earlier paper \cite{C2}, there 
is a canonically defined smooth map of orbifolds in the sense of \cite{C1},
$\mbox{Ev}: Z\rightarrow X$, such that the induced map between the 
underlying spaces is the evaluation map $[(f,z)]\mapsto f(z)$, 
$\forall f\in\M, z\in\Sigma\setminus\{z_\infty\}$. 

The map $\mbox{Ev}: Z\rightarrow X$ is a diffeomorphism of orbifolds 
onto $X\setminus\{p_0\}$. In fact, as in the proof of Lemma 4.3 in 
\cite{C2}, one can show that the differential of $\mbox{Ev}$ is invertible
and that the induced map between the underlying spaces is onto 
$X\setminus\{p_0\}$. It remains to see that the induced map of $\mbox{Ev}$
between the underlying spaces is injective. This is because: (1) for each 
$f\in\M$, the $J$-holomorphic curve $\mbox{Im }f$ is a quasi-suborbifold, 
and (2) for any $f_1,f_2\in\M$ which have different orbits in $\M^\dag
\equiv\M/\G$, the $J$-holomorphic curves $\mbox{Im }f_1,\mbox{Im }f_2$ 
intersect only at $p_0$. The former is proved in Lemma 3.3. To see the 
latter, suppose for simplicity that $f_1,f_2\in\M$ are not multiply covered. 
Then by the intersection formula,
the contribution of $p_0$ to the intersection product 
$\mbox{Im }f_1\cdot \mbox{Im }f_2$ is at least
$$
\frac{\frac{|G|}{2m}\cdot\frac{|G|}{2m}}{|G|}=\frac{|G|}{4m^2}
=c_1(E)\cdot c_1(E),
$$
which implies that $\mbox{Im }f_1,\mbox{Im }f_2$ can not intersect at
any other point. The discussion for the remaining cases is similar, so
we leave the details to the reader. Hence the claim.

Let $M_0\equiv \mbox{Ev}^{-1}(C_0)$ be the inverse image of $C_0$ in $Z$.
Then $M_0$ is a suborbifold in $Z$. Moreover, since for each $f\in\M$, 
the $J$-holomorphic curve $\mbox{Im }f$ intersects $C_0$
at exactly one point, we see that $M_0$ is a smooth section of the orbifold
fiber bundle $Z\rightarrow \M^\dag$. Consequently, we may regard $Z$ as an
orbifold complex line bundle over $M_0$. Note that under
$\mbox{Ev}:Z\rightarrow X$, $M_0$ is mapped diffeomorphically onto
$C_0\subset X$.

One can show, by an identical argument as in \cite{C2}, that there exists
a regular neighborhood $N_0$ of the singular point $p_0$ in $X$, such that
for any $f\in\M$, $\partial N_0$ intersects $\mbox{Im }f$
transversely at a simple closed loop. It follows easily that
$X\setminus int(N_0)$ is diffeomorphic to the associated unit disc bundle
of $Z\rightarrow M_0$ via the inverse of $\mbox{Ev}$, under which $C_0$
is mapped diffeomorphically onto the $0$-section $M_0$. Now observe that the
$s$-cobordism $W$ is diffeomorphic to $X\setminus int(N_0)$ with a
regular neighborhood of $C_0$ removed. It follows easily that $W$ is
diffeomorphic to the product $(\s^3/G)\times [0,1]$.

\hfill $\Box$

\vspace{2mm}

\noindent{\bf Proof of Lemma 3.9}

\vspace{2mm}

The basic observation here is that if a component $C_i$ in the
Poincar\'{e} dual of $c_1(E)$ has a relatively small self-intersection
$C_i\cdot C_i$, then one can easily show that $C_i=\mbox{Im }f$ for
some $f\in\M$. In particular, $\M$ is nonempty
when $c_1(E)\cdot c_1(E)=\frac{|G|}{4m^2}$ is sufficiently small.
On the other hand, in the cases where $c_1(E)\cdot c_1(E)=\frac{|G|}{4m^2}$
is not small, it turns out that $d(E)$, the dimension of the Seiberg-Witten
moduli space, is also considerably large, which allows us to break the
Poincar\'{e} dual of $c_1(E)$ into smaller pieces by requiring it to pass
through a certain number of specified points (cf. Remark 2.3).

\vspace{2mm}

{\bf Case 1}. $|G|<4m^2$. Let $\{C_i\}$ be the set of $J$-holomorphic curves
obtained by applying Theorem 2.2 (2) to $E$. The assumption
$|G|<4m^2$ has the following immediate consequences: (1)
$C_0$ is not contained in $\{C_i\}$ because $C_0\cdot C_0
=\frac{4m^2}{|G|}>1$ and $c_1(E)\cdot c_1(E)=\frac{|G|}{4m^2}<1$, and
(2) if let $C_i=r_i\cdot c_1(E)$ for some $0<r_i\leq 1$, then the
virtual genus
\begin{eqnarray*}
g(C_i) & = & \frac{1}{2}(r^2_i\cdot c_1(E)\cdot c_1(E)+r_i\cdot c_1(K_X)
\cdot c_1(E))+1\\
       & = & \frac{1}{2}(r^2_i\cdot 
\frac{|G|}{4m^2}-r_i\cdot \frac{m+1}{m})+1<1.
\end{eqnarray*}
As corollaries of (2), we note that for any $f_i:\Sigma_i\rightarrow X$
parametrizing $C_i$, $g_{|\Sigma_i|}=0$
because $g_{|\Sigma_i|}\leq g_{\Sigma_i}\leq g(C_i)<1$.
(Here $g_{|\Sigma_i|}$ is the genus of the underlying Riemann surface
of $\Sigma_i$.) Furthermore, note that $p_0\in\cup_i C_i$ because
in the present case, the representation $\rho:G\rightarrow\s^1$
defined in Lemma 3.6 is nontrivial, cf. Remark 2.3. If $C\in\{C_i\}$
is a compoment containing $p_0$, then $f^{-1}(p_0)$ consists of only
one point for any $f:\Sigma\rightarrow X$ parametrizing $C$. This is
because for any $z^\prime\in f^{-1}(p_0)$ with order $m^\prime\geq 1$,
the contribution from $z^\prime$ to $g_\Sigma$
is $\frac{1}{2}(1-\frac{1}{m^\prime})$, and $k_{z^\prime}\geq
\frac{1}{2m^\prime}(\frac{|G|}{m^\prime}-1)\geq \frac{1}{2m^\prime}$.
Hence the contribution from each point in $f^{-1}(p_0)$ to the right
hand side of the adjunction formula for $C$ is least $\frac{1}{2}$.
If there were more than one point in $f^{-1}(p_0)$, the right hand
side of the adjunction formula would be no less than $1$, which is
a contradiction to $g(C)<1$.

With these understood, note that $d(E)\geq 2$ by Lemma 3.8, so that we
may require that $\cup_i C_i$ also contains a smooth point $p\in C_0$.
It follows easily, since $C_0$ is not contained in $\{C_i\}$, that
$\{C_i\}$ consists of only one component, denoted
by $C$, which contains both $p\in C_0$ and $p_0$. Let $f:\Sigma
\rightarrow X$ be a parametrization of $C$. Then as we argued earlier,
$f^{-1}(p_0)$ consists of only one point, say $z_\infty$. Moreover,
$f^{-1}(C_0)$ also consists of only one point, say $z_0$, because
$C,C_0$ intersect at a smooth point $p$ and $C\cdot C_0=1$. It follows
easily that the link of $p$ in $C$ is homotopic in $\s^3/G$ to a regular
fiber of the Seifert fibration, which has homotopy class $\mu_{2m}^{-1} I
\in Z_{2m}$. On the other hand, $g_{|\Sigma|}=0$, so that the link of
$p$ in $C$ is homotopic in $W$ to the inverse of the link of $p_0$ in
$C$. Hence the link of $p_0$ in $C$ must have homotopy class $\mu_{2m}
I\in Z_{2m}$ (cf. Lemma 3.5), from which it is easily seen that
$f\in\M$. This proves that $\M$ is nonempty
when $|G|<4m^2$.

\vspace{2mm}

{\bf Case 2}. $|G|>4m^2$. The proof will be done case by case 
according to the type of $G$.

(1) $G=\langle Z_{2m},Z_{2m};\widetilde{D}_n,\widetilde{D}_n\rangle$
or $\langle Z_{4m},Z_{2m};\widetilde{D}_n,C_{2n}\rangle$. In this
case, note that $|G|>4m^2$ is equivalent to $m<n$. We start with
the following

\vspace{3mm}

\noindent{\bf Sublemma 3.10}{\em\hspace{3mm}
Let $C$ be a $J$-holomorphic curve which intersects $C_0$ at only
one singular point. If furthermore, {\em (1)} $C\cdot C_0<1$ when
the singular point in $C\cap C_0$ is of order $2$, and {\em (2)}
$C$ contains $p_0$ when the singular point in $C\cap C_0$ is of
order $n$. Then $C$ is the image of a member of $\M$.
}
\vspace{2mm}

\pf
Let $f:\Sigma\rightarrow X$ be a parametrization of $C$.

First, consider the case where the singular point in $C\cap C_0$,
say $p_1$, has order $2$. Note that $C\cdot C_0<1$ implies that
$f^{-1}(C_0)$ consists of only one point, say $z_0\in\Sigma$, which
has order $2$ in $\Sigma$, and in a local representative $(f_0,\rho_0)$
of $f$ at $z_0$, $\rho_0(\mu_2)=\mu_2$ and $f_0(z)=(a(z^l+\cdots),z)$
with $l$ odd if $a\neq 0$. In particular, the link of $p_1$
in $C$ is homotopic in $\s^3/G$ to the singular fiber of the
Seifert fibration at $p_1$. Now recall that $H_1(\s^3/G;\Z)
=\Z_{m}\oplus\Z_2\oplus\Z_2$ if $n$ is even, and $H_1(\s^3/G;\Z)=\Z_{4m}$
when $n$ is odd, where, if we let $x,y$ be the standard generators of
$\widetilde{D}_n$ with relations $x^2=y^n=(xy)^2=-1$, one of the factor
in $\Z_2\oplus\Z_2$ in the former case is generated by $x$ and the other
by $y$, and in the latter case, the generator of $Z_{4m}$ is the class of
$\mu_{2m}x$ or $\mu_{4m}x$, depending on whether $m$ is odd or even.
With this understood, note that the class in $H_1(\s^3/G;\Z)$ of the
link of $p_1$ in $C$ projects nontrivially to the $\Z_2$ factor generated
by $x$ in the former case, and is a generator of $Z_{4m}$ in the latter
case. It follows easily that $f^{-1}(p_0)$ is nonempty, and there must be
a $z_\infty\in f^{-1}(p_0)$, such that the pushforward of the link of
$z_\infty$ in $\Sigma$ under $f$ has a homology class in $\s^3/G$ which
projects nontrivially onto the $\Z_2$ factor generated by $x$ in the former
case, and is a generator of $Z_{4m}$ in the latter case. In any event, the
order of $z_\infty$ in $\Sigma$ must be $4m$ or less, and as argued in the
proof of Lemma 3.3, $C$ is a quasi-suborbifold, and is easily seen to be the
image of a member of $\M$.

Next we suppose that the singular point in $C\cap C_0$ is $p_3$, which has
order $n$. Note that $m<n$ implies that the normalized Seifert invariant
at $p_3$ is $(n,m)$. Let $(w_1,w_2)$ be a holomorphic coordinate system on
a local uniformizing system at $p_3$, where $C_0$ is given locally by $w_2=0$,
and the singular fiber of the Seifert fibration at $p_3$ is defined by
$w_1=0$, $|w_2|\equiv\mbox{constant}$, and the $\Z_n$-action is given
by $\mu_n\cdot (w_1,w_2)=(\mu_nw_1,\mu_n^m w_2)$. Let $f^{-1}(C_0)=
\{z_i\mid i=1,2,\cdots,k\}$ where each $z_i$ has order $m_i\geq 1$, and
let $(f_i,\rho_i)$ be a local representative of $f$ at $z_i$, where
$\rho_i(\mu_{m_i})=\mu_{m_i}^{r_i}$, with $0\leq r_i<m_i$, $r_i,m_i$
relatively prime, and $f_i(w)=(c_i(w^{l_i^\prime}+\cdots),w^{l_i})$
such that $c_i\neq 0$ unless $m_i=n$ and $l_i=1$. Note that $f_i$ being
$\rho_i$-equivariant implies that $l_i\equiv mr_i\pmod{m_i}$, and
when $c_i\neq 0$, $l_i^\prime\equiv r_i\pmod{m_i}$. By the intersection
formula, the contribution from $z_i$ to $C\cdot C_0$ is
$\frac{(n/m_i)l_i}{n}=\frac{l_i}{m_i}$. Hence
$C\cdot C_0=\sum_{i=1}^k\frac{l_i}{m_i}$, and the virtual genus of $C$ is
$$
g(C)=\sum_{i,j=1}^k\frac{l_il_j}{m_im_j}\cdot\frac{n}{2m}-\sum_{i=1}^k
\frac{l_i}{m_i}\cdot\frac{m+1}{2m}+1.
$$
Evidently, the contribution to $g(C)$ from each $z_i$ is
$$
L_i\equiv \frac{l_i^2}{m_i^2}\cdot\frac{n}{2m}-
\frac{l_i}{m_i}\cdot\frac{m+1}{2m},
$$
and the contribution from each unordered pair $[z_i,z_j]$, $i\neq j$, is
$$
L_{[i,j]}\equiv \frac{l_il_j}{m_im_j}\cdot\frac{n}{m}.
$$
On the other hand, the contribution of each $z_i$ to the right hand side
of the adjunction formula is
$$
R_i\equiv\frac{1}{2}(1-\frac{1}{m_i})+k_{z_i},
$$
and the contribution of each unordered pair $[z_i,z_j]$, $i\neq j$, is
$$
R_{[i,j]}\equiv k_{[z_i,z_j]}.
$$

In order to estimate $k_{z_i}$ and $k_{[z_i,z_j]}$, we next recall some
basic facts about the local self-intersection number and local
intersection number of $J$-holomorphic curves, cf. \cite{C2} and
the references therein.
\begin{itemize}
\item Let $C$ be a holomorphic curve in $\C^2$ parametrized by
$f(z)=(a(z^{l_1}+\cdots), z^{l_2})$, where $f:(D,0)\rightarrow
(\C^2,0)$ is from a disc $D\subset\C$ centered at $0$ such that
$f|_{D\setminus\{0\}}$ is embedded. Then the local self-intersection
number $C\cdot C\geq\frac{1}{2}(l_1-1)(l_2-1)$. Note that the above
inequality still makes sense even if $a=0$ in the formula for $f$,
in which case $l_1$ is undefined. This is because $l_2=1$ by the
assumption that $f|_{D\setminus\{0\}}$ is embedded.
\item Let $C, C^\prime$ be distinct holomorphic curves in $\C^2$ 
parametrized by $f(z)=(a(z^{l_1}+\cdots), z^{l_2})$ and
$f^\prime(z)=(a^\prime(z^{l_1^\prime}+\cdots), z^{l_2^\prime})$
respectively, where $f:(D,0)\rightarrow (\C^2,0)$, $f^\prime:(D,0)
\rightarrow (\C^2,0)$ are from a disc $D\subset\C$ centered at $0$ such
that $f|_{D\setminus\{0\}}$, $f^\prime|_{D\setminus\{0\}}$ are embedded.
Then the local intersection number $C\cdot C^\prime\geq\min (l_1l^\prime_2,
l_2l^\prime_1)$. Here $l_1=\infty$ (resp. $l_1^\prime=\infty$) if $a=0$
(resp. $a^\prime=0$).
\end{itemize}
With the preceding understood and by the definition in \cite{C2}, we
have
$$
k_{z_i}\geq\frac{1}{2m_i}((l_i-1)(l_i^\prime-1)
+(\frac{n}{m_i}-1)l_il_i^\prime),\;
k_{[z_i,z_j]}\geq\frac{1}{n}\cdot\frac{n}{m_i}\cdot\frac{n}{m_j}
\cdot\min(l_il_j^\prime,l_jl_i^\prime).
$$
(Note that the right hand side of the first inequality still makes sense
even when $l_i^\prime$ is undefined, because in this case, $l_i=1$ and
$n=m_i$ must be true.)

Next we shall compare $L_i$ with $R_i$ and $L_{[i,j]}$ with $R_{[i,j]}$.
To this end, we write $l_i^\prime=r_i+t_im_i$ and
$mr_i=l_i+s_im_i$. Here $t_i\geq 0$, and $s_i\geq 0$ if $l_i<m_i$.
When $l_i=m_i$, we must have $l_i=m_i=1$ and $r_i=0$. In this
case, $s_i=-1$ and $l_i^\prime=t_i\geq 1$. It follows easily that
$\min(l_il_j^\prime,l_jl_i^\prime)\geq \frac{l_il_j}{m}$, hence
$R_{[i,j]}\geq L_{[i,j]}$ for all $i\neq j$. To compare $L_i$ with $R_i$,
we note that
\begin{eqnarray*}
k_{z_i} & \geq & \frac{1}{2m_i}((l_i-1)(l_i^\prime-1)+
(\frac{n}{m_i}-1)l_il_i^\prime)\\
& = & \frac{1}{2m_i}(1-l_i-(r_i+t_im_i)+\frac{nl_i}{m_i}
(\frac{l_i+s_im_i}{m}+t_im_i))\\
& = & \frac{1}{2m_i}+\frac{nl_i^2}{2m_i^2 m}-\frac{l_i(m+1)}{2m_im}+
\frac{nl_i-m_i}{2m_i m}(s_i+t_im),\\
\end{eqnarray*}
which easily gives $R_i-L_i\geq \frac{1}{2}+\frac{nl_i-m_i}{2m_i m}
(s_i+t_im)$.

With the above estimates in hand, now observe that $f^{-1}(p_0)$ is not 
empty
by the assumption, so that, as we argued earlier, the contribution of
$f^{-1}(p_0)$ to the right hand side of the adjunction formula is at
least $\frac{1}{2}$. It follows easily that $f^{-1}(C_0)$ contains
only one point (i.e. $k=1$), and that either $m_1=n$ with $l_1=1$, or
$s_1+t_1m=0$, which means either $s_1=t_1=0$, or $m=t_1=-s_1=1$ with
$m_1=l_1=1$. Moreover, $g_{|\Sigma|}=0$, and $f^{-1}(p_0)$ contains
only one point, say $z_\infty$, with the contribution from $z_\infty$
to the right hand side of the adjunction formula being exactly
$\frac{1}{2}$. The last point particularly implies that the order
of $z_\infty$ is $\frac{1}{2}|G|=2mn$.

The case where $s_1=t_1=0$ but $nl_1\neq m_1$ or $l_1=m_1=1$ can be
ruled out as follows. Consider $s_1=t_1=0$ but $nl_1\neq m_1$ first.
Note that $l_1=mr_1$ must be true in this case. As we have seen earlier,
the link of $p_3$ in $C$ is homotopic in $\s^3/G$ to $l_1\cdot\frac{n}{m_1}
=\frac{r_1mn}{m_1}$ times of the singular fiber of the Seifert fibration
at $p_3$, which has order $2mn$. It is easily seen that the order $d$ of
the link of $p_3$ is either divisible by $m_1$, in which case $r_1$ is even,
or divisible by $2m_1$, in which case $r_1$ is odd. In any event, $d<2mn$
if $nl_1\neq m_1$. On the other hand, $g_{|\Sigma|}=0$ implies that the
link of $p_0$ in $C$ is homotopic, in the $s$-cobordism $W$ and through $C$,
to the inverse of the link of $p_3$ in $C$. But the homotopy class of the
link of $p_0$ has the same order in $\pi_1(\s^3/G)$ as the order of $z_\infty$,
which is $2mn$. This is a contradiction. The discussion for the case 
where $l_1=m_1=1$
is similar. In this case, the link of $p_3$ is homotopic in $\s^3/G$
to $n$ times of the singular fiber, hence has order at most $2m$, which is
less than $2mn$.  This is also a contradiction. Hence $m_1=n$ and $l_1=1$, 
and the adjunction
formula implies that $C$ is a quasi-suborbifold. Moreover, it is easily
seen that $C$ is the image of a member of $\M$. Hence
Sublemma 3.10.

\hfill $\Box$

Now back to the proof. Let $\{C_i\}$ be the $J$-holomorphic
curves which are obtained by applying Theorem 2.2 (2) to $E$.
Set $N\equiv \frac{d(E)}{2}$ when $m\neq 1$ and $N\equiv
\frac{d(E)}{2}-1$ when $m=1$. Then by Remark 2.3, we can
specify any $N$ distinct smooth points $q_1,\cdots,q_N\in
X\setminus C_0$ and require that $q_1,\cdots,q_N\in \cup_i C_i$,
and moreover, $p_0\in \cup_i C_i$. Now we let $q_1\equiv q_{1,j}$
be a sequence of points converging to a smooth point $q$ in
$C_0$, while keeping $q_2,\cdots, q_N$ fixed, and let $\{C_i^{(j)}\}$
be the corresponding sequence of (sets of) $J$-holomorphic curves.
By passing to a subsequence if necessary, we may assume that the number
of components in $\{C_i^{(j)}\}$ is independent of $j$, and each
$C_i^{(j)}$ is parametrized by a $J$-holomorphic map $f_{i,j}:
\Sigma_i\rightarrow X$ from an orbifold Riemann surface
independent of $j$ (note that the complex structure on
$\Sigma_i$ is allowed to vary). This follows readily from
the fact that $C_i^{(j)}\cdot C_0\geq\frac{1}{n}$, and that the virtual
genus $g(C_i^{(j)})$, hence the corresponding orbifold genus, is uniformly
bounded from above. By the Gromov compactness theorem (cf. \cite{CR}), each
$f_{i,j}$ converges to a cusp-curve $f^\prime_i:\Sigma^\prime_i\rightarrow X$.
The upshot here is that we can always manage to have $C_0$ contained
in $\cup_i\mbox{Im }f_i^\prime$, or else $\M$ is nonempty.
Accept this momentarily, we note as a consequence that
$$
C\cdot C_0\leq 1-C_0\cdot C_0=1-\frac{m}{n},
$$
where $C\subset \cup_i\mbox{Im }f_i^\prime$ is any component
containing $p_0$. By letting $q_2,\cdots$, $q_N$ converge to a smooth
point in $C_0$ one by one, we have at the end 
$$
C\cdot C_0\leq 1-N\cdot\frac{m}{n}
$$
for any component $C$ in the limiting cusp-curve that contains $p_0$.
Now observe that $N\geq\frac{\delta+1}{2}$ when $m\neq 1$, where
$n=\delta m+r$ with $0\leq r\leq m-1$, and $N\geq \frac{n-1}{2}$ when $m=1$.
It follows easily that $C\cdot C_0<\frac{1}{2}$ when $m\neq 1$ and
$C\cdot C_0\leq \frac{1}{2}+\frac{1}{2n}$ when $m=1$. Clearly, $C$ can
only intersect $C_0$ at one singular point, because otherwise, we would
have $C\cdot C_0\geq \frac{1}{2}+\frac{1}{n}$. By Sublemma 3.10,
$\M$ is nonempty.

It remains to show that $C_0\subset \cup_i\mbox{Im }f_i^\prime$.
Note that if the component $C^{(j)}\in\{C_i^{(j)}\}$ which contains
$q_{1,j}$ intersect $C_0$ at a singular point, which is the case
when $C^{(j)}\cdot C_0<1$, then it is clear that $C_0$ must be one
of the component in the limiting cusp-curve of $\{C^{(j)}\}$. If $C^{(j)}$
intersect $C_0$ at a smooth point $q_j$, then $\lim_{j\rightarrow\infty}
q_j=q=\lim_{j\rightarrow\infty}q_{1,j}$ must hold if $C_0$ is not
contained in the limiting curve, and moreover, the limiting curve must
contain only one nonconstant component, which intersect $C_0$ at $q$
transversely. In this case we let $q\equiv q_k$ be a sequence of smooth
points on $C_0$ which converges to the singular point $p_3\in C_0$ of
order $n$. Let $C_k$ be the corresponding $J$-holomorphic curves, which
we assume to be parametrized by $f_k:\Sigma\rightarrow X$ from a fixed
orbifold Riemann surface without loss of generality. If the limiting curve
$f^\prime:\Sigma^\prime\rightarrow X$ of $\{C_k\}$ intersect $C_0$
only at $p_3$, then $\M$ is nonempty by Sublemma 3.10. If the
limiting curve $f^\prime$ has a nonconstant component, denoted by $f_\nu\equiv
f^\prime|_{\Sigma_\nu}:\Sigma_\nu\rightarrow X$, which intersect $C_0$
at a singular point of order $2$, say $p_1$, then there must be a simple
closed loop $\gamma\subset\Sigma$ collapsed to $p_1$ during the
convergence. Note that there are two distinct points $z_\nu, z_\omega
\in\Sigma^\prime$, where $z_\nu\in\Sigma_\nu$ with $f_\nu(z_\nu)=p_1$,
which are the image of $\gamma$ under the collapsing $\Sigma\rightarrow
\Sigma^\prime$. Let $\Sigma_\omega$ be the component of $\Sigma^\prime$
which contains $z_\omega$ (here $\Sigma_\omega=\Sigma_\nu$ is allowed). 
Then one of the following must be true:
(a) either $f^\prime$ is nonconstant on $\Sigma_\omega$, or $f^\prime$
is constant on $\Sigma_\omega$ but there exists a component
$\Sigma_\omega^\prime$ and a point $z_\omega^\prime\neq z_\nu$
such that $f^\prime$ is nonconstant on $\Sigma_\omega^\prime$ and
$f^\prime(z_\omega^\prime)=p_1$, (b) the simple closed loop $\gamma$
bounds a sub-surface $\Gamma$ in $\Sigma$ which contains no orbifold
points and is of nonzero genus, such that $f^\prime$ is constant on
every component in the image of $\Gamma\subset\Sigma\rightarrow
\Sigma^\prime$. However, the latter case can be ruled out as follows.
According to the Gromov compactness theorem (cf. \cite{CR}), if
we fix a sufficiently small $\epsilon>0$, then there exists a regular
neighborhood of $\gamma$ in $\Sigma$, identified with $\gamma\times [-1,1]$,
such that (1) $\gamma\times\{-1\}$ is mapped to the link of $z_\nu$ of
radius $\epsilon$ in $\Sigma_\nu$ under $\Sigma\rightarrow\Sigma^\prime$,
(2) $f_k$ converges to $f^\prime$ in $C^\infty$ on $\gamma\times\{-1,1\}$,
(3) the diameter of $f_k(\gamma\times [-1,1])$ is less than $\epsilon$
when $k$ is sufficiently large, and (4) because $f^\prime$ is constant on
every component in the image of $\Gamma\subset\Sigma\rightarrow
\Sigma^\prime$, the diameter of $f_k(\Gamma\setminus\gamma\times [0,1))$
is also less than $\epsilon$ when $k$ is sufficiently large. Let
$U(10\epsilon)$ be the regular neighborhood of $p_1$ in $X$ of radius
$10\epsilon$. Then it is clear that when $k$ is sufficiently large,
$f_k(\Gamma)\subset U(10\epsilon)\setminus\{p_1\}$ and $f_k(\gamma)
=f_k(\partial\Gamma)$ is homotopic in $U(10\epsilon)\setminus\{p_1\}$
to the push-forward of the link of $z_\nu$ in $\Sigma_\nu$ under $f_\nu$.
But this is impossible because (1) $f_\nu$ is clearly not multiply
covered, and (2) the link of $p_1$ in $\mbox{Im }f_\nu$, which is
the push-forward of the link of $z_\nu$ in $\Sigma_\nu$ under $f_\nu$
because of (1), is not null-homologous in $U(10\epsilon)\setminus\{p_1\}
\cong\R\P^3\times (0,1]$. Hence the latter case (i.e. case (b)) is
ruled out. On the other hand, the former case (i.e. case (a))
is also impossible if $C_0$ is not contained in $\cup_i\mbox{Im }
f_i^\prime$. This is because each of $z_\nu$ and $z_\omega$ (or
$z_\omega^\prime$) will contribute at least $\frac{1}{2}$ to
$[f(\Sigma^\prime)]\cdot C_0$, and with the contribution from $p_3$,
we would have $[f(\Sigma^\prime)]\cdot C_0>1$, which is a contradiction. 
Hence we can always manage to have $C_0\subset\cup_i\mbox{Im }f_i^\prime$,
or else $\M$ is nonempty.

(2) $G=\langle Z_{2m},Z_{2m};\widetilde{T},\widetilde{T}\rangle$
or $\langle Z_{6m},Z_{2m};\widetilde{T},\widetilde{D}_2\rangle$.
Note that $|G|>4m^2$ is equivalent to $m<6$, which means that
$m=1$ or $5$ in the former case, and $m=3$ in the latter case.

First, observe that if $C$ is a $J$-holomorphic curve parametrized by
$f:\Sigma\rightarrow X$ such that $p_0\in C$, then each $z\in f^{-1}(p_0)$
will contribute at least $\frac{1}{2}+\frac{1}{6m}$ to the right hand side
of the adjunction formula. To see this, let $m_0$ be the order of $z$.
Then the total contribution of $z$ is
$$
\frac{1}{2}(1-\frac{1}{m_0})+k_z\geq \frac{1}{2}(1-\frac{1}{m_0})+
\frac{1}{2m_0}(\frac{24m}{m_0}-1)=\frac{1}{2}+\frac{1}{2m_0}
(\frac{24m}{m_0}-2).
$$
Note that $m_0\leq 6m$. Hence the claim.

We consider the case where $m=3$ or $5$ first. Let $\{C_i\}$ be the
$J$-holomorphic curves obtained by applying Theorem 2.2 (2) to $E$.
Note that $p_0\in\cup_i C_i$, and because $d(E)=2$, we can specify a
smooth point $q\in C_0$ and require that $q\in \cup_i C_i$. We claim
that there exists a $J$-holomorphic curve $\hat{C}$ such that
$\hat{C}\cdot C_0\leq\frac{1}{2}$ and $p_0\in\hat{C}$.

Let $C\in\{C_i\}$ be a component containing $p_0$. If $C$ does not
contain $q$, then the component in $\{C_i\}$ which contains $q$ must
be $C_0$. As a consequence, $C\cdot C_0\leq 1-C_0\cdot C_0=1-\frac{m}{6}
\leq\frac{1}{2}$.

Now suppose $q\in C$. Then $C$ must be the only component in $\{C_i\}$,
and $C,C_0$ intersect transversely at the smooth point $q$. We let
$q\equiv q_k$ be a sequence of smooth points on $C_0$ converging to
the singular point $p_3\in C_0$ of order $3$, and let $C_k$ be the
corresponding $J$-holomorphic curves, which we assume without loss of
generality to be parametrized by $f_k:\Sigma\rightarrow X$ from a fixed
orbifold Riemann surface. Note that $g_{|\Sigma|}=0$ because
$g(C_k)=\frac{5-m}{2m}+1\leq 1+\frac{1}{3}$, and because $C_k$ contains
$p_0$ so that $f^{-1}_k(p_0)$ contributes at least $\frac{1}{2}$ to
the right hand side of the adjunction formula. Similarly, $f^{-1}_k(p_0)$
contains at most two points. Let $z_0=f^{-1}_k(q_k)$. Note that $z_0$
is a regular point of $\Sigma$.

By the Gromov compactness theorem, a subsequence of $\{f_k\}$, after
reparametrization if necessary, will converge either in $C^\infty$ to
$f:\Sigma\rightarrow X$, or to a cusp-curve $f:\Sigma^\prime\rightarrow X$.
If the convergence is in $C^\infty$, then $f$ must be multiply covered,
because otherwise, we will have $k_{z_0}\geq\frac{1}{2}(\frac{3}{1}-1)=1$,
and together with the contribution of $f^{-1}(p_0)$ which is at least
$\frac{1}{2}$, it would imply that the right hand side of the adjunction 
formula is no less than $1+\frac{1}{2}$, which is greater than the left 
hand side $g(\mbox{Im }f)=\frac{5-m}{2m}+1\leq 1+\frac{1}{3}$. This is 
a contradiction. For a multiply covered $f$, it is clear that 
$\hat{C}\equiv\mbox{Im }f$
is a $J$-holomorphic curve which obeys $\hat{C}\cdot C_0\leq\frac{1}{2}$
and $p_0\in \hat{C}$.

Now suppose $f_k$ converges to a cusp-curve $f:\Sigma^\prime\rightarrow X$.
For technical convenience, we shall regard $z_0$ as a marked point so
that $z_0$ will not lie on a collapsing simple closed loop during the
Gromov compactification (cf. \cite{CR}).
Let $\Sigma_\nu$ be the component of $\Sigma^\prime$ which contains $z_0$.
First, we consider the case where $f$ is nonconstant over $\Sigma_\nu$.
Note that under this assumption, if $\mbox{Im }f|_{\Sigma_\nu}\neq C_0$,
one can easily show, because $z_0\in\Sigma_\nu$ is a regular point, that
$[f(\Sigma_\nu)]\cdot C_0=1$. This implies that $\Sigma_\nu$ is the only
component of $\Sigma^\prime$ over which $f$ is nonconstant. In particular,
$p_0\in\mbox{Im }f|_{\Sigma_\nu}$ because $\Sigma^\prime$ is connected.
Then as we argued in the preceding paragraph, $f$ must be multiply covered
over $\Sigma_\nu$, and $\mbox{Im }f|_{\Sigma_\nu}$ is the $J$-holomorphic
curve that we are looking for. If $\mbox{Im }f|_{\Sigma_\nu}=C_0$, then
$\hat{C}\cdot C_0\leq 1-\frac{m}{6}\leq\frac{1}{2}$ for any (nonconstant)
component $\hat{C}$ in $\mbox{Im }f$ that contains $p_0$. Now suppose $f$
is constant over $\Sigma_\nu$. Then it follows easily, because
$g_{|\Sigma|}=0$ and because
no simple closed loops bounding a disc $D\subset\Sigma$ will collapse if
$f_k(D)$ lies in the complement of $C_0$ and $p_0$ (cf. case (2) in the
proof of Lemma 3.4), it follows easily that each $f_k^{-1}(p_0)$ contains
exactly two points $z_{\infty}^{(1)}$, $z_{\infty}^{(2)}$, and there are
simple closed loops $\gamma_1,\gamma_2\subset\Sigma$, each bounding a disc
$D\subset\Sigma$, such that (1) $D$ contains exactly one of
$z_{\infty}^{(1)}, z_{\infty}^{(2)}$ but not $z_0$, (2) no simple closed
loops in $D$ collapsed, (3) $D$ is mapped to a component $\Sigma_\omega$
under $\Sigma\rightarrow\Sigma^\prime$ such that $f|_{\Sigma_\omega}
\neq\mbox{constant}$. It is clear that one of these
two components of $\mbox{Im }f$, which all contains $p_0$, will have
intersection product with $C_0$ no greater than $\frac{1}{2}$.

Let $C$ be a $J$-holomorphic curve such that $p_0\in C$ and
$C\cdot C_0\leq\frac{1}{2}$, which we have just shown to exist, and let
$f:\Sigma\rightarrow X$ be a parametrization of $C$. Note that $g(C)<1$,
so that $g_{|\Sigma|}=0$ and $f^{-1}(p_0)$ consists of only one point
$z_\infty$. Moreover, it follows easily from $C\cdot C_0\leq\frac{1}{2}$
that $f^{-1}(C_0)$ consists of only one point also, and either
$C\cdot C_0=\frac{1}{2}$ or $\frac{1}{3}$, and if $p_i$ of order $a_i$
for some $i=1,2$ or $3$ is the singular point where $C,C_0$ intersect,
then the link of $p_i$ in $C$ is homotopic in $\s^3/G$ to the singular
fiber of the Seifert fibration at $p_i$, which has order $2ma_i$ in
$\pi_1(\s^3/G)=G$. It implies that the order of $z_\infty$ is also
$2ma_i$ because $g_{|\Sigma|}=0$, and the adjunction formula implies
that $C$ is a quasi-suborbifold, which is easily seen to be the image
of a member of $\M$.

It remains to consider the case where $m=1$. We will show in this case
that there is also a $J$-holomorphic curve $C$ such that $p_0\in C$
and $C\cdot C_0\leq\frac{1}{2}$. To see this, let $\{C_i\}$ be the
$J$-holomorphic curves obtained by applying Theorem 2.2 (2) to
$E$. Since $d(E)=8$ in the present case, we can specify any $3$
distinct smooth points $q_1,q_2,q_3$, where $q_1\in C_0$ and $q_2,q_3\in
X\setminus C_0$, and require that they are contained in $\cup_i C_i$,
and moreover, we require $p_0\in \cup_i C_i$ also. Let $q_2
\equiv q_{2,j}$ be a sequence of points converging to a smooth point
$q_2^\prime\in C_0$ such that $q_2^\prime\neq q_1$, and we denote by
$\{C_i^{(j)}\}$ the corresponding sequence of $J$-holomorphic curves.
If $C_0\in \{C_i^{(j)}\}$, then as we have seen earlier, the components
containing $q_{2,j}$ will converge to a cusp-curve during which a component
will be split off which goes to $C_0$. At this stage, the component
$\hat{C}$ in the limiting cusp-curve which contains $p_0$ obeys
$\hat{C}\cdot C_0\leq 1-2C_0\cdot C_0$. We claim that this is also true
even when $C_0$ is not contained in $\{C_i^{(j)}\}$. To see this, note
that under this assumption, there is only one component, denoted by
$C^{(j)}$, in $\{C_i^{(j)}\}$, which intersects $C_0$ transversely at
$q_1$. Let $f_j:\Sigma\rightarrow X$ be a parametrization of $C^{(j)}$,
which is from a fixed orbifold Riemann surface by passing to a
subsequence if necessary. Since $q_2^\prime=\lim_{j\rightarrow\infty}
q_{2,j}\neq q_1$, $f_j$ has to converge to a cusp-curve $f:\Sigma^\prime
\rightarrow X$. Let $\Sigma_\nu$ be a component of $\Sigma^\prime$
such that $f_\nu\equiv f|_{\Sigma_\nu}$ is nonconstant and $q_1\in
\mbox{Im }f_\nu$. Then since $q_1$ is a smooth point, and $\mbox{Im }f$
contains another smooth point $q_2^\prime=\lim_{j\rightarrow\infty}
q_{2,j}\neq q_1$, we see that $\mbox{Im }f_\nu=C_0$ must be true.
If the degree of  $f_\nu:\Sigma_\nu\rightarrow X$ is at least $2$,
then the claim is clearly true. Suppose $f_\nu:\Sigma_\nu\rightarrow X$
is of degree $1$. Then $\Sigma_\nu$ contains three orbifold points
$z_1,z_2,z_3$, with $f_\nu(z_i)=p_i$ for $i=1,2,3$. Note that $z_1,z_2,
z_3$ must be the result of collapsing $3$ simple closed loops under
$\Sigma\rightarrow\Sigma^\prime$. Consequently, there are components
$\Sigma_i$ of $\Sigma^\prime$, where $i=1,2$ and $3$, such that
$\Sigma_i\neq\Sigma_\nu$ and there are $z_i^\prime\in\Sigma_i$ satisfying
$f(z_i^\prime)=p_i$ (here $\Sigma_i=\Sigma_{i^\prime}$
is allowed). The key observation is that one of $C^{(i)}\equiv
\mbox{Im }f_{\Sigma_i}$ must be either $C_0$ or constant, because 
otherwise, one has $c_1(E)\cdot C_0\geq \sum_{i=1}^3 C^{(i)}\cdot C_0
\geq\sum_{i=1}^3\frac{1}{a_i}>1$ (here $a_i$ is the order of $p_i$, with
$a_1=2$, $a_2=a_3=3$), which is a contradiction. But as we have seen earlier,
none of $C^{(i)}$ is constant. (Suppose $C^{(i)}$ is constant for
some $i=1,2$ or $3$, then $z^\prime_i$ must be the image of a collapsing
simple closed loop $\gamma$ which bounds a sub-surface $\Gamma\subset
\Sigma$ such that for sufficiently large $j$, $f_j(\Gamma)$ lies in
a regular neighborhood of $p_i$ with $p_i$ removed, which is diffeomorphic
to the product of a lens space with $(0,1]$. But on the other hand,
$f_j(\gamma)=f_j(\partial\Gamma)$ is homotopic in the punctured
neighborhood of $p_i$ to the link of $p_i$ in $C_0$, which is not
null-homologous in the punctured neighborhood. This is a contradiction.)
Hence one of $C^{(i)}$ is $C_0$, and the claim follows. Now we let $q_3$
converge to a smooth point in $C_0$, and at the end, it follows easily that
the component ${C}$ in the limiting cusp-curve which contains $p_0$ obeys
${C}\cdot C_0\leq 1-3C_0\cdot C_0=\frac{1}{2}$.

We shall prove that such a $J$-holomorphic curve $C$ is the image
of a member of $\M$ also. Let $f:\Sigma\rightarrow X$ be
a parametrization of $C$. We write $C=\frac{r}{2}\cdot c_1(E)$ for
some $0<r\leq 1$. Then $g(C)=\frac{1}{2}(r^2\cdot\frac{6}{4}
-r\cdot\frac{1+1}{2})+1\leq \frac{1}{4}+1$. Now observe that each
$z\in f^{-1}(C_0)$ will contribute at least $\frac{1}{2}(1-\frac{1}{2})=
\frac{1}{4}$ to the right hand side of the adjunction formula, and
each $z^\prime\in f^{-1}(p_0)$ will contribute at least
$\frac{1}{2}+\frac{1}{6m}=\frac{1}{2}+\frac{1}{6}$.
It follows easily that $f^{-1}(p_0)$ consists of only one point,
and $g_{|\Sigma|}=0$. Finally, observe that $f^{-1}(C_0)$ also
consists of only one point, because each point in $f^{-1}(C_0)$
will contribute at least $\frac{1}{3}$ to $C\cdot C_0$. It is easily
seen from the adjunction formula that $C$ is a quasi-suborbifold, and
moreover, it is the image of a member of $\M$.

(3) $G=\langle Z_{2m},Z_{2m};\widetilde{O},\widetilde{O}\rangle$.
In this case, $|G|>4m^2$ implies that $m=1,5,7$ or $11$. The proof is
similar, although modification is needed at a few places.

First of all, observe that the largest order of an element in $G$ is
$8m$, so that if $C$ is a $J$-holomorphic curve parametrized by
$f:\Sigma\rightarrow X$ with $p_0\in C$, then any point in $f^{-1}(p_0)$
will contribute at least $\frac{1}{2}+\frac{1}{4m}$ to the right hand
side of the adjunction formula.

Consider the case where $m\neq 1$ first. Let $\{C_i\}$ be the
$J$-holomorphic curves obtained by applying Theorem 2.2 (2) to $E$,
where $p_0\in\cup_i C_i$, and a specified smooth point $q\in C_0$
is also contained in $\cup_i C_i$. This time we claim that there exists
a $J$-holomorphic curve $\hat{C}$ such that $p_0\in\hat{C}$ and
$\hat{C}\cdot C_0\leq\frac{7}{12}$, which is clearly true when $q$
is not contained in the same component in $\{C_i\}$ with $p_0$.

Now suppose $C\in\{C_i\}$ is a component containing both $p_0$ and $q$.
Then $C$ must be the only component in $\{C_i\}$, and $C,C_0$ intersect
transversely at the smooth point $q$. We let $q\equiv q_k$ be a sequence
of smooth points on $C_0$ converging to the singular point $p_3\in C_0$
of order $4$, and let $C_k$ be the corresponding $J$-holomorphic curves
parametrized by $f_k:\Sigma\rightarrow X$ from a fixed orbifold Riemann
surface. Note that this time $g_{|\Sigma|}\leq 1$, because $g(C_k)
=\frac{11-m}{2m}+1\leq 1+\frac{3}{5}$, and moreover, $g_{|\Sigma|}=1$
only when $m=5$ and $f^{-1}(p_0)$ consists of only one point. In general,
$f^{-1}_k(p_0)$ contains at most two points. Let $z_0=f^{-1}_k(q_k)$,
which is a regular point of $\Sigma$.

The proof goes in the same way as in the preceding case, except
at the end, we need to consider the following scenario caused by
the possibility that $g_{|\Sigma|}=1$. More concretely, let
$f:\Sigma^\prime\rightarrow X$ be the limiting cusp-curve of $f_k$,
and let $\Sigma_\nu$ be the component of $\Sigma^\prime$ containing
$z_0$. We need to consider the situation where $g_{|\Sigma|}=1$ and
$f$ is constant on $\Sigma_\nu$. Note that if $f$ is constant on
$\Sigma_\nu$, then one of the following must be true: (a) $\Sigma_\nu$
is an orbifold Riemann sphere obtained 
from collapsing two simple closed loops in $\Sigma$, and (b) $\Sigma_\nu$
is an orbifold torus obtained from collapsing one simple closed loop in
$\Sigma$. Let $\Sigma_\omega$ be the component of $\Sigma^\prime$ which
contains $z_\infty$. Then in the latter case, $\Sigma_\omega$ is obtained
by collapsing one simple closed loop, and hence $f$ must be nonconstant
on $\Sigma_\omega$. It can be easily shown that in this case there is a
$J$-holomorphic curve $C$ such that $p_0\in C$ and $C\cdot C_0\leq
\frac{7}{12}$. Now consider the former case. If $f$ is constant on
$\Sigma_\omega$, then $\Sigma_\omega$ must be obtained from collapsing
two simple closed loops, and it is easy to see that there will be at
least two components of $\Sigma^\prime$ over which $f$ is nonconstant.
It is clear that one of these components will give the $J$-holomorphic
curve that we are looking for. Suppose $f$ is nonconstant on
$\Sigma_\omega$. Then one can easily show that we are done in either
one of the following cases: there is a constant component other than
$\Sigma_\nu$, in which case $\Sigma_\omega$ is obtained from collapsing
one simple closed loop, or there is a nonconstant component other than
$\Sigma_\omega$, which will break away with at least $\frac{5}{12}$ or
$\frac{1}{2}$ of the homology. So it remains to consider the case where
$\Sigma_\omega$ is the only component of $\Sigma^\prime$ other than
$\Sigma_\nu$, in which case $\Sigma_\omega$ is an orbifold Riemann sphere
obtained from collapsing two simple closed loops. We are done if
$f_\omega\equiv f|_{\Sigma_\omega}$ is multiply covered, so we assume
that $f_\omega$ is not multiply covered. Set $C_\omega\equiv\mbox{Im }
f_\omega$. Then note that $f^{-1}_\omega(C_0)=\{z_0^{(1)},z^{(2)}_0\}$,
both of which are sent to $p_3$ of order $4$ under $f_\omega$. By the
intersection formula as we have seen earlier, the contribution of each
$z_0^{(i)}$, $i=1,2$, to $C_\omega\cdot C_0$ can be written as
$\frac{l_i}{m_i}$, where $m_i$ is the order of $z_0^{(i)}$ and $l_i,m_i$
are relatively prime. It is clear that either $\frac{l_1}{m_1}
=\frac{l_2}{m_2}=\frac{1}{2}$, or one of them is $\frac{1}{4}$ and the
other is $\frac{3}{4}$, because $C_\omega\cdot C_0=1$. In the former
case, the right hand side of the adjunction formula receives at least $2$
for the contribution from $f^{-1}_\omega(C_0)$, which is more than the
left hand side $g(C_\omega)=\frac{11-m}{2m}+1=1+\frac{3}{5}$. This is a 
contradiction. As for the latter case, note that when $m=5$, the
normalized Seifert invariant at $p_3$ is $(4,1)$ (recall the relation
$m=12b+6+4b_2+3b_3$). Assume without loss of generality that
$\frac{l_1}{m_1}=\frac{3}{4}$. Then it follows easily that a local
representative $(f_1,\rho_1)$ of $f_\omega$ at $z_0^{(1)}$ must obey
$\rho_1(\mu_4)=\mu_4^3$ and $f_1(w)=(c(w^l+\cdots),w^3)$ for some
$c\neq 0$ and some positive interger $l$ satisfying $l\equiv 3\pmod{4}$.
With this understood, it follows that
$$
k_{z_0^{(1)}}\geq\frac{1}{2m_1}(l-1)(3-1)\geq\frac{1}{2}
\mbox{ and } k_{[z_0^{(1)},z_0^{(2)}]}\geq\frac{1}{4}
\min(l,3l^\prime)\geq\frac{3}{4}, 
$$
which is easily seen a contradiction to the adjunction formula. In any
event, there exists a $J$-holomorphic curve $C$ such that $p_0\in C$
and $C\cdot C_0\leq\frac{7}{12}$.

Next we show that such a $J$-holomorphic curve $C$ is the image of a
member of $\M$. First, note that if $C,C_0$ intersect at
more than one point, then $C\cap C_0=\{p_2,p_3\}$ of order $3$ and $4$,
and one can easily show that this is impossible using the adjunction
formula. (In this case, $g(C)=\frac{1}{2}((\frac{7}{12})^2\cdot
\frac{12}{m}-\frac{7}{12}\cdot\frac{m+1}{m})+1\leq\frac{127}{120}$,
but the right hand side of the adjunction formula is at least
$\frac{1}{2}(1-\frac{1}{3})+\frac{1}{2}(1-\frac{1}{4})
+\frac{1}{2}+\frac{1}{4m}>\frac{29}{24}$.) Second, suppose $C,C_0$
intersect only at the singular point $p_i$ of order $a_i$. Then one
can easily show, as we did earlier, that $C$ is the image of a member
of $\M$ if $C\cdot C_0=\frac{1}{a_i}$. With this understood,
it remains to rule out the case where $C$, $C_0$ intersect at $p_3$ of
order $4$ but $C\cdot C_0=\frac{1}{2}$. To this end, let
$f:\Sigma\rightarrow X$ be a parametrization of $C$. Then there are
two possibilities: either $f^{-1}(p_3)$ consists of two points, or
it contains only one point. In any case, one can easily show that the
contribution of $f^{-1}(p_3)$ to the right hand side of the adjunction
formula is at least $\frac{1}{2}$. With the contribution of at least
$\frac{1}{2}+\frac{1}{4m}$ from $f^{-1}(p_0)$, the right hand side of
the adjunction formula is greater than $1$. But on the left hand side,
$g(C)=\frac{1}{2}((\frac{1}{2})^2\cdot\frac{12}{m}-\frac{1}{2}\cdot
\frac{m+1}{m})+1=\frac{5-m}{4m}+1\leq 1$, which is a contradiction. Hence $C$
is the image of a member of $\M$, and the case where $m\neq 1$
is done.

Finally, consider the case of $m=1$. Let $\{C_i\}$ be the
$J$-holomorphic curves obtained by applying Theorem 2.2 (2) to
$E$. This time $d(E)=14$, so we can specify any $6$ distinct smooth
points $q_1,q_2,\cdots,q_6$, where $q_1\in C_0$ and $q_2,\cdots,q_6\in
X\setminus C_0$, and require that they are contained in $\cup_i C_i$,
and moreover, we require $p_0\in \cup_i C_i$ also. We let $q_2,\cdots,
q_6$ converge to a smooth point $q\neq q_1$ in $C_0$ one by one, and
as we have argued earlier, we obtain at the end a $J$-holomorphic curve
$C$ such that $p_0\in C$ and $C\cdot C_0\leq 1-6C_0\cdot C_0
=\frac{1}{2}$. However, in order to show that such a curve $C$ is the
image of a member of $\M$, we actually need to obtain a sharper
estimate that $C\cdot C_0<\frac{1}{2}$. To this end, for each $k=2,
\cdots,6$, we let $\{C_i^{(k)}\}$ be the non-$C_0$ components in the
limiting cusp-curve as $q_k$ converges to a smooth point in $C_0$, and
let $\alpha_k\equiv\sum_i n_i^{(k)} C_i^{(k)}\cdot C_0$, where $n_i^{(k)}$
is the multiplicity of $C_i^{(k)}$. If the said sharper estimate does not
hold, then we must have $\alpha_2=\frac{10}{12}$, $\alpha_3=\frac{9}{12}$,
$\alpha_4=\frac{8}{12}$, $\alpha_5=\frac{7}{12}$, and $\alpha_6
=\frac{6}{12}=\frac{1}{2}$. We will show that this is impossible.
To see this, note that $\alpha_4=\frac{8}{12}=\frac{2}{3}$. It follows
easily that every component in $\{C_i^{(4)}\}$ must intersect $C_0$ only
at the singular point $p_2$ of order $3$. On the other hand, as $q_5$
converges to a smooth point in $C_0$, there must be a component $\Sigma_\nu$
in the limiting cusp-curve $f:\Sigma^\prime\rightarrow X$ such that
$\mbox{Im }f|_{\Sigma_\nu}=C_0$. The key point here is that the degree of
$f|_{\Sigma_\nu}:\Sigma_\nu\rightarrow C_0$ is at least $2$. Suppose
to the contrary that the degree is $1$. Then $\Sigma_\nu$ must contain
three orbifold points $z_1,z_2,z_3$ such that $f(z_i)=p_i$ for $i=1,2$
and $3$. Moreover, the points $z_1,z_3$, which are sent to the singular
points $p_1,p_3$ of order $2$ and $4$ under the map $f$, must be the images
of some collapsing simple closed loops. This implies that for each of
$p_1$ and $p_3$, there exists a nonconstant component in $\{C_i^{(5)}\}$
which intersects $C_0$ at it. But this is impossible because it would
imply that $\alpha_5\geq\frac{1}{2}+\frac{1}{4}=\frac{9}{12}$, which is a 
contradiction. Hence there exists a $J$-holomorphic curve $C$ such
that $p_0\in C$ and $C\cdot C_0<\frac{1}{2}$. It is easy to show that
$C$ is the image of a member of $\M$. Hence the
case where $m=1$.

(4) $G=\langle Z_{2m},Z_{2m};\widetilde{I},\widetilde{I}\rangle$.
In this case, $|G|>4m^2$ implies that $m=1,7,11,13,17, 19,23$ or $29$.
We shall divide them into three groups: $m\geq 11$, $m=7$, and $m=1$.

First of all, since each element of $G$ has order no greater than $10m$,
we see that for any $J$-holomorphic curve $C$ parametrized by $f$, a point
in $f^{-1}(p_0)$ will contribute at least $\frac{1}{2}+\frac{1}{2m}$ to
the right hand side of the adjunction formula.

Consider first the case where $m\geq 11$. By a similar argument, one can
show that there exists a $J$-holomorphic curve $C$ such that $p_0\in C$
and $C\cdot C_0\leq\frac{19}{30}$. The only part in the proof that is
not so straightforward is to rule out the possibility, in the case where
$m=11$ or $13$, of having a $J$-holomorphic curve $C^\prime$ which obeys
(1) $p_0\in C^\prime$, (2) $C^\prime\cdot C_0=1$, (3) $C^\prime$ is
parametrized by $f:\Sigma\rightarrow X$ such that $f^{-1}(C_0)=\{z_1,z_2\}$
with both $f(z_1)$, $f(z_2)$ being the singular point $p_3$ of order $5$.
First, suppose $m=11$. In this case, the normalized Seifert invariant
at $p_3$ is $(5,1)$. It follows easily that a local representative
$(f_i,\rho_i)$ of $f$ at $z_i$, where $i=1,2$, must be of the form
$\rho_i(\mu_5)=\mu_5^{l_i}$, $f_i(w)=(c_i(w^{l_i^\prime}+\cdots),w^{l_i})$,
where $l_i^\prime\equiv l_i\pmod{5}$ when $c_i\neq 0$, which is the case
unless $l_i=1$. Moreover, $l_1+l_2=5$ because $C^\prime\cdot C_0=1$.
With these data, one can easily show that $k_{z_1}+k_{z_2}\geq\frac{2}{5}$
and $k_{[z_1,z_2]}\geq \frac{4}{5}$. Consequently, the right hand side
of the adjunction formula is at least $\frac{1}{2}(1-\frac{1}{5})+
\frac{1}{2}(1-\frac{1}{5})+\frac{2}{5}+\frac{4}{5}+\frac{1}{2}
=2+\frac{1}{2}$. But the left hand side is $g(C^\prime)=\frac{1}{2}
(\frac{30}{m}-\frac{m+1}{m})+1=\frac{29-m}{2m}+1=\frac{20}{11}$, which is
a contradiction. As for the case where $m=13$, the normalized Seifert
invariant at $p_3$ is $(5,3)$. By a similar argument, one can show that
in this case, the right hand side of the adjunction formula is greater
than $2$, which is also a contradiction.

The next order of business is to show that a $J$-holomorphic curve $C$
with $p_0\in C$ and $C\cdot C_0\leq\frac{19}{30}$ must be the image of
a member of $\M$. First of all, observe that $C,C_0$ can
not intersect at more than one point. This is because if otherwise,
the two points in the intersection must be $p_2,p_3$ of order $3$ and
$5$ because $C\cdot C_0\leq\frac{19}{30}$. But in this case,
the right hand side of the adjunction formula
is greater than $\frac{1}{2}(1-\frac{1}{3})+\frac{1}{2}(1-\frac{1}{5})+
\frac{1}{2}=1+\frac{7}{30}$, while the left hand side is $g(C)=
\frac{1}{2}((\frac{8}{15})^2\cdot\frac{30}{m}-\frac{8}{15}\cdot
\frac{m+1}{m})+1=1+\frac{32}{11}\cdot\frac{1}{30}$, which is a contradiction.
Second, if $C,C_0$ intersect only at the singular point $p_i$ of
order $a_i$ for some $i=1,2$ or $3$ and $C\cdot C_0=\frac{1}{a_i}$,
then $C$ is the image of a member of $\M$ as we argued
earlier. With these understood, it is easy to see that there are only
two other possibilities that we need to rule out: $C\cdot C_0=\frac{2}{5}$
or $C\cdot C_0=\frac{3}{5}$, where in both cases, $C,C_0$ intersect only
at $p_3$. Let $f:\Sigma\rightarrow X$ be a parametrization of $C$.
First, it is fairly easy to rule out the possibility that $f^{-1}(p_3)$
and $f^{-1}(p_0)$ may contain more than one point. Moreover, observe also
that $g_{|\Sigma|}=0$. Now let $z_0=f^{-1}(p_3)$ and $z_\infty=
f^{-1}(p_0)$. Note that in both cases, the order of $z_0$ is $5$.
If $C\cdot C_0=\frac{2}{5}$, then as we have seen earlier, the link of
$p_3$ in $C$ is homotopic in $\s^3/G$ to $2$ times of the singular fiber
of the Seifert fibration at $p_3$, which has order $10m$ in $G$. This
implies, since $g_{|\Sigma|}=0$, that the order of $z_\infty$ is
no greater than $5m$. As a consequence, the right hand side of the
adjunction formula is no less than
$ \frac{1}{2}(1-\frac{1}{5})+\frac{1}{2}(1-\frac{1}{5m})+\frac{1}{10m}
(\frac{120m}{5m}-1)=\frac{11}{5m}+\frac{9}{10}
$.
But on the left hand side, $g(C)=\frac{1}{2}((\frac{2}{5})^2\cdot
\frac{30}{m}-\frac{2}{5}\cdot\frac{m+1}{m})+1=\frac{11}{5m}+\frac{4}{5}$,
which is a contradiction. The case where $C\cdot C_0=\frac{3}{5}$ is 
more involved.
First, let $h\equiv\mu_{2m} I\in Z_{2m}$ and let $x,y$ be the generators
of $\widetilde{I}$ with relations $x^2=y^5=(xy)^3=-1$. Then the homotopy
class of the singular fiber of the Seifert fibration at $p_3$ is represented
by $\gamma^{-1}$, where $\gamma=h^{-t}y^s$ for some positive integers
$s,t$ satisfying $sm-5t=1$. The action of $\gamma$ on $\C^2$ is given,
in suitable coordinates, by $\gamma\cdot (z_1,z_2)=(\mu_{10m}z_1,
\mu_{10m}^k z_2)$ with $k\equiv -sm-5t\pmod{10m}$. With these understood,
observe that the link of $p_3$ in $C$ is homotopic in $\s^3/G$ to $3$
times of the singular fiber of the Seifert fibration at $p_3$, and
consequently, since $g_{|\Sigma|}=0$, there are holomorphic coordinates
$z_1,z_2$ on a local uniformizing system at $p_0$, such that a local
representative of $f$ at $z_\infty$ is given by $(f_\infty,\rho_\infty)$,
with $\rho_\infty(\mu_{10m})$ acting by $\rho_\infty(\mu_{10m})\cdot
(z_1,z_2)=(\mu_{10m}^3 z_1,\mu_{10m}^{3k} z_2)$, and $f_\infty(z)=
(c_1(z^{l_1}+\cdots),c_2(z^{l_2}+\cdots))$, where both $c_1,c_2\neq 0$
because both $l_1,l_2$ can not be $1$, and because of this, we have
relations $l_1\equiv 3\pmod{10m}$ and $l_2\equiv 3k\equiv -3(sm+5t)
\pmod{10m}$. The latter particularly implies that $l_2\geq 3$. With
these understood, we have $k_{z_\infty}\geq\frac{1}{20m}((3-1)(3-1)+
(\frac{120m}{10m}-1)\cdot 3^2)=\frac{103}{20m}$. With this estimate,
it is easy to see that the right hand side of the adjunction formula
is at least $\frac{9}{10}+\frac{51}{10m}$. However, the left hand
side is $g(C)=\frac{1}{2}((\frac{3}{5})^2\cdot\frac{30}{m}-\frac{3}{5}
\cdot\frac{m+1}{m})+1=\frac{51}{10m}+\frac{7}{10}$, which is a contradiction.
This finishes the proof that a $J$-holomorphic curve $C$
with $p_0\in C$ and $C\cdot C_0\leq\frac{19}{30}$ must be the image of
a member of $\M$, and the case where $m\geq 11$ is done.

For the next case where $m=7$, we begin with the following observation.
Let $C_i$, $i=0,1,\cdots,k$, be $J$-holomorphic curves with
multiplicities $n_i$ such that $\sum_{i=0}^k n_i C_i\cdot C_0=1$.
Here $C_i$ with $i=0$ stands for the distinguished $J$-holomorphic curve
$C_0$, and we allow $n_0=0$, which simply means that $C_0$ is not included.
Note that on the one hand, $\sum_{i\neq 0} n_i C_i\cdot C_0=
1-n_0 C_0\cdot C_0=1-\frac{7n_0}{30}$, and on the other hand,
$\sum_{i\neq 0}n_i C_i\cdot C_0=\frac{c_1}{2}+\frac{c_2}{3}+\frac{c_3}{5}$
for some non-negative integers $c_1,c_2,c_3$, where at least one of
them must be $0$ because $\frac{1}{2}+\frac{1}{3}+\frac{1}{5}>1$.
It follows easily that either $n_0=0$, or $n_0=2$ with
$\sum_{i\neq 0}n_i C_i\cdot C_0=\frac{1}{3}+\frac{1}{5}$. With this
understood, note that $d(E)=4$ when $m=7$, so that we can specify any $2$
distinct smooth points $q_1,q_2$, with $q_1\in C_0$ and $q_2\in X\setminus
C_0$, such that $q_1,q_2$ are contained in the $J$-holomorphic
curves $\{C_i\}$ obtained by applying Theorem 2.2 (2) to $E$,
where we note that $p_0\in\cup_i C_i$ also. It follows easily from
the preceding observation that as $q_2$ converges to a smooth
point $q_2^\prime\in C_0$ where $q_2^\prime\neq q_1$, $\{C_i\}$
will have only one component, which intersects $C_0$ transversely at $q_1$.
Let $C$ be the $J$-holomorphic curve and let $f:\Sigma\rightarrow X$ be a
parametrization of $C$. First, the virtual genus of $C$ is
$g(C)=\frac{1}{2}(\frac{30}{7}-\frac{7+1}{7})+1=2+\frac{4}{7}$,
from which it follows that $g_{|\Sigma|}\leq 2$. However, we shall
need a sharper estimate that $g_{|\Sigma|}\leq 1$, which is
obtained as follows. Observe that by the adjunction formula, if
$g_{|\Sigma|}=2$, then $f^{-1}(p_0)$ must contain only one point,
denoted by $z_\infty$, which must have order $10m=70$. Moreover,
let $(f_\infty,\rho_\infty)$ be a local representative of $f$ at
$z_\infty$, then $f_\infty$ must also be embedded. Now suppose the action
of $\rho_\infty(\mu_{70})$ on a local uniformizing system at $p_0$ is
given by $\rho_\infty(\mu_{70})\cdot (z_1,z_2)=(\mu_{70}^{m_1}z_1,
\mu_{70}^{m_2}z_2)$ in some holomorphic coordinates $z_1,z_2$.
Write $f_\infty(z)=(c_1(z^{l_1}+\cdots),c_2(z^{l_2}+\cdots))$ where
$c_1\neq 0$ (resp. $c_2\neq 0$) unless $l_2=1$ (resp. $l_1=1$). Then
for any $i=1,2$, we have $l_i\equiv m_i\pmod{70}$ as long as $c_i\neq 0$.
It follows easily, since $f_\infty$ is embedded, that one of $m_1,m_2$ must
equal $1$. On the other hand, the index formula for the
linearization $D\underline{L}$ at $f$ gives rise to
$$
c_1(TX)\cdot [f(\Sigma)]+2-2g_{|\Sigma|}-\frac{m_1+m_2}{70}=
\frac{7+1}{7}+2-2g_{|\Sigma|}-\frac{m_1+m_2}{70}\in\Z.
$$
If we write $\rho_\infty(\mu_{70})=h^k y^l$, then $\frac{m_1+m_2}{70}
=\frac{k}{7}$ because $\widetilde{I}\subset SU(2)$. Here $h=\mu_{14}
I\in Z_{2m}=Z_{14}$ and $y\in\widetilde{I}$ with eigenvalues $\mu_{10},
\mu_{10}^{-1}$, and without loss of generality, we assume that
$0\leq k\leq 6$. As a consequence, we obtain
$k=1$ and $\rho_\infty(\mu_{70})=hy^l$. It follows easily that
$m_1\equiv 5+7l\pmod{70}$ and $m_2\equiv 5-7l\pmod{70}$, and from
this one can easily check that $m_1\neq 1$, $m_2\neq 1$ for any
$l$. This is a contradiction, hence $g_{|\Sigma|}\leq 1$. With this in
hand, we let $q_2\equiv q_{2,j}$ be a sequence of points converging
to a smooth point $q_2^\prime\neq q_1$ on $C_0$, and denote by $C_j$
the corresponding $J$-holomorphic curves, and by $f_j:\Sigma\rightarrow
X$ a parametrization of $C_j$, which is assumed to be from a fixed
orbifold Riemann surface without loss of generality. As we argued
earlier, $f_j$ will converge to a cusp-curve $f:\Sigma^\prime
\rightarrow X$ such that a component $\Sigma_\nu$ of $\Sigma^\prime$
is mapped to $C_0$ under $f$, over which $f$ is a map of degree
at least $2$. By the observation made at the beginning of this paragraph,
we see that the degree of $f|_{\Sigma_\nu}:\Sigma_\nu\rightarrow X$
is exactly $2$, and moreover, the remaining component(s) in the limiting
cusp-curve must intersect $C_0$ at exactly two singular points, $p_2$
of order $3$ and $p_3$ of order $5$, each contributing $\frac{1}{3}$
and $\frac{1}{5}$ to the intersection product. Now we observe that
since $f|_{\Sigma_\nu}:\Sigma_\nu\rightarrow X$ is of degree $2$, there
must be at least two orbifold points in $\Sigma_\nu$, one is of order
$3$ and the other of order $5$, which are all obtained by collapsing
simple closed loops in $\Sigma$. Furthermore, if $\Sigma_\nu$ contains
exactly $2$ orbifold points, then we must also have $g_{|\Sigma_\nu|}
\neq 0$. Since $g_{|\Sigma|}\leq 1$, it is not hard to see that
$g_{|\Sigma|}$ must equal $1$ and each $f^{-1}_j(p_0)$ consists of
two points, and that there are exactly two non-$C_0$ components, denoted
by $C_1,C_2$, in the limiting cusp-curve $\mbox{Im }f$, such that
$p_0\in C_1\cap C_2$ and $C_1\cdot C_0=\frac{1}{3}$, $C_2\cdot C_0=
\frac{1}{5}$ (or the other way). It follows easily that both $C_1,C_2$
are the image of a member of $\M$. Hence the case where
$m=7$.

Finally, we consider the case where $m=1$. Since $d(E)=32$, we can
specify any $15$ distinct smooth points $q_1,q_2,\cdots,q_{15}$,
where $q_1\in C_0$ and $q_2,\cdots,q_{15}\in X\setminus C_0$, such that
the $J$-holomorphic curves $\{C_i\}$ obtained from Theorem 2.2 (2)
contain these points as well as the singular point $p_0$. We then
let $q_k$, $k=2,\cdots,15$, converge one by one to a smooth point in
$C_0$ which is different than $q_1$. If we denote by $\alpha_k$ the
intersection product with $C_0$ of the non-$C_0$ components
(counted with multiplicity) in the limiting cusp-curve at each stage,
then we have, as in the previous cases, that $\alpha_2\leq\frac{28}{30}$
and $\alpha_{k}-\alpha_{k+1}\geq\frac{1}{30}$ for $k=2,\cdots,14$. The
key observation for the present case is that each
$\alpha_k=\frac{c_1^{(k)}}{2}+\frac{c_2^{(k)}}{3}+\frac{c_3^{(k)}}{5}$
for some integers $c_i^{(k)}\geq 0$, $i=1,2,3$,
where at least one of $c_i^{(k)}$ is zero. With this understood, note,
for instance, that $\frac{23}{30}$ can not be written in the above form,
and therefore it can not be realized as $\alpha_k$ for any $k$. In fact,
a simple inspection like this shows that the following is the only
possibility for the values of $\alpha_k$:
$$
\begin{array}{l}
\alpha_2=\frac{28}{30}, \cdots, \alpha_6=\frac{24}{30},\alpha_7
=\frac{22}{30},\cdots,\alpha_9=\frac{20}{30},\alpha_{10}=\frac{18}{30},\\
\alpha_{11}=\frac{16}{30},
\alpha_{12}=\frac{15}{30},\alpha_{13}=\frac{12}{30},\alpha_{14}=
\frac{10}{30} \mbox{ and } \alpha_{15}=\frac{6}{30}=\frac{1}{5}.
\end{array}
$$
(In fact, using the adjunction formula, one can explicitly recover
the process of degeneration of the $J$-holomorphic curves, i.e.
understanding how at each stage a component carrying the correct amount
of homology splits off during the convergence. But these details
are not needed here for the proof, so we leave them to the reader
as an exercise.)
In particular, we obtain at the last stage a $J$-holomorphic curve $C$
such that $p_0\in C$ and $C\cdot C_0=\frac{1}{5}$. It follows easily
that $C$ is the image of a member of $\M$. Hence the case
where $m=1$.

\vspace{2mm}

The proof of Lemma 3.9 is thus completed.

\hfill $\Box$

\sectioni{Proof of Taubes' theorems for $4$-orbifolds}

For the assertions in Theorem 2.2 (1), observe that Taubes' proof
(cf. e.g. \cite{HT}) works in the orbifold setting without changing a word.

The rest of this section is occupied by a proof of Theorem 2.2 (2). 
We shall follow the proof of Taubes in \cite{T3}, indicating along 
the way where modifications to Taubes' proof are necessary in the orbifold 
setting, and how to implement them.

\vspace{2mm}

{\em Basic estimates.} Section 2 in Taubes \cite{T3} is concerned with
the following estimates:
\begin{itemize}
\item $|\alpha|\leq 1+zr^{-1}$
\item $|\beta|^2\leq zr^{-1}((1-|\alpha|^2)+r^{-2})$
\item $|P_{\pm}F_a|\leq (4\sqrt{2})^{-1}r(1+zr^{-1/2})(1-|\alpha|^2)+z$
\item $|\nabla_a\alpha|^2+r|\nabla_A^\prime\beta|^2\leq zr(1-|\alpha|^2)+z$
\end{itemize}
Here $z>0$ is a constant solely determined by $c_1(E)$ and the Riemannian
metric, and $r$ is sufficiently large. The principal tool
for obtaining these estimates is to apply the maximum principle to the
various differential inequalities derived from the Seiberg-Witten equations.
Another important ingredient is the total energy bound in Lemma 2.6
of \cite{T3}:
$$
|\frac{r}{4}\int_X|1-|\alpha|^2|-2\pi [\omega]\cdot c_1(E)|\leq zr^{-1}
$$
These are all valid in the orbifold setting. However, in the estimate for
$|P_{-}F_a|$ (specifically $(2.35)$ in the proof of Lemma 2.7 in \cite{T3}),
Green's function for the Laplacian $d^\ast d$ is also involved.
Here additional care is needed in the orbifold case because even on a
compact, closed Riemannian orbifold, the injectivity radius at each point
is not uniformly bounded from below by a positive constant due to the
presence of orbifold points.

Green's function for the Laplacian on orbifolds is discussed in Appendix B.
Given that, let's recall that the part in the proof of Lemma 2.7 in Taubes
\cite{T3} which involves Green's function is to derive the following estimate
(cf. $(2.35)$ in \cite{T3}) for the function $q_1^\prime$:
$$
q_1^\prime\leq \frac{z\cdot R\cdot \sup(|P_{-}F_a|)}{r^{1/2}}
$$
where $q_1^\prime$ satisfies $\frac{1}{2}d^\ast dq^\prime_1+\frac{r}{4}
|\alpha|^2 q_1^\prime=R\cdot \sup(|P_{-}F_a|)\cdot |1-|\alpha|^2|$. In the
present case, we apply Theorem 1 in Appendix B to $q_1^\prime$,
$$
q_1^\prime(x)=Vol(X)^{-1}\int_X q_1^\prime+\int_X G(x,\cdot)\Delta q_1^\prime.
$$
Now observe that in the first term, $\int_X q_1^\prime$ is bounded by
$$
\int_X|1-|\alpha|^2|q_1^\prime+\int_X|\alpha|^2q_1^\prime\leq
z(\frac{\sup(q_1^\prime)}{r}+\frac{R\cdot\sup(|P_{-}F_a|)}{r^2})
$$
because $|\frac{r}{4}\int_X|1-|\alpha|^2|-2\pi [\omega]
\cdot c_1(E)|\leq zr^{-1}$ and $\frac{r}{4}\int_X|\alpha|^2q_1^\prime=
\int_X R\cdot \sup(|P_{-}F_a|)\cdot |1-|\alpha|^2|$. As for the second
term, suppose $q_1^\prime(x_0)=\sup(q_1^\prime)$ for some $x_0\in X$, and
recall Theorem 1 (3) in Appendix B that one may write
$G(x_0,y)=G_0(x_0,y)+G_1(x_0,y)$.
Thus $\int_X G(x_0,\cdot) \Delta q_1^\prime$ is bounded by
\begin{eqnarray*}
& &\int_X G_0(x_0,\cdot)(2R\cdot\sup(|P_{-}F_a|)\cdot |1-|\alpha|^2|) \\
& + & \int_X G_1(x_0,\cdot)(2R\cdot\sup(|P_{-}F_a|)\cdot |1-|\alpha|^2|).
\end{eqnarray*}
The last term above is bounded by $\frac{z\cdot
R\cdot\sup(|P_{-}F_a|)}{r}$,
and for the first term, recall that there is a uniformizing system
$(\widehat{U},G,\pi)$ such that $G_0(x_0,y)$ is supported in $\pi(\widehat{U})$
and $G_0\circ\pi$ equals $\sum_{h\in G}\widehat{G_0}(h\cdot\hat{x_0},\hat{y})$
for some $\hat{x_0}\in\pi^{-1}(x_0)$ with $\widehat{G_0}(\hat{x_0},\hat{y})$
satisfying $|\widehat{G_0}(\hat{x_0},\hat{y})|\leq
\frac{z}{\hat{d}(\hat{x_0},\hat{y})^{2}}$. Moreover, $\widehat{U}$
contains a closed ball of radius $\delta_0>0$ centered at $\hat{x_0}$.
Now when $r\geq \delta_0^{-4}$, the first term $\int_X G_0(x_0,\cdot)(2R
\cdot\sup(|P_{-}F_a|)\cdot |1-|\alpha|^2|)$, which equals
$$
\int_{\widehat{U}}\widehat{G_0}(\hat{x_0},\cdot)((2R
\cdot\sup(|P_{-}F_a|)\cdot |1-|\alpha|^2|)\circ\pi),
$$
is bounded by
$$
z\cdot R\cdot \sup(|P_{-}F_a|)\cdot r^{-1/2}+ z\cdot R\cdot
\sup(|P_{-}F_a|)\cdot r^{1/2}\int_X |1-|\alpha|^2|
$$
by writing the integration over $\widehat{U}$ as the sum of integration
over the closed ball $B_{\hat{x_0}}(r^{-1/4})\subset \widehat{U}$ of
radius $r^{-1/4}\leq \delta_0$ centered at $\hat{x_0}$ and its complement
in $\widehat{U}$ as in \cite{T3}. It is easily seen that the estimate for
$q_1^\prime$ follows immediately for $r$ sufficiently large.

\vspace{2mm}

{\em Monotonicity formula and refined estimate for $|P_{-}F_a|$.}
Recall that the monotonicity formula in Section 3 of Taubes \cite{T3}
is for the purpose of estimating the growth rate of the local energy
$\frac{r}{4}\int_B|1-|\alpha|^2|$, where $B$ is a geodesic ball of radius
$s$ centered at a given point, against the radius of the ball $s$. In this
part of the argument, the radius $s$ is required to satisfy an inequality
$\frac{1}{2r^{1/2}}\leq s\leq \frac{1}{z}$. Thus again, because the
injectivity radius has no positive uniform lower bound on the orbifold
$X$, a reformulation for the definition of local energy is needed.

More precisely, we shall fix the set $\U$ of finitely many uniformizing
systems and the constant $\delta_0>0$ as described in Theorem 1 (3) in
Appendix B. Given that, for any $p\in X$, we choose a uniformizing system
$(\widehat{U},G,\pi)\in\U$ for $p$ as described therein, and define the
local energy at $p$ to be
$$
\E(p,s)=\frac{r}{4}\int_B|1-|\alpha|^2|,
$$
where $B$ is the geodesic ball of radius $s\leq \delta_0$ in $\widehat{U}$
centered at some $\hat{p}\in\pi^{-1}(p)$, and by abusing the notation,
the function $|1-|\alpha|^2|\circ\pi$ on $\widehat{U}$ is still denoted
by $|1-|\alpha|^2|$. It is easily seen that $\E(p,s)$ is well-defined,
i.e., $\E(p,s)$ is independent of the choice of $(\widehat{U},G,\pi)\in\U$
and $\hat{p}\in\pi^{-1}(p)$.

With the preceding understood, the relevant argument in Taubes \cite{T3}
can be quoted to establish the corresponding estimates in the present case:
\begin{itemize}
\item $\E(p,s)\leq zs^2$ for all $p\in X$, and
\item $\E(p,s)\geq\frac{1}{z+1} s^2$ when $|\alpha(p)|<1/2$,
\end{itemize}
where $z>0$ is a constant depending only on $c_1(E)$ and the Riemannian metric,
$r$ is sufficiently large, and $\frac{1}{2r^{1/2}}\leq s\leq \frac{1}{z}$.
(cf. Prop. 3.1 in \cite{T3}.)

Now we discuss the refined estimate for $|P_{-}F_a|$ (cf. Prop. 3.4 in
\cite{T3}). Here the argument involves Green's function as well as a
ball covering procedure using geodesic balls. Hence Taubes' original
proof in \cite{T3} needs to be modified in the present case.

Recall that the key to the refinement is Lemma 3.5 in Taubes \cite{T3}
where a smooth function $u$ on $X$ is constructed which obeys
\begin{itemize}
\item [{(1)}] $|u|\leq z$.
\item [{(2)}] $\frac{1}{2}d^\ast du\geq r$ where $|\alpha|<1/2$.
\item [{(3)}] $|d^\ast du|\leq z\cdot r$.
\end{itemize}
Here $z>0$ is a constant depending only on $c_1(E)$ and the Riemannian
metric. The strategy for the present case is to follow the proof of
Lemma 3.5 in Taubes \cite{T3} to construct, for each uniformizing system
$(\widehat{U},G,\pi)\in\U$, a function $u_{\widehat{U}}$ on $X$, and
define $u\equiv \sum u_{\widehat{U}}$.

To be more concrete, let $(\widehat{U},G,\pi)$ be any element in $\U$.
Recall that (cf. Theorem 1 in Appendix B) $\widehat{U}$ is a geodesic ball
of radius $\delta(p_i)$, the injectivity radius at $p_i$ for some
$p_i\in X$, and $G=G_{p_i}$, the isotropy group at $p_i$. Moreover,
the open subset $\widehat{U}^\prime\subset\widehat{U}$ is the concentric
ball of radius $\delta_i=N^{-1}\delta(p_i)$, and if we denote by
$\widehat{U}_0$ the concentric ball of radius $N^{-1}\delta_i
=N^{-1}\cdot\mbox{radius}(\widehat{U}^\prime)$, then the set
$\{\pi(\widehat{U}_0)\}$ is an open cover of $X$. Here $N$ is a fixed
integer no less than $12=3\cdot\dim X$.

With the preceding understood, let $V$ be the region in $\widehat{U}_0$
where $|\alpha|<1/2$. Then Lemma 3.6 in Taubes \cite{T3} is valid here.
To be more precise, there is a set $\{B_i\}$ of geodesic balls in
$\widehat{U}$ of radius $r^{-1/2}\leq\delta_0$ having the following
properties: (1) each $B_i$ is centered at a point $\hat{p}_i\in V$,
the region in $\widehat{U}_0$ where $|\alpha|<1/2$, (2) $\{B_i\}$ covers
$V$, (3) the number of balls, $\#\{B_i\}$, is bounded by $z\cdot r$ as
$r$ grows, (4) the concentric balls of only half radius (i.e.
$\frac{1}{2}r^{-1/2}$) are disjoint, and furthermore in the present
case, (5) the set of centers $\{\hat{p}_i\}$ of the balls is invariant
under the action of $G$.

Now observe that Lemma 3.7 in Taubes \cite{T3} is valid for the set
of balls $\{B_i\}$. Thus there exists a set of concentric balls
$\{\widetilde{B_i}\}$ of radius $z\cdot r^{-1/2}$ for some constant $z>1$
such that $\mbox{Volume}((\widehat{U}\setminus V^\prime)\cap\widetilde{B_i})
\geq\mbox{Volume}(B_i)$, where $V^\prime$ is the region in $\widehat{U}$
where $|\alpha|<3/4$. Here we choose $r$ sufficiently large so that
each $\widetilde{B_i}$ is contained in the ball of radius $\delta=\delta_0+
\mbox{radius}(\widehat{U}_0)$ which has the same center of $\widehat{U}_0$.

As in Taubes \cite{T3}, we let $s_i,\tilde{s}_i$ be the characteristic
functions of $B_i$ and $(\widehat{U}\setminus V^\prime)\cap\widetilde{B_i}$.
Then as in \cite{T3}, there is a $\kappa_i$, with bound $z^{-1}<\kappa_i
<z$, such that
$$
\int_{\widehat{U}}(s_i-\kappa_i\tilde{s}_i)=0.
$$
Note that the function $\sum_i (s_i-\kappa_i\tilde{s}_i)$ on $\widehat{U}$
is invariant under the action of $G$ and is compactly supported in
the ball of radius $\delta=\delta_0+\mbox{radius}(\widehat{U}_0)$ which has
the same center of $\widehat{U}_0$. Thus $\sum_i (s_i-\kappa_i\tilde{s}_i)$
descends to a function $f_{\widehat{U}}$ on $X$ by defining $f_{\widehat{U}}
\equiv 0$ outside $\pi(\widehat{U})$. With the preceding understood,
the function $u_{\widehat{U}}$ is the unique solution to
$$
\frac{1}{2}d^\ast du_{\widehat{U}}=r\cdot f_{\widehat{U}}
\mbox{ and } \int_X u_{\widehat{U}}=0.
$$
(By suitably smoothing $f_{\widehat{U}}$, as indicated in \cite{T3},
one may arrange to have $u_{\widehat{U}}$ smooth.)

The following properties of $u_{\widehat{U}}$ are straightforward as
in \cite{T3}:
\begin{itemize}
\item $|d^\ast du_{\widehat{U}}|\leq z\cdot r$.
\item $\frac{1}{2}d^\ast du_{\widehat{U}}\geq 0$ where $|\alpha|<1/2$,
and $\frac{1}{2}d^\ast du_{\widehat{U}}\geq r$ in $\pi(V)$.
\end{itemize}

Thus to furnish Lemma 3.5 in Taubes \cite{T3} with $u\equiv
\sum u_{\widehat{U}}$, it suffices to show that
$$
|u_{\widehat{U}}|\leq z
$$
for a constant $z>0$ which is independent of $r$.

To this end, we invoke Theorem 1 in Appendix B to obtain
$$
u_{\widehat{U}}(p)=2r\cdot\int_X G(p,\cdot) f_{\widehat{U}}=
2r\cdot \int_X G_0(p,\cdot) f_{\widehat{U}}
+2r\cdot\int_X G_1(p,\cdot) f_{\widehat{U}}.
$$
Note that $G_1(p,q)$ is $C^1$ on $X\times X$, so that
$$
2r\cdot\int_X G_1(p,\cdot) f_{\widehat{U}}\leq z_1\cdot\frac{2r}{|G|}
\sum_i\mbox{Volume}(\widetilde{B_i})\leq
z_1\cdot\frac{2r}{|G|}\cdot\#\{B_i\}\cdot\frac{z_2}{r^2}\leq z_3,
$$
which is an $r$-independent constant. As for
$2r\cdot \int_X G_0(p,\cdot) f_{\widehat{U}}$, note that
$\int_X G_0(p,\cdot) f_{\widehat{U}}=0$ if $p\in X\setminus
\pi(\widehat{U}^\prime)$ because $\{q\mid G_0(p,q)\neq 0\}
\subseteq \{q\mid d(p,q)\leq (4+1)\delta_0\}$ (cf. Theorem 1 in
Appendix B). Hence by fixing a $\hat{p}\in\pi^{-1}(p)$ for any given
$p\in\pi(\widehat{U}^\prime)$, we have
$$
2r\cdot\int_X G_0(p,\cdot) f_{\widehat{U}}=\sum_i 2r\cdot\int_{\widehat{U}}
\widehat{G_0}(\hat{p},\cdot)\cdot(s_i-\kappa_i\tilde{s}_i).
$$
If we set $u_i(\hat{p})=2r\cdot\int_{\widehat{U}}
\widehat{G_0}(\hat{p},\cdot)\cdot(s_i-\kappa_i\tilde{s}_i)$, then as
in Taubes \cite{T3}, $u_i$ satisfies
$$
\begin{array}{c}
|u_i(\hat{p})|\leq z \mbox{ when } \hat{d}(\hat{p},\hat{p}_i)
\leq \frac{z}{r^{1/2}}, \mbox{ and }\\
|u_i(\hat{p})|\leq \frac{z}{r^{3/2}\hat{d}(\hat{p},\hat{p}_i)^3}
\mbox{ when } \hat{d}(\hat{p},\hat{p}_i)>\frac{z}{r^{1/2}}.
\end{array}
$$
Recall that, here, $\hat{p}_i\in V$ is the center of the ball $B_i$.

To complete the proof, we observe that Lemma 3.8 in Taubes \cite{T3}
is valid here, that is, for any $\hat{p}\in\widehat{U}^\prime$,
the number $N(n)$ of balls in $\{B_i\}$ whose center $\hat{p}_i$ satisfies
$\hat{d}(\hat{p},\hat{p}_i)\leq n\cdot r^{-1/2}$ obeys $N(n)\leq z\cdot n^2$.
(Here $n$ is any positive integer.) If we let $\Omega(n)$ be the set of
indices $i$ for the balls $B_i$ whose center $\hat{p}_i$ obeys
$(n-1)\cdot r^{-1/2}<\hat{d}(\hat{p},\hat{p}_i)\leq n\cdot r^{-1/2}$,
then as in \cite{T3},
\begin{eqnarray*}
|2r\cdot\int_X G_0(p,\cdot) f_{\widehat{U}}| & \leq & \sum_i|u_i(\hat{p})|
=\sum_{n\geq 1}\sum_{i\in\Omega(n)}|u_i(\hat{p})|\\
& \leq &
z_1+\sum_{n\geq 1}z_2\cdot\frac{N(n)-N(n-1)}{n^3}\leq z_3
\end{eqnarray*}
for a constant $z_3>0$ which is independent of $r$.

\vspace{2mm}

{\em The local structure of $\alpha^{-1}(0)$ and exponential decay
estimates.}
The discussion in Section 4 of Taubes \cite{T3} extends to
the present case almost word by word, except for the exponential decay
estimates
$$
|q(x)|\leq z\cdot r\cdot\exp(-\frac{1}{z}r^{1/2}d(x,\alpha^{-1}(0)))
$$
for $q\in\{r(1-|\alpha|^2),r^{3/2}\beta,F_a,r^{1/2}\nabla_a\alpha,
r\nabla_A^\prime\beta\}$, where $d$ is the distance function.

The part that needs modification is the construction of a comparison
function $h$ (cf. $(4.19)$ in \cite{T3}) which obeys
\begin{itemize}
\item $\frac{1}{2}d^\ast dh+\frac{r}{32}h\geq 0$ where
$d(\cdot,\alpha^{-1}(0))\geq zr^{-1/2}$.
\item $h\geq r$ where $d(\cdot,\alpha^{-1}(0))=zr^{-1/2}$.
\item $h\leq z_1\cdot r\cdot\exp (-\frac{1}{z_1}r^{1/2}
d(\cdot,\alpha^{-1}(0)))$ where $d(\cdot,\alpha^{-1}(0))\geq zr^{-1/2}$.
\end{itemize}
Here $z,z_1>1$ are $r$-independent constants.

Modification is needed here at least for one reason: the construction
of $h(x)$ involves a ball covering argument by geodesic balls of radius
of size $r^{-1/2}$, along with the local energy growth rate estimates,
i.e., Prop. 3.1 in \cite{T3}. On the other hand, observe that the
construction of comparison function in \cite{T3} does not require the
compactness of the underlying manifold. The compactness enters only when
the maximum principle is applied. Hence Taubes' argument in
\cite{T3} should in principle work here also.

More concretely, we shall construct for each uniformizing system
$(\widehat{U},G,\pi)\in\U$ (cf. Theorem 1 in Appendix B) a smooth
function $h_{\widehat{U}}>0$ on $X$ which obeys
\begin{itemize}
\item $\frac{1}{2}d^\ast dh_{\widehat{U}}+\frac{r}{32}h_{\widehat{U}}
\geq 0$ where $d(\cdot,\alpha^{-1}(0)\cap\pi(\widehat{U}^\prime))
\geq zr^{-1/2}$.
\item $h_{\widehat{U}}\geq r$ where $d(\cdot,\alpha^{-1}(0)\cap
\pi(\widehat{U}^\prime))=zr^{-1/2}$.
\item $h_{\widehat{U}}\leq z_1\cdot r\cdot\exp (-\frac{1}{z_1}r^{1/2}
d(\cdot,\alpha^{-1}(0)\cap\pi(\widehat{U}^\prime)))$ where
$d(\cdot,\alpha^{-1}(0)\cap\pi(\widehat{U}^\prime))\geq zr^{-1/2}$.
\end{itemize}
Here $z,z_1>1$ are $r$-independent constants. Accepting this, we may
take $h\equiv \sum h_{\widehat{U}}$ for the comparison function.

To define $h_{\widehat{U}}$, use Lemma 3.6 in \cite{T3} to find a
maximal set $\{\hat{p}_i\}\subset\alpha^{-1}(0)\cap\widehat{U}^\prime$
such that (1) the geodesic balls with centers $\{\hat{p}_i\}$ and radius
$r^{-1/2}$ are disjoint, (2) the set $\{\hat{p}_i\}$ is invariant under
the action of $G$. Then set $\hat{h}\equiv \sum_i H_i$ where
$$
H_i(\hat{p})=\frac{\rho(\hat{d}(\hat{p},\hat{p}_i))}
{\hat{d}(\hat{p},\hat{p}_i)^2}\exp(-\frac{1}{c}r^{1/2}\cdot
\hat{d}(\hat{p},\hat{p}_i))+c\cdot\exp(-\frac{\delta_0}{2c}r^{1/2}).
$$
Here $\hat{d}$ is the distance function on $\widehat{U}$, $\rho(t)$ is
a fixed cut-off function which equals $1$ for
$t\leq\frac{\delta_0}{2}$ and equals zero for $t\geq\delta_0$.
Moreover, $r$ is sufficiently large, and $c>1$ is a fixed, sufficiently
large, $r$-independent constant. (cf. $(4.17)$ in \cite{T3}.) Note
that $\hat{h}\equiv \sum_i H_i$ is smooth, positive, and invariant under
the action of $G$, and is constant outside a compact subset in $\widehat{U}$.
Hence $\hat{h}$ descends to a smooth, positive function on $X$, which is
defined to be $h_{\widehat{U}}$. The claimed properties of $h_{\widehat{U}}$
follow essentially as in Taubes \cite{T3}. (cf. Lemma 4.6 and $(4.18)$
in \cite{T3}.)

\vspace{2mm}

{\em Convergence to a current} (Section 5 of Taubes \cite{T3}). First of
all, note that generalization of the basic theory of currents on smooth
manifolds (cf. e.g. \cite{dR}) to the orbifold setting is straightforward.
In particular, note that a differential form or a differentiable chain
in an orbifold (as introduced in \cite{C2}) naturally defines a current.

Having said this, for any given sequence of solutions
$(a_n,\alpha_n,\beta_n)$ to the Seiberg-Witten equations with
the values of the parameter $r$ unbounded, we define as in Taubes
\cite{T3} a sequence of currents $\F_n$ by
$$
\F_n(\eta)=\frac{\sqrt{-1}}{2\pi}\int_X F_{a_n}\wedge\eta,\;
\forall \eta\in\Omega^2(X).
$$
As in \cite{T3}, the mass norm of $\{\F_n\}$ is uniformly bounded, thus
there is a subsequence, still denoted by $\{\F_n\}$ for simplicity,
which weakly converges to a current $\F$, namely,
$$
\F(\eta)=\lim_{n\rightarrow\infty}\F_n(\eta),\;\forall\eta\in\Omega^2(X).
$$
The current $\F$ is closed, and is Poincar\'{e} dual to $c_1(E)$
in the sense that
$$
\F(\eta)=c_1(E)\cdot [\eta]
$$
for all closed $2$-forms $\eta$.

As for the support of $\F$, which, by definition, is the intersection of
all the closed subsets of $X$ such that the evaluation of $\F$ on any
$2$-form supported in the complement of the closed subset is zero, we
proceed as follows. We fix the set $\U$ of finitely many uniformizing
systems in Theorem 1 of Appendix B, and for each $(\widehat{U},G,\pi)\in
\U$, we run Taubes' argument on $\widehat{U}$. More concretely, for each
integer $N\geq 1$ and each index $n$ with $r_n>z^2\cdot (256)^N$, we
find a maximal set $\Lambda^\prime_n(N)_{\widehat{U}}$ of disjoint
geodesic balls in $\widehat{U}$ with centers in
$\alpha_n^{-1}(0)\cap\mbox{closure}(\widehat{U}^\prime)$
and radius $16^{-N}$ such that the centers of the balls in
$\Lambda^\prime_n(N)_{\widehat{U}}$ are invariant
under the action of $G$, and for any two uniformizing systems
$(\widehat{U}_i,G_i,\pi_i)\in\U$, $i=1,2$, the centers of the balls in
$\Lambda^\prime_n(N)_{\widehat{U}_i}$ which are in the domain or range
of a transition map between the two uniformizing systems are invariant
under the transition map. Then proceeding as in Taubes \cite{T3},
we find a nested set $\{U(N)_{\widehat{U}}\}_{N\geq 1}$ for each
$(\widehat{U},G,\pi)\in\U$, which satisfies
$$
\hat{d}(U(N+1)_{\widehat{U}},\widehat{U}\setminus U(N)_{\widehat{U}})
\geq \frac{3}{2}16^{-N}.
$$
We define $C_{\widehat{U}}\equiv \bigcap_N U(N)_{\widehat{U}}$, and
define $C$ to be the set of orbits of $\bigsqcup C_{\widehat{U}}\subset
\bigsqcup\widehat{U}$ in $X$. It is clear as in \cite{T3} that the
support of $\F$ is contained in $C$, and $\F$ is of type $1-1$.

As for the Hausdorff measure of $C$, first of all, we say that a subset
of $X$ has a finite $m$-dimensional Hausdorff measure if for every
uniformizing system of $X$, the inverse image of the subset in that
uniformizing system has a finite $m$-dimensional Hausdorff measure.
Equivalently, a subset of $X$ has a finite $m$-dimensional Hausdorff
measure if the inverse image of the subset in
$\widehat{U}^\prime$ for each $(\widehat{U},G,\pi)\in\U$
has a finite $m$-dimensional Hausdorff measure. Having said
this, the subset $C$ has a finite $2$-dimensional Hausdorff measure
because each $C_{\widehat{U}}$ does, as argued in Taubes \cite{T3}.

Finally, the local intersection number. We simply apply the relevant
definition and discussion in Taubes \cite{T3} to $C_{\widehat{U}}$
in $\widehat{U}$ for each $(\widehat{U},G,\pi)\in\U$.

\vspace{2mm}

{\em Representing $\F$ by $J$-holomorphic curves.}
Section 6 of Taubes \cite{T3} deals with the regularity of the subset
$C$ in the manifold case, where the main conclusions are: (1) each regular
point in $C$ has a neighborhood which is an embedded, $J$-holomorphic disc
(cf. Lemma 6.11 in \cite{T3}), (2) the singular points in $C$ are
isolated (cf. Lemmas 6.17, 6.18 in \cite{T3}), and their complement in $C$
is diffeomorphic to a disjoint union of $[1,\infty)\times S^1$ when
restricted in a small neighborhood of each singular point (courtesy of
Lemma 6.3 in \cite{T3}). The arguments for these results
are local in nature, hence applicable to $C_{\widehat{U}}\subset
\widehat{U}$ for each $(\widehat{U},G,\pi)\in\U$.

With the preceding understood, particularly, that the subset
$C_{\widehat{U}}\cap\widehat{U}^\prime$ has the said regularity
properties for each $(\widehat{U},G,\pi)\in\U$, we now analyze
the subset $C$ of $X$, which is the set of orbits of $\bigsqcup
C_{\widehat{U}}\subset \bigsqcup\widehat{U}$ in $X$. To this end,
note that the isotropy subgroup $G_{\hat{p}}$ of a point
$\hat{p}\in\widehat{U}$ falls into two types: type A if $G_{\hat{p}}$
fixes a complex line through $\hat{p}$ (which is in fact a
$J$-holomorphic submanifold), or type B if $G_{\hat{p}}$ only fixes
$\hat{p}$ itself. By the unique continuation property of $J$-holomorphic
curves (cf. e.g. \cite{McDS}), it follows easily that there is a subset
$C_{\widehat{U},s}\subset C_{\widehat{U}}\cap\widehat{U}^\prime$
of isolated points, such that the complement of $C_{\widehat{U},s}$,
denoted by $C_{\widehat{U},r}$, consists of regular points and is modeled
on a disjoint union of $[1,\infty)\times S^1$ in a small neighborhood of
each point in $C_{\widehat{U},s}$, and furthermore, $C_{\widehat{U},r}$
is the disjoint union of $C_{\widehat{U},r}^{(1)}$ and
$C_{\widehat{U},r}^{(2)}$, where $C_{\widehat{U},r}^{(1)}$ consists
of points of trivial isotropy subgroup and $C_{\widehat{U},r}^{(2)}$
consists of points of type A isotropy subgroups. In particular,
$C_{\widehat{U},r}$ is a $J$-holomorphic submanifold in $\widehat{U}$.
It is easy to see that the quotient space $C_{\widehat{U},r}/G\subset X$
has the structure of an open Riemann surface with a set of isolated points
removed, and since $\{\pi(\widehat{U}^\prime)\}$ is an open cover of $X$,
there is a subset $C_0\subset C$ which has the structure of a closed
Riemann surface with a set of isolated, hence finitely many (since $C$
is compact) points removed. The restriction of the inclusion
$C\hookrightarrow X$ to each component $C_{0,i}$ of $C_0$ extends to
a continuous map $f_i:\Sigma_i\rightarrow X$, none of which is multiply
covered, where $\Sigma_i$ is the closed Riemann surface obtained by
closing up $C_{0,i}$. We define $C_i\equiv f_i(\Sigma_i)$. Clearly
$C=\cup_i C_i$. Note that $C_0$ is the disjoint union of $C_0^{(1)}$ and
$C_0^{(2)}$, where the former is covered by $\{C_{\widehat{U},r}^{(1)}\}$
and the latter by $\{C_{\widehat{U},r}^{(2)}\}$. The set $\{C_i\}$ is
correspondingly a disjoint union of two subsets, $\{C_i^{(1)}\}$ and
$\{C_i^{(2)}\}$. We will show next that $\{C_i^{(1)}\}$ are type
I $J$-holomorphic curves and $\{C_i^{(2)}\}$ are type II $J$-holomorphic
curves (in the sense of \cite{C2}).
Moreover, there are positive integers $n_i$ such that
$$
c_1(E)=\sum_i n_i\cdot PD(C_i).
$$
(Note that for a type II $J$-holomorphic curve, the Poincar\'{e} dual
$PD(C)$ differs from the usual one by a factor, see \cite{C2} for
details.)

To this end, consider the \'{e}tale topological groupoid $\Gamma$ that
defines the orbifold
structure on $X$ whose space of units is $\bigsqcup \widehat{U}^\prime$.
There is a canonical orbispace structure on each $C_{0,i}$,
making it into a suborbispace $f_i^\prime:C_{0,i}\hookrightarrow X$, which
is obtained by restricting $\Gamma$ to $\bigsqcup C_{\widehat{U},i}$, where
$C_{\widehat{U},i}$ is the inverse image of $C_{0,i}$ under
$\pi:\widehat{U}^\prime\rightarrow\widehat{U}^\prime/G$ (cf. \cite{C1}).
Because $C_{\widehat{U},r}$ (the inverse image of $C_0$ in
$\widehat{U}^\prime$) is modeled on a disjoint union of
$[1,\infty)\times S^1$ in a small neighborhood
of each point in $C_{\widehat{U},s}$ (the inverse image of $C\setminus
C_0$ in $\widehat{U}^\prime$), it is easily seen that the orbispace structure
on $C_{0,i}$ extends uniquely to define an orbifold structure on $\Sigma_i$,
making it into an orbifold Riemann surface. Moreover, the map
$f_i^\prime:C_{0,i}\hookrightarrow X$ extends uniquely to a map
$\hat{f}_i:\Sigma_i\rightarrow X$ between orbifolds in the sense of
\cite{C1}, which is $J$-holomorphic and defines $C_i$ as a $J$-holomorphic
curve in $X$ in the sense of \cite{C2}. Clearly $\{C_i^{(1)}\}$ are of type I
and $\{C_i^{(2)}\}$ are of type II according to the definitions in \cite{C2}.

The positive integers $n_i$ are assigned to $C_i$ as follows. At each point
$\hat{p}\in C_{\widehat{U},i}$, there is an embedded $J$-holomorphic disc
$D$ in $\widehat{U}$ which intersects $C_{\widehat{U},i}$ transversely at
$\hat{p}$. It is shown in \cite{T3} that $\lim_{n\rightarrow\infty}
\frac{\sqrt{-1}}{2\pi}\int_D F_{a_n}$,  denoted by $n(\hat{p})$, 
exists and is a positive integer. Moreover, $n(\hat{p})$ is locally constant,
hence it depends on $C_i$ only. We define $n_i\equiv n(\hat{p})$, $\forall
\hat{p}\in C_{\widehat{U},i}$.
(cf. Prop 5.6 and the discussion before Lemma 6.9 in \cite{T3}.)

With $n_i$ so defined, now for any $2$-form $\eta$ on $X$, we write
$\eta=\sum\eta_{\widehat{U}}$ by a partition of unity subordinate to
the cover $\{\pi(\widehat{U})\}$, and observe that
\begin{eqnarray*}
\F(\eta) & = & \lim_{n\rightarrow\infty}\F_n(\eta)=
\lim_{n\rightarrow\infty}(\sum \frac{1}{|G|}\int_{\widehat{U}^\prime}
\frac{\sqrt{-1}}{2\pi}F_{a_n}\wedge\eta_{\widehat{U}})\\
         & = & \sum_i n_i\cdot (\sum \frac{1}{|G|}\int_{C_{\widehat{U},i}}
\eta_{\widehat{U}})=\sum_i n_i\cdot\int_{\Sigma_i}\hat{f}^\ast_i\eta.
\end{eqnarray*}
Thus $c_1(E)=\sum_i n_i\cdot PD(C_i)$.

Finally, as in Taubes \cite{T3}, if a subset $\Omega\subset X$ is
contained in $\alpha_n^{-1}(0)$ for all $n$, then $\Omega$ is contained
in $C=\cup_i C_i$ also.

\newpage

\centerline{\bf Appendix A:
Dimension of the Seiberg-Witten Moduli Space}

\vspace{5mm}

We begin with a brief review on the index theorem in Kawasaki \cite{Ka}.

Let $X$ be an orbifold (compact and connected), and $P$ be an elliptic
operator over $X$. In order to state the index theorem, we first
introduce the space $\widetilde{X}\equiv\{(p,(g)_{G_p})\mid p\in X,
g\in G_p\}$, where $G_p$ is the isotropy group at $p$, and $(g)_{G_p}$
is the conjugacy class of $g$ in $G_p$. The following properties of
$\widetilde{X}$ are easily verified (cf. \cite{Ka}, compare also \cite{CR}).

\begin{itemize}
\item $\widetilde{X}$ has a canonical orbifold structure, with a
canonical map $i:\widetilde{X}\rightarrow X$: let $(V_p,G_p)$ be a local
uniformizing system at $p\in X$, then $(V_p^g,Z_p(g))$,
where $V_p^g\subset V_p$ is the
fixed-point set of $g\in G_p$ and $Z_p(g)\subset G_p$ is the centralizer
of $g$, is a local uniformizing system at $(p,(g)_{G_p})\in\widetilde{X}$,
and the map $i:\widetilde{X}\rightarrow X$ is defined by the collection
of embeddings $\{(V_p^g,Z_p(g))\hookrightarrow (V_p,G_p)\mid p\in X,g\in
G_p\}$.
\item $\widetilde{X}$ is a disjoint union of compact, connected orbifolds
of various dimensions, containing the orbifold $X$ as the component
$\{(p,(1)_{G_p})\mid p\in X\}$: $\widetilde{X}=
\bigsqcup_{(g)\in T}X_{(g)}$ with $X_{(1)}=X$, where $T=\{(g)\}$ is the
set of equivalence classes of $(g)_{G_p}$ with the equivalence relation
$\sim$ defined as follows: $(h)_{G_q}\sim (g)_{G_p}$ if $q$ is contained
in a local uniformizing system centered at $p$ and $h\mapsto g$ under
the natural injective homomorphism $G_q\rightarrow G_p$ which is
defined only up to conjugation by an element of $G_p$.
\end{itemize}

We remark that the orbifold structure on $\widetilde{X}$ is in a more
general sense that the group action in a local uniformizing
system is not required to be effective. For such an orbifold, we shall
use the following convention: the fundamental class of the
orbifold, whenever it exists, equals the fundamental class of the
underlying space divided by the order of the isotropy group at a
smooth point (cf. \cite{CR} and \S 2 in \cite{C2}).

Next we describe the characteristic classes involved in the index theorem.
Let $u=[\sigma(P)]$ be the class of the principal symbol of the elliptic
operator $P$ in the K-theory of $TX$. Then the pullback of $u$ via the
differential of the map $i:X_{(g)}\rightarrow X$, denoted by $u_{(g)}$,
is naturally decomposed as $\oplus_{0\leq\theta<2\pi} u_{(g),\theta}$ where
$u_{(g),\theta}$ is the restriction of $u_{(g)}$ to the
$\exp(\sqrt{-1}\theta)$-eigenbundle of $g\in (g)$. We set
$$
\mbox{ch}_{(g)} u_{(g)}\equiv\sum_\theta\exp(\sqrt{-1}\theta)
\mbox{ch}\;u_{(g),\theta}\in H^\ast_c(TX_{(g)};\C).
$$

On the other hand, the normal bundle $N_p^g$ of $V_p^g\hookrightarrow V_p$
patches up to define an orbifold vector bundle $N_{(g)}\rightarrow X_{(g)}$,
and the decomposition $N_p^g=N_p^g(-1)\oplus_{0<\theta<\pi}N_p^g(\theta)$,
where $N_p^g(-1)$, $N_p^g(\theta)$ are the $(-1)$-eigenspace and
$\exp(\sqrt{-1}\theta)$-eigenspace of $g$ respectively, defines a
natural decomposition of orbifold vector bundles $N_{(g)}=
N_{(g)}(-1)\oplus_{0<\theta<\pi}N_{(g)}(\theta)$.

Now let $R,S_\theta$ be the characteristic classes over $\C$ of the orthogonal
group and unitary group, which are defined by the power series
\begin{eqnarray*}
& &\{\prod_i(\frac{1+\exp x_i}{2})(\frac{1+\exp (-x_i)}{2})\}^{-1},\\
& &\{\prod_i(\frac{1-\exp(y_i+\sqrt{-1}\theta)}{1-\exp(\sqrt{-1}\theta)})
(\frac{1-\exp (-y_i-\sqrt{-1}\theta)}{1-\exp(-\sqrt{-1}\theta)})\}^{-1}
\end{eqnarray*}
respectively. We set
$$
I_{(g)}\equiv R(N_{(g)}(-1))\prod_{0<\theta<\pi} \frac{S_\theta(N_{(g)}
(\theta))\tau(TX_{(g)}\otimes_\R\C)}{\det(1-(g)|_{N_{(g)}})}
\in H^\ast(X_{(g)};\C),
$$
where $\tau=\prod_i x_i(1-\exp(-x_i))^{-1}$ is the Todd class, and
$\det(1-(g)|_{N_{(g)}})$ is the constant function on $X_{(g)}$ which
equals $\det(1-g|_{N_p^g})$ at $(p,(g)_{G_p})\in X_{(g)}$. Note that
when $X$ is almost complex, $N_{(g)}$ is an orbifold complex vector
bundle, and there is a compatible decomposition $N_{(g)}=
\oplus_{0<\theta<2\pi}N_{(g)}(\theta)$. In this case, it is easily
seen that
$$
I_{(g)}=\prod_{0<\theta<2\pi}\frac{S_\theta(N_{(g)}(\theta))
\tau(TX_{(g)}\otimes_\R\C)}{\det(1-(g)|_{N_{(g)}})}\in H^\ast
(X_{(g)};\C).
$$

\vspace{3mm}

\noindent{\bf Theorem}\hspace{3mm} (Kawasaki \cite{Ka})
{\em
$$
{\em index }\;P=\sum_{(g)\in T}(-1)^{\dim X_{(g)}}\langle
{\em ch}_{(g)} u_{(g)}\cdot I_{(g)},[TX_{(g)}]\rangle
$$
where $u=[\sigma(P)]$. {\em(}Here the orientation for the fundamental
class $[TX_{(g)}]$ is given according to the {\em(}now standard{\em)}
convention in Atiyah-Singer \cite{AS}.{\em)}
}

\vspace{3mm}

For the rest of the appendix, we shall consider specifically
the case where $X$ is an almost complex 4-orbifold (which is in
the classical sense that the local group actions are effective),
and where $P$ is either the Dirac operator associated to a $Spin^\C$
structure on $X$, or the de Rham operator, or the signature operator.

First, the index of the Dirac operator associated to a $Spin^\C$ structure
on $X$. Recall that the almost complex structure of $X$ defines a canonical
$Spin^\C$ structure $S_{+}\oplus S_{-}$, $S_{+}=\I\oplus K_X^{-1}$, $S_{-}=TX$,
where $\I$ is the trivial orbifold complex line bundle, $K_X$ is the canonical
bundle, and $TX$ is the tangent bundle which, with the given almost
complex structure, is viewed as an orbifold $\C^2$-bundle. Any other $Spin^\C$
structure has the form $(S_{+}\oplus S_{-})\otimes E$ for
an orbifold complex line bundle $E$ over $X$.

Let $P=P_{Dirac}^E:C^\infty(S_{+}\otimes E)\rightarrow C^\infty(S_{-}
\otimes E)$ be the Dirac operator. We shall determine the contribution
to the index of $P$ from each component $X_{(g)}$ of $\widetilde{X}$.
Let $l=\dim_{\C} X_{(g)}$. Then $l=0,1,2$, where $X_{(g)}=p/G_p$ for a
singular point $p\in X$ when $l=0$, $X_{(g)}$ is 2-dimensional,
pseudoholomorphic when $l=1$, and $X_{(g)}=X_{(1)}=X$ when
$l=2$. In any event, the orbifold principal $U(2)$-bundle associated to
the almost complex structure reduces to an orbifold principal
$U(l)\times U(2-l)$-bundle when restricted to $X_{(g)}$ via the map
$i:X_{(g)}\rightarrow X$, and there is an orbifold principal $H$-bundle
$F$ over $X_{(g)}$, $H\equiv U(l)\times U(2-l)\times U(1)\subset U(3)$,
such that $TX_{(g)}=F\times_H\C^l$ and $E|_{X_{(g)}}=F\times_H\C$,
where $\C^l,\C$ are $H$-modules via $\C^l=\C^l\times \{0\}\subset\C^3$
and $\C=\{0\}\times\C\subset\C^3$. Moreover, let
$M_{+}=(\I\oplus \Lambda^2\C^2)\otimes\C$, $M_{-}=\C^2\otimes\C$
be the $H$-modules where $\I$ is the 1-dimensional trivial module,
$\C^2=\C^2\times\{0\}\subset\C^3$ and $\C=\{0\}\times\C\subset\C^3$.
Then $S_{+}\otimes E|_{X_{(g)}}=F\times_H M_{+}$ and
$S_{-}\otimes E|_{X_{(g)}}=F\times_H M_{-}$, so that $u_{(g)}$, the pullback
of the symbol class of $P$ via the differential of the map
$i:X_{(g)}\rightarrow X$, is an elliptic symbol class associated to
the $H$-structure (cf. Atiyah-Singer \cite{AS}).

There is a linear action of $g\in (g)$ on $M_{+}$ and $M_{-}$
associated to the bundles $TX|_{X_{(g)}},E|_{X_{(g)}}$.
Let $M_{+}=\oplus_{0\leq\theta< 2\pi} M_{+,\theta}$,
$M_{-}=\oplus_{0\leq\theta< 2\pi} M_{-,\theta}$ be the
corresponding decompositions into $\exp(\sqrt{-1}\theta)$-eigenspaces.
Let $\psi:H^\ast(X_{(g)};\C)\rightarrow H^\ast_c(TX_{(g)};\C)$ be
the Thom isomorphism. Then the contribution to the index of $P$
from $X_{(g)}$ is
\begin{eqnarray*}
&   & (-1)^{2l}\langle\mbox{ch}_{(g)} u_{(g)}\cdot I_{(g)},
[TX_{(g)}]\rangle\\
& = & (-1)^{l}\langle\psi^{-1}(\mbox{ch}_{(g)} u_{(g)})
\cdot I_{(g)},[X_{(g)}]\rangle\\
& = & (-1)^{l}\frac{\sum_\theta (\exp(\sqrt{-1}\theta)
\mbox{ch}\;M_{+,\theta}-\exp(\sqrt{-1}\theta)\mbox{ch}\;M_{-,\theta})}
{x_1\cdots x_l}(F)I_{(g)}[X_{(g)}]
\end{eqnarray*}
where $x_1\cdots x_l=1$ when $l=0$.

In the above formula, the symbol class of the Dirac operator contributes
through the modules $M_{+},M_{-}$. Similarly, in order to determine the
index formulae for the other geometric differential operators on $X$,
it suffices to write down the corresponding modules.

Let's look at the de Rham operator $d+d^\ast$, whose index is the Euler
characteristic $\chi(X)$. The modules are $N_{+}=(\I\oplus\Lambda^2\R^4\oplus
\Lambda^4\R^4)\otimes_\R\C$ and $N_{-}=(\R^4\oplus\Lambda^3\R^4)
\otimes_\R\C$. Because of the almost complex structure, $\R^4\otimes_\R\C=
\C^2\oplus\overline{\C}^2$, and if set
$W=\Lambda^2(\C_1\oplus\overline{\C_2})$ where $\C^2=\C_1\oplus\C_2$,
then we may rewrite $N_{+}=4\I\oplus\Lambda^2\C^2\oplus\Lambda^2
\overline{\C}^2\oplus W\oplus \overline{W}$, $N_{-}=2(\C^2\oplus
\overline{\C}^2)$.

For the signature operator whose index is $\mbox{sign}(X)$, the signature
of $X$, the modules are $Q_{+}=\Lambda^2_{+}\R^4\otimes_\R\C$ and
$Q_{-}=\Lambda^2_{-}\R^4\otimes_\R\C$. With the almost complex structure,
we may rewrite $Q_{+}=\I\oplus\Lambda^2\C^2\oplus\Lambda^2\overline{\C}^2$
and $Q_{-}=\I\oplus W\oplus\overline{W}$.

With these preparations, we give a formula in the following proposition for 
the dimension $d(E)$ of the moduli space of Seiberg-Witten equations 
associated to the $Spin^\C$-structure given by an orbifold complex line 
bundle $E$. 

\vspace{4mm}

\noindent{\bf Proposition}\hspace{3mm}
{\em $d(E)=2\cdot {\em index}\;P^E_{Dirac}
-\frac{1}{2}(\chi(X)+{\em sign}(X))=I_0+I_1+I_2$ where
$$
I_0=\sum_{\{(g)|\dim_{\C}X_{(g)}=0\}}\frac{2(\exp(\sqrt{-1}
\theta_{E,(g)})-1)}{(1-\exp(-\sqrt{-1}\theta_{1,(g)}))(1-
\exp(-\sqrt{-1}\theta_{2,(g)}))}[X_{(g)}]
$$
\begin{eqnarray*}
I_1 & = &\sum_{\{(g)|\dim_{\C}X_{(g)}=1\}}(\frac{2\exp(\sqrt{-1}
\theta_{E,(g)})c_1(E)}{1-\exp(-\sqrt{-1}\theta_{(g)})}\\
 & & +\frac{(\exp(\sqrt{-1}\theta_{E,(g)})-1)c_1(TX_{(g)})}
{1-\exp(-\sqrt{-1}\theta_{(g)})}\\
    &  & -\frac{2\exp(-\sqrt{-1}\theta_{(g)})
(\exp(\sqrt{-1}\theta_{E,(g)})-1)c_1(N_{(g)})}
{(1-\exp(-\sqrt{-1}\theta_{(g)}))^2})[X_{(g)}]
\end{eqnarray*}
and
$$
I_2=(c_1^2(E)-c_1(E)c_1(K_X))[X].
$$
In the above equations, $\exp(\sqrt{-1}\theta_{E,(g)})$
denotes the eigenvalue of $g\in (g)$ acting on the orbifold complex line
bundle $E|_{X_{(g)}}$, $\exp(\sqrt{-1}\theta_{1,(g)})$ and 
$\exp(\sqrt{-1}\theta_{2,(g)})$ denote the eigenvalues of $g\in (g)$ acting
on the normal bundle $N_{(g)}$ of $X_{(g)}$ when $\dim_{\C}X_{(g)}=0$,
and $\exp(\sqrt{-1}\theta_{(g)})$ denotes the eigenvalue of $g\in (g)$
acting on $N_{(g)}$ when $\dim_{\C}X_{(g)}=1$. Moreover, $X_{(g)}$ and $X$
are given with the canonical orientation as almost complex orbifolds. 
}

\vspace{3mm}

\pf
We first consider the case when the orbifold complex line bundle $E=\I$
is trivial, and show that $2\cdot\mbox{index}\;P^{\I}_{Dirac}
-\frac{1}{2}(\chi(X)+\mbox{sign}(X))=0$.

To this end, set $M_{+}^0=\I\oplus\Lambda^2\C^2$ and $M_{-}^0=\C^2$,
and let $M_{+}^0=\oplus_{0\leq\theta< 2\pi}M_{+,\theta}^0$ and $M_{-}^0
=\oplus_{0\leq\theta< 2\pi} M_{-,\theta}^0$ be the decompositions into
$\exp(\sqrt{-1}\theta)$-eigenspaces. Then the contribution to 
$\frac{1}{2}(\chi(X)+\mbox{sign}(X))$ that comes from $X_{(g)}$ is
\begin{eqnarray*}
(-1)^{l}\frac{\sum_\theta (\exp(\sqrt{-1}\theta)
\mbox{ch}\;M_{+,\theta}^0-\exp(\sqrt{-1}\theta)\mbox{ch}\;M_{-,\theta}^0)}
{x_1\cdots x_l}(F_0)I_{(g)}[X_{(g)}] +\\
(-1)^{l}\frac{\sum_\theta (\exp(-\sqrt{-1}\theta)
\mbox{ch}\;\overline{M_{+,\theta}^0}-\exp(-\sqrt{-1}\theta)\mbox{ch}\;
\overline{M_{-,\theta}^0})}{x_1\cdots x_l}(F_0)I_{(g)}[X_{(g)}],\\
\end{eqnarray*}
where $l=\dim_{\C}X_{(g)}$, $F_0$ is the orbifold principal
$U(l)\times U(2-l)$-bundle over $X_{(g)}$ which is the reduction
when restricted to $X_{(g)}$ of the orbifold principal $U(2)$-bundle
over $X$ associated to the almost complex structure.

Observe that only the terms of degree $2l$ in
$$
\sum_\theta (\exp(\sqrt{-1}\theta)
\mbox{ch}\;M_{+,\theta}^0-\exp(\sqrt{-1}\theta)\mbox{ch}\;M_{-,\theta}^0)
$$
could possibly make a contribution, which is invariant under $x_i\mapsto
-x_i$. Moreover, $\theta\mapsto -\theta$ under $(g)\mapsto (g^{-1})$,
and $I_{(g)}=I_{(g^{-1})}$ under the identification $X_{(g)}=X_{(g^{-1})}$.
Hence the following two expressions 
\begin{eqnarray*}
(-1)^{l}\frac{\sum_\theta (\exp(-\sqrt{-1}\theta)
\mbox{ch}\;\overline{M_{+,\theta}^0}-\exp(-\sqrt{-1}\theta)\mbox{ch}\;
\overline{M_{-,\theta}^0})}{x_1\cdots x_l}(F_0)I_{(g)}[X_{(g)}],\\
(-1)^{l}\frac{\sum_\theta (\exp(\sqrt{-1}\theta)
\mbox{ch}\;M_{+,\theta}^0-\exp(\sqrt{-1}\theta)\mbox{ch}\;M_{-,\theta}^0)}
{x_1\cdots x_l}(F_0)I_{(g^{-1})}[X_{(g^{-1})}]\\
\end{eqnarray*} 
are equal, from which it follows easily that 
$$
2\cdot\mbox{index}\;P_{Dirac}^{\I}-\frac{1}{2}(\chi(X)+\mbox{sign}(X))=0.
$$

For the general case, notice that
$$
2\cdot \mbox{index}\;P_{Dirac}^E
-\frac{1}{2}(\chi(X)+\mbox{sign}(X))=2(\mbox{index}\;P_{Dirac}^E
-\mbox{index}\;P_{Dirac}^{\I}).
$$
The formula follows easily from direct evaluation of the right hand side.

\hfill $\Box$

\vspace{3mm}

\noindent{\bf Proof of Lemma 3.8}

\vspace{2mm}

The dimension $d(E)$ of the Seiberg-Witten moduli space corresponding to $E$
equals
$$
2\cdot \mbox{index }P_{Dirac}^E-\frac{1}{2}(\chi(X)+\mbox{sign }(X)).
$$
According to the proposition, it is the sum of 
$c_1(E)\cdot c_1(E)-c_1(K_X)\cdot
c_1(E)$ with a term contributed by the singular point $p_0$, which
can be written as $\frac{1}{|G|}\sum_{g\in G,g\neq 1}\chi(g)$, with
$$
\chi(g)=\frac{2(\rho(g)-1)}{(1-\exp(-\sqrt{-1}\theta_{1,g}))
(1-\exp(-\sqrt{-1}\theta_{2,g}))},
$$
where $\rho:G\rightarrow \s^1$ is the representation given in
Lemma 3.6, and $\exp(\sqrt{-1}\theta_{i,g})$, $i=1,2$, are the
two eigenvalues of $g\in G\subset U(2)$. The evaluation of this term
constitutes the main task in the proof, which will be done
case by case according to the type of $G$.

\vspace{2mm}

(1) $G=\langle Z_{2m},Z_{2m};\widetilde{D}_n,\widetilde{D}_n\rangle$.
Let $h=\mu_{2m} I\in Z_{2m}$, and let $x,y$ be generators of $\widetilde{D}_n$
satisfying $x^2=y^n=(xy)^2=-1$. Then $G\setminus\{1\}=
\Lambda_1\sqcup\Lambda_2\sqcup\Lambda_3$, where $\Lambda_1=\sqcup_{l=0}^{n-1}
\Lambda_1^{(l)}$ with $\Lambda_1^{(l)}=\{h^k y^l\mid k=1,2,\cdots,2m-1\}$,
$\Lambda_2=\{h^k xy^l\mid k=0,1,\cdots,2m-1,l=0,1,\cdots,n-1\}$, and
$\Lambda_3=\{y^l\mid l=1,2,\cdots,n-1\}$. Note that $\chi(g)=0$ for any
$g\in\Lambda_3$.

We first calculate $\sum_{g\in\Lambda_1}\chi(g)$. To this end, we set,
for each $s=1,2,\cdots,2n$, $S_{l,s}(t)\equiv
\sum_{k=1}^{2m-1}\frac{(\mu_{2m}^k\mu_{2n}^{-l})^s}
{1-\mu_{2m}^{-k}\mu_{2n}^{-l} t}$. Introduce $[s]$ such that
$s\equiv [s] \pmod{2m}, 0\leq [s]\leq 2m-1$. Then
\begin{eqnarray*}
S_{l,s}(t) & = &
\sum_{k=1}^{2m-1}(\mu_{2m}^k\mu_{2n}^{-l})^s\sum_{j=0}^\infty
(\mu_{2m}^{-k}\mu_{2n}^{-l})^j t^j\\ &= &
\sum_{j=0}^\infty
\mu_{2n}^{-ls-lj}\sum_{k=1}^{2m-1}(\mu_{2m}^{s-j})^k t^j\\
            & = &
-\sum_{j=0}^\infty \mu_{2n}^{-ls-lj} t^j+ 2m \mu_{2n}^{-ls-l[s]}t^{[s]}
\sum_{j=0}^\infty (\mu_{2n}^{-l} t)^{2mj}\\
& = & -\frac{\mu_{2n}^{-ls}}{1-\mu_{2n}^{-l}
t}+2m\frac{\mu_{2n}^{-ls-l[s]} t^{[s]}}{1-(\mu_{2n}^{-l} t)^{2m}}.
\end{eqnarray*}
We consider separately when $l=0$ or $l\neq 0$.
\begin{eqnarray*}
\sum_{s=1}^{2n}S_{0,s}(1) & = & \sum_{s=1}^{2n}S_{0,s}(t)|_{t=1}
=\sum_{s=1}^{2n}(-\frac{1}{1-t}+\frac{2m
t^{[s]}}{1-t^{2m}})|_{t=1}\\
& = & \sum_{s=1}^{2n}\frac{-\sum_{i=0}^{2m-1} t^i +2m
t^{[s]}}{1-t^{2m}}|_{t=1}=\sum_{s=1}^{2n}
\frac{-\sum_{i=0}^{2m-1} i +2m [s]}{-2m}\\
& = & \sum_{s=1}^{2n}(\frac{1}{2}(2m-1)-[s]).
\end{eqnarray*}
For each $l\neq 0$, note that
$\sum_{s=1}^{2n}(\mu_{2n}^{-l})^s=0$,  so that $\sum_{s=1}^{2n}
S_{l,s}(1)=\sum_{s=1}^{2n}\frac{2m\mu_{2n}^{-ls-l[s]}}
{1-(\mu_{2n}^{-l})^{2m}}$.
Introduce $j_s$ which is uniquely defined by $0\leq j_s\leq n-1$
and $s+[s]+2mj_s\equiv 0\pmod{2n}$. Then
\begin{eqnarray*}
S(t) &\equiv &\sum_{s=1}^{2n}\sum_{l=1}^{n-1}\frac{\mu_{2n}^{-ls-l[s]}}
{1-\mu_{2n}^{-2ml}t}=\sum_{s=1}^{2n}\sum_{l=1}^{n-1}\sum_{j=0}^\infty
(\mu_{2n}^{-2ml})^j t^j\\
& = & \sum_{s=1}^{2n}\sum_{j=0}^\infty\sum_{l=1}^{n-1}
(\mu_{2n}^{-(s+[s]+2mj)})^l
t^j=\sum_{s=1}^{2n}(\sum_{j=0}^\infty -t^j+n t^{j_s}
\sum_{j=0}^\infty t^{nj})\\
& = & \sum_{s=1}^{2n}(-\frac{1}{1-t}+\frac{nt^{j_s}}{1-t^n}).
\end{eqnarray*}
Similarly, $S(1)=\sum_{s=1}^{2n}(\frac{1}{2}(n-1)-j_s)$, and
$$
\sum_{l=1}^{n-1}\sum_{s=1}^{2n}S_{l,s}(1)=2m
S(1)=\sum_{s=1}^{2n}(m(n-1)-2mj_s).
$$
With these preparations, now observe that if we set $S_l\equiv
\frac{1}{2}\sum_{g\in\Lambda_1^{(l)}}\chi(g)$, then
\begin{eqnarray*}
S_l & = & \sum_{k=1}^{2m-1}\frac{(\mu_{2m}^k\mu_{2n}^{-l})^{2n}-1}
{(1-\mu_{2m}^{-k}\mu_{2n}^{-l})(1-\mu_{2m}^{-k}\mu_{2n}^l)}\\
& = & \sum_{s=1}^{2n}\sum_{k=1}^{2m-1}\frac{(\mu_{2m}^k\mu_{2n}^{-l})^s}
{1-\mu_{2m}^{-k}\mu_{2n}^{-l}}=\sum_{s=1}^{2n}S_{l,s}(1).
\end{eqnarray*}
Hence
$$
\sum_{g\in\Lambda_1}\chi(g)=2\sum_{l=0}^{n-1}S_l=\sum_{s=1}^{2n}
((2mn-1)-2([s]+2mj_s)).
$$

Next we calculate $\sum_{g\in\Lambda_2}\chi(g)$. First of all,
\begin{eqnarray*}
\sum_{g\in\Lambda_2}\chi(g) & = & \sum_{l=0}^{n-1}\sum_{k=0}^{2m-1}
\frac{2((\mu_{2m}^k)^{2n}(-1)^n-1)}{(1-\mu_{2m}^{-k}\sqrt{-1})
(1-\mu_{2m}^{-k}(\sqrt{-1})^{-1})}\\
& = & \sum_{s=1}^{2n}\sum_{l=0}^{n-1}\sum_{k=0}^{2m-1}
\frac{2(\mu_{2m}^k\sqrt{-1})^s}{1-\mu_{2m}^{-k}\sqrt{-1}}.
\end{eqnarray*}
Set $S_s(t)\equiv \sum_{k=0}^{2m-1}\frac{(\mu_{2m}^k\sqrt{-1})^s}
{1-(\mu_{2m}^{-k}\sqrt{-1})t}$, $s=1,2,\cdots,2n$. Then
\begin{eqnarray*}
S_s(t) & = & \sum_{k=0}^{2m-1}(\mu_{2m}^k\sqrt{-1})^s\sum_{j=0}^\infty
(\mu_{2m}^{-k}\sqrt{-1})^j t^j\\
& = & \sum_{j=0}^\infty\sum_{k=0}^{2m-1}
(\mu_{2m}^{s-j})^k (\sqrt{-1})^{s+j} t^j\\
& = & 2m \sum_{j=0}^\infty (\sqrt{-1})^{s+[s]+2mj}\cdot t^{[s]+2mj}\\
& = & \frac{2m (\sqrt{-1})^{s+[s]}\cdot t^{[s]}}{1-(\sqrt{-1}\cdot
t)^{2m}},\\
\end{eqnarray*}
and
$$
\sum_{g\in\Lambda_2}\chi(g)=n\sum_{s=1}^{2n}2S_s(1)=2mn\sum_{s=1}^{2n}
(\sqrt{-1})^{s+[s]}.
$$

In order to evaluate $\sum_{g\in\Lambda_1}\chi(g)$ and
$\sum_{g\in\Lambda_2}\chi(g)$, we consider the cases where $m>n$
and $m<n$ separately.

Suppose $m>n$. In this case, we have $[s]=s$ for any $s=1,2,\cdots,2n$.
Furthermore, $s=[s]$ and $m,n$ being relatively prime imply
that $s\mapsto j_s$ is a surjective, two to one correspondence from
$\{1,2,\cdots,2n\}$ to $\{0,1,\cdots,n-1\}$. It then follows easily
from these observations that
$$
\sum_{g\in\Lambda_1}\chi(g)=4n(m-n-1), \sum_{g\in\Lambda_2}\chi(g)=0.
$$
Hence
\begin{eqnarray*}
d(E) & = & c_1(E)\cdot c_1(E)-c_1(K_X)\cdot c_1(E)+
\frac{1}{|G|}\sum_{g\neq 1}\chi(g)\\
& = & \frac{n}{m}+\frac{m+1}{m}+\frac{1}{4mn}\cdot 4n(m-n-1)=2.
\end{eqnarray*}

Now consider the case where $m<n$. We introduce $\delta,r$ satisfying
$n=\delta m+r$, $0\leq r\leq m-1$. Then a simple inspection shows that
$$
\sum_{g\in\Lambda_2}\chi(g)=2mn((-1)^\delta-1).
$$
In order to evaluate $\sum_{g\in\Lambda_1}\chi(g)$, we introduce, for
each $s=1,2,\cdots,2n$, $k_s$ which obeys $s+[s]+2mj_s=2nk_s$. Then
one can easily check that $k_s$ satisfies $1\leq k_s\leq m$. Now observe
that for any $l=0,1,\cdots,2\delta-1$, $s\mapsto k_s$ is injective,
hence surjective, if $lm+1\leq s\leq lm+m$. Summing up the equations
$s+[s]+2mj_s=2nk_s$ from $s=1$ to $s=2\delta m$, we have
$$
\sum_{s=1}^{2\delta m} ([s]+2m j_s)
=\sum_{s=1}^{2\delta m}2nk_s-\sum_{s=1}^{2\delta m} s
=\delta m(2nm+2r-1).
$$
If $m\neq 1$, we need to consider the rest of the values
of $s$, $s=2\delta m+1,\cdots,2\delta m+2r$. For this part, observe
that $j_r=0$, $j_{2r}=\delta$, and for any $1\leq s\leq r-1$,
we have the relation $2m(j_s+j_{2r-s})=2n(k_s+k_{2r-s}-2)$, which implies
easily that $j_s+j_{2r-s}=n$. Thus
$$
\sum_{s=2\delta m+1}^{2n}([s]+2mj_s)=r(2r+1)+2m(rn-n+\delta).
$$
Putting things altogether, we have
$$
\sum_{g\in\Lambda_1}\chi(g)=4mn-4n(r+1).
$$
Finally, when $m<n$, we have
\begin{eqnarray*}
d(E) & = & \frac{n}{m}+\frac{m+1}{m}+\frac{1}{4mn}(4mn-4n(r+1)
+2mn((-1)^\delta-1))\\
& = &\delta+2+\frac{1}{2}((-1)^\delta-1).
\end{eqnarray*}

(2) $G=\langle Z_{4m},Z_{2m};\widetilde{D}_n,C_{2n}\rangle$.
Let $h=\mu_{4m} I\in Z_{4m}$, and let $x,y$ be generators of
$\widetilde{D}_n$ satisfying $x^2=y^n=(xy)^2=-1$. Introduce
$\bar{h}=h^2$, $\bar{x}=hx$, and $\bar{y}=y$. Then $G\setminus\{1\}=
\Lambda_1\sqcup\Lambda_2\sqcup\Lambda_3$, where $\Lambda_1=
\{\bar{h}^k \bar{y}^l\mid k=1,2,\cdots,2m-1,l=0,1,\cdots,n-1\}$, $\Lambda_2
=\{\bar{h}^k \bar{x}\bar{y}^l\mid k=0,1,\cdots,2m-1,l=0,1,\cdots,n-1\}$, and
$\Lambda_3=\{\bar{y}^l\mid l=1,2,\cdots,n-1\}$. Again, we have $\chi(g)=0$
for any $g\in\Lambda_3$.

Note that $\sum_{g\in\Lambda_1}\chi(g)$ is the same as in the previous case,
so we only need to evaluate $\sum_{g\in\Lambda_2}\chi(g)$, for which
a similar calculation shows that
$$
\sum_{g\in\Lambda_2}\chi(g)=2mn\sum_{s=1}^{2n}\mu_{4m}^{s-[s]}
(\sqrt{-1})^{s+[s]}.
$$
A simple inspection, with the fact that $m$ is even this time, shows that
$\sum_{g\in\Lambda_2}\chi(g)=0$ when $m>n$, and when $m<n$, 
$$
\sum_{g\in\Lambda_2}\chi(g)=2mn ((-1)^\delta-1),
\mbox{ where } n=\delta m+r, 0\leq r\leq m-1.
$$
By the same calculation, $d(E)=2$ if $m>n$, and $d(E)=\delta+2+\frac{1}{2}
((-1)^\delta-1)$ if $m<n$.

(3) $G=\langle Z_{2m},Z_{2m};\widetilde{T},\widetilde{T}\rangle$.
Let $h=\mu_{2m} I\in Z_{2m}$, and let $x,y$ be generators of $\widetilde{T}$
satisfying $x^2=y^3=(xy)^3=-1$. Then $\sum_{g\neq 1}\chi(g)=S_0+S_1+S_2$,
where $S_0=\sum_{k=1}^{2m-1}\chi(h^k)$,
$$
S_1=\sum_{g=y,xy,x^{-1}yx,yx}\sum_{l=1}^2\sum_{k=0}^{2m-1}\chi(h^kg^l),
$$
and
$$
S_2=\sum_{g=x,y^{-1}xy,y^{-2}xy^2}\sum_{k=0}^{2m-1}\chi(h^k g).
$$
Let $[s]$ be defined by $s\equiv [s]$ and $0\leq [s]\leq 2m-1$. Then
a similar calculation shows that
\begin{eqnarray*}
S_0 & = & \sum_{s=1}^{12}(2m-1-2[s])\\
S_1 & = & 16m\sum_{l=1}^2\sum_{s=1}^{12}\frac{\mu_6^{(s+[s])l}}
{1-\mu_3^{ml}}\\
S_2 & = & 6m\sum_{s=1}^{12}\mu_4^{s+[s]}.\\
\end{eqnarray*}
When $m>6$, we have $S_0=24(m-7)$, $S_1=S_2=0$, so that
$$
d(E)=\frac{6}{m}+\frac{m+1}{m}+\frac{1}{24m}\cdot 24(m-7)=2.
$$
When $m<6$, then either $m=1$ or $m=5$. For $m=1$, $S_0=S_1=S_2=0$,
and $d(E)=8$. For $m=5$, $S_0=12$, $S_1=0$, and $S_2=-60$, which
gives $d(E)=2$.

(4) $G=\langle Z_{6m},Z_{2m};\widetilde{T},\widetilde{D}_2\rangle$.
Let $h=\mu_{2m} I\in Z_{2m}$, and let $x,y$ be generators of $\widetilde{T}$
satisfying $x^2=y^3=(xy)^3=-1$. Introduce $\bar{h}=h^3$, $\bar{x}=x$,
and $\bar{y}=hy$. Then $\sum_{g\neq 1}\chi(g)=S_0+S_1+S_2$,
where $S_0=\sum_{k=1}^{2m-1}\chi(\bar{h}^k)$,
$$
S_1=\sum_{g=\bar{y},\bar{x}\bar{y},
\bar{x}^{-1}\bar{y}\bar{x},\bar{y}\bar{x}}
\sum_{l=1}^2\sum_{k=0}^{2m-1}\chi(\bar{h}^kg^l),
$$
and
$$
S_2=\sum_{g=\bar{x},\bar{y}^{-1}\bar{x}\bar{y},
\bar{y}^{-2}\bar{x}\bar{y}^2}\sum_{k=0}^{2m-1}\chi(\bar{h}^k g).
$$
A similar calculation, with the fact that $m$ is divisible by $3$,
shows that
\begin{eqnarray*}
S_0 & = & \sum_{s=1}^{12}(2m-1-2[s])\\
S_1 & = & 16m\sum_{l=1}^2\sum_{s=1}^{12}
\frac{\mu_6^{(s+[s])l}\mu_{6m}^{(s-[s])l}}
{1-\mu_3^{-l}}\\
S_2 & = & 6m\sum_{s=1}^{12}\mu_4^{s+[s]}.\\
\end{eqnarray*}
When $m>6$, we have $S_0=24(m-7)$, $S_1=S_2=0$, so that $d(E)=2$.
When $m<6$, then $m=3$, and in this case, $S_0=S_2=0$ and $S_1=-96$,
which also gives $d(E)=2$.

(5) $G=\langle Z_{2m},Z_{2m};\widetilde{O},\widetilde{O}\rangle$.
Let $h=\mu_{2m} I\in Z_{2m}$, and let $x,y$ be generators of $\widetilde{O}$
satisfying $x^2=y^4=(xy)^3=-1$. Recall that $\widetilde{O}$ is the union
of three cyclic subgroups of order $8$ generated by $y,xyx$ and $y^2x$,
four cyclic subgroups of order $6$ generated by $xy,yx,y^3xy^2$ and
$y^2xy^3$, and six cyclic subgroups of order $4$ generated by
$x,yxy^3,y^2xy^2,y^3xy^2x,xy^2xy^3$ and $y^2xy^3x$, where these
subgroups only intersect at $\{1,-1\}$. Consequently, we have
$\sum_{g\neq 1}\chi(g)=S_0+S_1+S_2+S_3$, where $S_0=\sum_{k=1}^{2m-1}
\chi(h^k)$,
\begin{eqnarray*}
S_1 & = &\sum_g\sum_{l=1}^3\sum_{k=0}^{2m-1}\chi(h^kg^l), \mbox{ where }
g \mbox{ has order } 8\\
S_2 & = & \sum_g\sum_{l=1}^2\sum_{k=0}^{2m-1}\chi(h^kg^l),\mbox{ where }
g \mbox{ has order } 6\\
S_3 & = & \sum_g\sum_{k=0}^{2m-1}\chi(h^kg),\mbox{ where }
g \mbox{ has order } 4.\\
\end{eqnarray*}
A similar calculation shows that
\begin{eqnarray*}
S_0 & = & \sum_{s=1}^{24}(2m-1-2[s])\\
S_1 & = & 12m\sum_{l=1}^3\sum_{s=1}^{24}
\frac{\mu_8^{(s+[s])l}}{1-\mu_4^{ml}}\\
S_2 & = & 16m\sum_{l=1}^2\sum_{s=1}^{24}
\frac{\mu_6^{(s+[s])l}}{1-\mu_3^{ml}}\\
S_3 & = & 12m \sum_{s=1}^{24}\mu_4^{s+[s]}.\\
\end{eqnarray*}
When $m>12$, $S_0=48(m-13)$ and $S_1=S_2=S_3=0$, so that $d(E)=2$.
When $m<12$, then $m=1,5,7$ or $11$. A direct calculation 
shows that $d(E)=2$ in all this cases except for $m=1$, for which
$d(E)=14$. Below we list the results of $S_0,S_1,S_2$ and $S_3$ for
the sake of records.
\begin{itemize}
\item $m=1$: $S_0=S_1=S_2=S_3=0$.
\item $m=5$: $S_0=16$, $S_1=-240$, $S_2=-160$, and $S_3=0$.
\item $m=7$: $S_0=20$, $S_1=84$, $S_2=-224$, and $S_3=-168$.
\item $m=11$: $S_0=36$, $S_1=132$, $S_2=0$, and $S_3=-264$.
\end{itemize}

(6) $G=\langle Z_{2m},Z_{2m};\widetilde{I},\widetilde{I}\rangle$.
Let $h=\mu_{2m} I\in Z_{2m}$, and let $x,y$ be generators of $\widetilde{I}$
satisfying $x^2=y^5=(xy)^3=-1$. Then $\sum_{g\neq 1}\chi(g)=S_0+S_1+S_2+S_3$,
where $S_0=\sum_{k=1}^{2m-1}\chi(h^k)$, and 
\begin{eqnarray*}
S_1 & = &\sum_g\sum_{l=1}^4\sum_{k=0}^{2m-1}\chi(h^kg^l),\;
g \mbox{ is one of the six elements of order } 10\\
S_2 & = & \sum_g\sum_{l=1}^2\sum_{k=0}^{2m-1}\chi(h^kg^l),\;
g \mbox{ is one of the ten elements of order } 6\\
S_3 & = & \sum_g\sum_{k=0}^{2m-1}\chi(h^kg),\;
g \mbox{ is one of the fifteen elements of order } 4.\\
\end{eqnarray*}
A similar calculation shows that
\begin{eqnarray*}
S_0 & = & \sum_{s=1}^{60}(2m-1-2[s])\\
S_1 & = & 24m\sum_{l=1}^4\sum_{s=1}^{60}
\frac{\mu_{10}^{(s+[s])l}}{1-\mu_5^{ml}}\\
S_2 & = & 40m\sum_{l=1}^2\sum_{s=1}^{60}
\frac{\mu_6^{(s+[s])l}}{1-\mu_3^{ml}}\\
S_3 & = & 30m \sum_{s=1}^{60}\mu_4^{s+[s]}.\\
\end{eqnarray*}
When $m>30$, $S_0=120(m-31)$ and $S_1=S_2=S_3=0$, so that $d(E)=2$.
When $m<30$, then $m=1,7,11,13,17,19,23$ or $29$. A direct calculation
shows that $d(E)=2$ for all cases except for $m=1$, in which case $d(E)=32$,
and for $m=7$, in which case $d(E)=4$. We record the calculation for
$S_0,S_1,S_2$ and $S_3$ below.
\begin{itemize}
\item $m=1$: $S_0=S_1=S_2=S_3=0$.
\item $m=7$: $S_0=32$, $S_1=-672$, $S_2=-560$, and $S_3=0$.
\item $m=11$: $S_0=64$, $S_1=-1584$, $S_2=-880$, and $S_3=0$.
\item $m=13$: $S_0=128$, $S_1=-1248$, $S_2=-1040$, and $S_3=0$.
\item $m=17$: $S_0=156$, $S_1=-816$, $S_2=0$, and $S_3=-1020$.
\item $m=19$: $S_0=308$, $S_1=912$, $S_2=-1520$, and $S_3=-1140$.
\item $m=23$: $S_0=420$, $S_1=0$, $S_2=0$, and $S_3=-1380$.
\item $m=29$: $S_0=108$, $S_1=1392$, $S_2=0$, and $S_3=-1740$.
\end{itemize}

\hfill $\Box$

\newpage

\centerline{\bf Appendix B:
Green's Function for the Laplacian on Orbifolds}

\vspace{5mm}

We shall follow the relevant discussion in Chapter 4 of Aubin \cite{Au}
for Green's function on compact Riemannian manifolds. 

Let $(X,g)$ be a compact, closed, oriented $n$-dimensional Riemannian
orbifold. For any $p,q\in X$, we define the distance between $p$ and $q$,
denoted by $d(p,q)$, to be the infinimum of the lengths of all piecewise
$C^1$ paths connecting $p$ and $q$. Then $(X,d)$ is a complete metric
space. Moreover, there is a geodesic $\gamma$ between $p,q$ such that
$d(p,q)=\mbox{length}(\gamma)$. (A $C^1$ path $f:[a,b]\rightarrow X$ is
called a (parametrized) geodesic in $X$ if $f$ is locally lifted to a
geodesic in a uniformizing system.) Observe that when $p,q$ are both in
a uniformized open subset $U$ and the geodesic $\gamma$ with
$d(p,q)=\mbox{length}(\gamma)$ is also contained in $U$, then $\gamma$
may be lifted to a geodesic in $\widehat{U}$, where $(\widehat{U},G)$
uniformizes $U$. This implies that in the said circumstance,
$$
d(p,q)=\min_{\{h_1,h_2\in G\}} \hat{d}(h_1\cdot\hat{p},h_2\cdot\hat{q})
$$
where $\hat{p},\hat{q}$ are any inverse image of $p,q$ in $\widehat{U}$,
and $\hat{d}$ is the distance function on $\widehat{U}$ (note that
$\widehat{U}$ has a Riemannian metric canonically determined by $g$).
On the other hand, for each $p\in X$, there is a $\delta(p)>0$, called
the injectivity radius at $p$, such that for any $0<\delta\leq\delta(p)$,
the set $U_p(\delta)=\{q\in X\mid d(p,q)<\delta\}$ is uniformized by
$(\widehat{U_p(\delta)},G_p)$ where $\widehat{U_p(\delta)}$ is a convex
geodesic ball of radius $\delta$ centered at the inverse image of $p$
and $G_p$ is the isotropy group at $p$ acting linearly on
$\widehat{U_p(\delta)}$. We point out that $\delta(p)\rightarrow 0$
as $p\rightarrow q$ for any $q$ with $|G_q|>|G_p|$. In particular, there is
no positive uniform lower bound for the injectivity radius on an orbifold
with a nonempty orbifold-point set.

\vspace{4mm}

\noindent{\bf Theorem 1}{\em\hspace{2mm}
Let $\Delta=d^\ast d$ be the Laplacian on $(X,g)$ and let
$n\equiv\dim X\geq 2$. There exists $G(p,q)$, a Green's function for
the Laplacian which has the following properties:
\begin{itemize}
\item [{(1)}] For all $C^2$ functions $\varphi$ on $X$,
$$
\varphi(p)=\mbox{Vol}(X)^{-1}\int_X\varphi+\int_X G(p,\cdot)\Delta\varphi,
$$
where $\mbox{Vol}(X)$ is the volume of $X$.
\item [{(2)}] $G(p,q)$ is a smooth function on $X\times X$ minus the diagonal.
\item [{(3)}] There exists a decomposition $G(p,q)=G_0(p,q)+G_1(p,q)$ such that
\begin{itemize}
\item $G_1(p,q)$ is continuous in both variables and $C^2$ in $q$.
\item There exist a $\delta_0>0$ and a set $\U$ of finitely many
uniformizing systems on $X$ with the following significance: For any
$p\in X$, there is a uniformizing system $(\widehat{U},G,\pi)\in\U$ and
a $G$-invariant open subset $\widehat{U}^\prime\subset\widehat{U}$,
such that {\em (i)} $p\in\pi(\widehat{U}^\prime)$, {\em (ii)} the
support of the function $q\mapsto G_0(p,q)$ is contained in
$\pi(\widehat{U})$ {\em (}more precisely, $\{q\mid G_0(p,q)\neq 0\}\subseteq
\{q\mid d(p,q)\leq (n+1)\delta_0${\em )},
and {\em (iii)} $\widehat{U}$ contains the closed ball of radius $\delta_0$
centered at any $\hat{p}\in\pi^{-1}(p)$. Moreover, the function $G_0(p,q)$ is
the descendant of $\sum_{h\in G}\widehat{G_0}(h\cdot\hat{p},\hat{q})$
for some function $\widehat{G_0}(\hat{p},\hat{q})$, which is, {\em (i)}
continuous on $\widehat{U}^\prime\times\widehat{U}$ minus the subset
$\{(\hat{p},\hat{q})\mid \hat{p}=\hat{q}\}$, {\em (ii)} $C^2$ in $\hat{q}$,
{\em (iii)} invariant under the diagonal action of $G$,
and {\em (iv)} satisfying the following estimates for a constant $z>0$:
$$
\begin{array}{c}
|\widehat{G_0}(\hat{p},\hat{q})|\leq z(1+|\log\hat{d}(\hat{p},\hat{q})|)
\hspace{2mm} \mbox{ for } n=2 \hspace{2mm} \mbox{ and } \\
|\widehat{G_0}(\hat{p},\hat{q})|\leq\frac{z}{\hat{d}(\hat{p},\hat{q})^{n-2}}
\hspace{2mm} \mbox{ for } n>2, \mbox{ with } \\
|\nabla_{\hat{q}}\widehat{G_0}(\hat{p},\hat{q})|\leq
\frac{z}{\hat{d}(\hat{p},\hat{q})^{n-1}},
|\nabla_{\hat{q}}^2\widehat{G_0}(\hat{p},\hat{q})|\leq
\frac{z}{\hat{d}(\hat{p},\hat{q})^{n}}.
\end{array}
$$
{\em (}Here $\hat{d}$ is the distance function on $\widehat{U}$.{\em )}
\item [{(4)}] There exists a constant $C$ such that $G(p,q)\geq C$.
Green's function is defined up to a constant, so one may arrange
so that $G(p,q)\geq 1$.
\item [{(5)}] The map $q\mapsto \int_X G(p,q)$ is constant. One may choose to
have $\int_X G(p,q)=0$.
\item [{(6)}] Green's function is symmetric: $G(p,q)=G(q,p)$.
\end{itemize}
\end{itemize}
}
\vspace{2mm}

\pf
Choose finitely many points $p_i\in X$ such that $X=\bigcup_i
U_{p_i}(N^{-1}\delta_i)$ with $\delta_i\equiv N^{-1}\delta(p_i)$, where
$\delta(p_i)$ is the injectivity radius at $p_i$ and $N$ is any fixed
integer which is no less than $3n$ (recall $n=\dim X$).
The set $\U$ is chosen to be the set of uniformizing system of
$U_{p_i}(\delta(p_i))$. Set $\delta_0\equiv \min_i\delta_i$.

Now given any $p\in X$, suppose $p$ is contained in $U^\prime\equiv
U_{p_i}(\delta_i)$ for some $i$. We denote by $(\widehat{U}^\prime,G^\prime,
\pi^\prime)$ the uniformizing system of $U^\prime$, and by
$(\widehat{U},G,\pi)$ the uniformizing system of
$U\equiv U_{p_i}(\delta(p_i))$,
which is an element of $\U$ by definition. Note that $G^\prime=G=G_{p_i}$.
With these understood, we define a function $\widehat{H_0}(\hat{p},\hat{q})$
on $\widehat{U}^\prime\times\widehat{U}$ for $\hat{p}\neq\hat{q}$, such that
$$
\begin{array}{c}
\widehat{H_0}(\hat{p},\hat{q})=-(2\pi)^{-1}\rho(r)\log r \hspace{3mm}
\mbox{ for } n=2 \hspace{3mm} \mbox{ and }\\
\widehat{H_0}(\hat{p},\hat{q})=[(n-2)\omega_{n-1}]^{-1}\rho(r)r^{2-n}
\hspace{3mm} \mbox{ for } n>2,
\end{array}
$$
where $r=\hat{d}(\hat{p},\hat{q})$, $\rho(r)$ is a fixed cut-off function
which equals zero when $r\geq\delta_0$, and $\omega_{n-1}$ is the volume
of the unit sphere in $\R^n$. It is clear that $\widehat{H_0}(\hat{p},\hat{q})$
is invariant under the diagonal action of $G$. We define
$$
\widehat{H}(\hat{p},\hat{q})=\sum_{h\in G}\widehat{H_0}(h\cdot\hat{p},\hat{q}),
$$
which is invariant under the action of $G\times G$. Let $H(p,q)$ be the
descendant of $\widehat{H}(\hat{p},\hat{q})$, which is defined on
$U^\prime\times U$ for $p\neq q$. We extend $H(p,q)$ over $q\in X$ by zero.
Now observe that if we use a different element of $\{U_{p_i}(\delta_i)\}$
for $U^\prime$, we end up with the same function $q\mapsto H(p,q)$. Hence
we obtain a function $H(p,q)$ on $X\times X$ minus the diagonal. Note
that $\{q\mid H(p,q)\neq 0\}\subseteq \{q\mid d(p,q)\leq\delta_0\}$.

Define $\Gamma_1(p,q)=-\Delta_{q}H(p,q)$ and for each $1\leq j\leq n$
define inductively
$$
\Gamma_{j+1}(p,q)=\int_X\Gamma_j(p,\cdot)\Gamma_1(\cdot,q).
$$
We note that $\{q\mid\Gamma_j(p,q)\neq 0\}\subseteq
\{q\mid d(p,q)\leq j\cdot\delta_0\}$. Now suppose $p\in U^\prime$.
If we let $\widehat{\Gamma}_1(\hat{p},\hat{q})=-\Delta_{\hat{q}}
\widehat{H_0}(\hat{p},\hat{q})$ and for each $1\leq j\leq n$, define
$$
\widehat{\Gamma}_{j+1}(\hat{p},\hat{q})=\int_{\widehat{U}}
\widehat{\Gamma}_j(\hat{p},\cdot)\widehat{\Gamma}_1(\cdot,\hat{q})
$$
inductively, then each $\widehat{\Gamma}_j$ is invariant
under the diagonal action of $G$. Moreover, each $\Gamma_j(p,q)$ is
the descendant of $\sum_{h\in G}\widehat{\Gamma}_j(h\cdot\hat{p},\hat{q})$.

As in Aubin \cite{Au} (cf. Prop. 4.12 in \cite{Au}),
$\widehat{\Gamma}_{n+1}(\hat{p},\hat{q})$ is $C^1$ on $\widehat{U}^\prime
\times\widehat{U}$. Hence $\Gamma_{n+1}(p,q)$ is $C^1$ on $X\times X$.
Fix each $p\in X$, we solve the Laplacian equation
$$
\Delta_{q}F(p,q)=\Gamma_{n+1}(p,q)-Vol(X)^{-1}.
$$
(Note that for each $p\in X$, $\int_X (\Gamma_{n+1}(p,\cdot)-Vol(X)^{-1})=0$,
cf. \cite{Au}.) Then we define
$$
G(p,q)=H(p,q)+\sum_{j=1}^n\int_X \Gamma_j(p,\cdot)H(\cdot,q)+ F(p,q),
$$
and by adding an appropriate constant to $F(p,q)$, we arrange to have
for all $p\in X$
$$
\int_X G(p,\cdot)=0.
$$

One can argue as in Aubin \cite{Au} that $G(p,q)$ is a Green's function
for the Laplacian which has the properties described in the theorem.
Particularly, in the decomposition
$$
G(p,q)=G_0(p,q)+G_1(p,q)
$$
in Thoerem 1 (3),
we have $G_0(p,q)=H(p,q)+\sum_{j=1}^n\int_X \Gamma_j(p,\cdot)H(\cdot,q)$
and $G_1(p,q)=F(p,q)$. (Note that $\{q\mid G_0(p,q)\neq 0\}
\subseteq \{q\mid d(p,q)\leq (n+1)\delta_0$.) Moreover, for any
$p\in U^\prime$, $G_0(p,q)$ is clearly the descendant of
$\sum_{h\in G}\widehat{G_0}(h\cdot\hat{p},\hat{q})$ where
$$
\widehat{G_0}(\hat{p},\hat{q})=\widehat{H_0}(\hat{p},\hat{q})+\sum_{j=1}^n
\int_{\widehat{U}}\widehat{\Gamma}_j(\hat{p},\cdot)
\widehat{H_0}(\cdot,\hat{q}).
$$
The estimates
$$
\begin{array}{c}
|\widehat{G_0}(\hat{p},\hat{q})|\leq z(1+|\log\hat{d}(\hat{p},\hat{q})|)
\hspace{2mm} \mbox{ for } n=2 \hspace{2mm} \mbox{ and } \\
|\widehat{G_0}(\hat{p},\hat{q})|\leq\frac{z}{\hat{d}(\hat{p},\hat{q})^{n-2}}
\hspace{2mm} \mbox{ for } n>2, \mbox{ with } \\
|\nabla_{\hat{q}}\widehat{G_0}(\hat{p},\hat{q})|\leq
\frac{z}{\hat{d}(\hat{p},\hat{q})^{n-1}},
|\nabla_{\hat{q}}^2\widehat{G_0}(\hat{p},\hat{q})|\leq
\frac{z}{\hat{d}(\hat{p},\hat{q})^{n}}
\end{array}
$$
in Thoerem 1 (3) follow immediately from the definition of
$\widehat{H_0}(\hat{p},\hat{q})$.

\hfill $\Box$

\newpage

\centerline{\bf Appendix C: 
Proof of Lemma 1.4}

\vspace{5mm}

Various transversality theorems for harmonic forms on a compact, closed 
Riemannian manifold were proved in Honda \cite{Ho}. The machinery 
developed therein 
can be properly adapted to deal with the present situation.

First of all, we recall the relevant discussion in \cite{Ho} regarding the case
of self-dual harmonic forms on a (compact, closed) $4$-manifold. Suppose $M$ is
a smooth $4$-manifold with $b_2^{+}(M)\neq 0$. Let $\text{Met}^l(M)$ be the 
space of $C^l$-H\"{o}lder metrics on $M$ for a sufficiently large 
non-integer $l$, let $Q^{+}$ be the space of pairs $(\omega,g)$ where 
$g\in\text{Met}^l(M)$ and $\omega$ is a nontrivial self-dual $g$-harmonic 
form, 
and let $\bigwedge^{2,+}\rightarrow \text{Met}^l(M)\times M$ be the vector
bundle whose fiber at $(g,x)$ is $\bigwedge_g^{2,+}(T_x^\ast M)$, the space
of $2$-forms at $x$ which is self-dual with respect to $g$. Then the
transversality of the following evaluation map 
$$
ev:Q^{+}\times M\rightarrow \bigwedge^{2,+}, \hspace{2mm}
((\omega,g),x)\mapsto (\omega(x), (g,x))
$$
to the zero section was studied in \cite{Ho}. The relevant results are
summarized below. For any $((\omega,g),x)\in Q^{+}\times M$ where
$\omega(x)=0$, it was shown that the differential 
$$
(ev_x)_\ast: T_{(\omega,g)}Q^{+}\rightarrow \bigwedge_g^{2,+}(T_x^\ast M),
\hspace{2mm} (v,h)\mapsto v(x)
$$
is surjective. Here $T_{(\omega,g)}Q^{+}$ is the tangent space of $Q^{+}$
at $(\omega,g)$, which consists of pairs $(v,h)$, where 
$h\in C^l(Sym^2(T^\ast M))$ and $v$ is a $2$-form, self-dual at $x$ with
respect to $g$ and satisfying the equation 
$$
\Delta_g v+ \frac{d}{dt}(\Delta_{g+th})|_{t=0}\omega=0.
$$
(Here $\Delta_g$ is the Laplacian associated to a metric $g$.) As a 
consequence, for a generic metric a nontrivial self-dual harmonic form
has only regular zeroes, which consist of a disjoint union of embedded
circles in $M$. 

In order to adapt the argument to the present situation, we recall some
relevant details about the surjectivity of $(ev_x)_\ast:
T_{(\omega,g)}Q^{+}\rightarrow \bigwedge_g^{2,+}(T_x^\ast M)$.
Suppose $(v,h)\in T_{(\omega,g)}Q^{+}$ and $v$ is orthogonal to the
space of $g$-harmonic $2$-forms. Then one can solve for $v$ from $h$ by
$$
v(x)=-(\Delta_g)^{-1}(\frac{d}{dt}(\Delta_{g+th})|_{t=0}\omega)(x)=
\pm\int_{M}\langle dd^\ast G_g(x,y),\ast(D_h\ast)\omega(y)\rangle_g,
$$
where $G_g(x,y)$ is the Green's function for $\Delta_g$, and
$D_h\ast$ is shorthand for $\frac{d}{dt}(\ast_{g+th})|_{t=0}$. 
(Here the Hodge star $\ast$ and the integration are with respect to
the metric $g$.) For any $x\in M$ with $\omega(x)=0$, one considers the map
$\Psi_x:C^l(Sym^2(T^\ast M))\rightarrow \bigwedge_g^{2,+}(T_x^\ast M)$ where
$$
\Psi_x: h\mapsto v(x)=\pm\int_{M}\langle dd^\ast G_g(x,y),
\ast(D_h\ast)\omega(y)\rangle_g.
$$
Then it is clear that the surjectivity of $(ev_x)_\ast$ is a consequence of
that of $\Psi_x$. 

To explain the basic ingredients in the proof of surjectivity of $\Psi_x$,
we first introduce some notations. For any $0\neq u\in\R^4$, let 
$R_u:\R^4\rightarrow \R^4$ be the reflection in $u$, and let
$\bigwedge^2 R_u:\bigwedge^2(\R^4)\rightarrow \bigwedge^2(\R^4)$ be the
induced isomorphism. For any skew-symmetric $4\times 4$ matrix $A$ and 
symmetric $4\times 4$ matrix $H$, let 
$
i_A(H)=HA+AH-\frac{1}{2}tr(H)\cdot A.
$  
Then the proof of surjectivity of $\Psi_x$ goes roughly as follows.

\begin{itemize}
\item Assume $g$ is flat. Then for any $y\neq x$ nearby, a
direct calculation shows that 
$\langle dd^\ast G_g(x,y),\cdot \rangle_g
=\frac{C}{|y-x|^4}\cdot\bigwedge^2 R_{y-x}$
for a constant $C\neq 0$, and 
that $\ast(D_h\ast)\omega(y)=i_{\omega(y)}(h)$ where $\omega(y)$, $h$ are
regarded canonically as a skew-symmetric and a symmetric $4\times 4$ matrix 
respectively. Moreover, one can verify that $\bigwedge^2 R_u:
\bigwedge^{2,\pm}(\R^4)\rightarrow\bigwedge^{2,\mp}(\R^4)$, and that when 
$\omega(y)\neq 0$ and is self-dual, the image of $i_{\omega(y)}$ is the 
space of 
anti-self-dual $2$-forms. The surjectivity of $\Psi_x$ follows in this case
by letting $h$ be a $\delta$-function like element centered at $y$. 
\item In general, use the asymptotic expension of the Green's
function $G_g(x,y)$ near the diagonal, whose leading term is the Green's 
function for a flat metric, to reduce the proof to the previous case. 
\end{itemize}

With these preparations, we now give a proof of Lemma 1.4.

\vspace{2mm}

Let $W$ be an oriented smooth $s$-cobordism of elliptic $3$-manifolds
as in Lemma 1.4. We attach a semi-cylinder $[0,+\infty)\times (\s^3/G)$
to the positive end of $W$ and cone-off the negative end by $\B^4/G$.
The resulting space, denoted by $\hat{W}$, is an orbifold with one isolated
singular point $p_0$ and a cylindrical end over $\s^3/G$. We shall fix
a Riemannian metric $g_0$ on $\hat{W}$, which is flat near $p_0$ and
is the product metric $dt^2+h_0$ on the semi-cylinder 
$[0,+\infty)\times (\s^3/G)$. Here $h_0$ stands for the standard metric 
on $\s^3/G$ which has a constant sectional curvature of $1$. Fix a
$T>1$ and a sufficiently large non-integer $l$, we will consider 
$\text{Met}^l_T(\hat{W})$, the space of $C^l$-H\"{o}lder metrics 
on $\hat{W}$ which equals $g_0$ on $[T,+\infty)\times (\s^3/G)$.

We fix an identification $\R^4=\C^2=\h$
so that $G$ as a subgroup of $\phi(\s^1\times\s^3)$ is canonically
regarded as a subgroup of $U(2)$. Let $z_1,z_2$ be the standard coordinates
on $\C^2$, and let $\tilde{\alpha}$ be the pull-back of the $1$-form
$\sqrt{-1}\sum_{i=1}^2(z_id\bar{z}_i-\bar{z}_idz_i)$ to $\s^3$. Then
it is easy to check that $\tilde{\alpha}$ obeys $d\tilde{\alpha}
=2(\ast\tilde{\alpha})$ with respect to the standard metric on $\s^3$,
and consequently, $d(\exp(2t)\tilde{\alpha})$ is a self-dual $2$-form
on $\R\times\s^3$ with respect to the corresponding product metric.
Note that $\tilde{\alpha}$ is invariant under the action of $G$. Let
$\alpha$ be the descendant of $\tilde{\alpha}$ to $\s^3/G$. 

\vspace{3mm}

\noindent{\bf Proposition}\hspace{2mm}
{\em
For each $g\in \text{Met}^l_T(\hat{W})$, there is a unique self-dual
$g$-harmonic $2$-form $\omega_g$ which has the following properties.
\begin{itemize}
\item [{(1)}] On $[T,+\infty)\times (\s^3/G)$, $\omega_g-d(\exp(2t)\alpha)
=d\alpha_t$ where $\alpha_t$ is a $1$-form on $\s^1/G$ such that
$\alpha_t$ and $\frac{d}{dt}\alpha_t$ converge to zero exponentially fast 
as $t\rightarrow +\infty$.
\item [{(2)}] For a generic $g$, $\omega_g$ has only regular zeroes
in the complement of the singular point $p_0$ and 
$[T,+\infty)\times (\s^3/G)$. 
\item [{(3)}] For a generic $g$ which is sufficiently close to $g_0$ near
$p_0$, $\omega_g(p_0)\neq 0$. 
\end{itemize}
}

Assuming the validity of the proposition, we obtain the $2$-form
$\omega$ claimed in Lemma 1.4 as follows. Consider the $g_0$-harmonic
$2$-form $\omega_{g_0}$ first. Since $\omega_{g_0}-d(\exp(2t)\alpha)$
converges to zero exponentially fast as $t\rightarrow +\infty$, and
note that $d(\exp(2t)\alpha)$ is a symplectic form, there is a $\tau_0>0$
such that $\omega_{g_0}$ is symplectic on 
$[\tau_0,+\infty)\times (\s^3/G)$. We fix a $T\geq \tau_0$, then 
for all $g\in\text{Met}^l_T(\hat{W})$ sufficiently close to $g_0$,
$\omega_g$ is symplectic on $[T,+\infty)\times (\s^3/G)$. We pick such
a $g$ which is generic. Then by (2) and (3) of the proposition, $\omega_g$
has only regular zeroes, which are in the complement of $p_0$ and 
$[T,+\infty)\times (\s^3/G)$. In particular, $\omega_g(p_0)\neq 0$. 
Let $\hat{\omega}_g$ be the $2$-form obtained from $\omega_g$ by replacing
$\omega_g=d(\exp(2t)\alpha+\alpha_t)$ with 
$d(\exp(2t)\alpha+(1-\rho_\tau)\alpha_t)$ on the cylindrical end of $\hat{W}$,
where $\rho_\tau$ is a cut-off function for a sufficiently large 
$\tau\geq T+2$ which equals $1$ on $t\geq \tau$ and equals $0$ on 
$t\leq \tau-1$. Then $\hat{\omega}_g=d(\exp(2t)\alpha)$ on the cylindrical
end of $\hat{W}$ where $t\geq\tau$, and $\hat{\omega}_g(p_0)={\omega}_g(p_0)
\neq 0$, so that by the equivariant Darboux' theorem $\hat{\omega}_g$
is equivalent near $p_0$ to the standard symplectic form on $\B^4/G$. 
It is clear that
$\hat{\omega}_g$ yields a $2$-form $\omega$ on the $s$-cobordism $W$
as described in Lemma 1.4.

\vspace{3mm}

\noindent{\bf Proof of Proposition}

\vspace{3mm}

For any $g\in \text{Met}^l_T(\hat{W})$,
we denote by $\bigwedge^{2,+}_g$ and $\bigwedge^1_g$ the associated bundle
of self-dual $2$-forms and $1$-forms on $\hat{W}$ respectively.
Consider the subspaces $\E_g$, $\F_g$ of the (weighted) Sobolev spaces
$H^2_1(\bigwedge^{2,+}_g)$, $H^2_0(\bigwedge^1_g)$,
where $\E_g$ is the closure of self-dual $2$-forms which
equal $dt\wedge\alpha_t+\ast_3\alpha_t$ on the cylindrical end with
$d^\ast_3\alpha_t=0$ for $t\geq T$, and $\F_g$ is the closure of
co-closed $1$-forms which can be written as $f_tdt+\beta_t$ on the
cylindrical end with $f_t=0$ and $d^\ast_3\beta_t=0$ when $t\geq T$.
(Here $\ast_3$ is the Hodge star and $d_3^\ast$ is the co-exterior
differential on $\s^3/G$, both with respect to the standard metric $h_0$.)
Then the differential operator $\ast_g d:\E_g\rightarrow\F_g$, which is of
form $\frac{d}{dt}-\ast_3 d_3$ on the cylindrical end, defines a Fredholm
operator, cf. \cite{LM}. (Note that there exist no harmonic $1$-forms
on $\s^3/G$, so that we can choose $\delta=0$ in the weight factor
$\exp(\delta t)$ of the weighted Sobolev spaces. )

\vspace{3mm}

\noindent{\bf Lemma 1}{\em\hspace{2mm}
$\ast_g d:\E_g\rightarrow\F_g$ has a trivial kernel and cokernel.
}

\vspace{3mm}

Assuming the validity of Lemma 1, we obtain the self-dual $g$-harmonic
form $\omega_g$ for each $g\in \text{Met}^l_T(\hat{W})$ as follows. Let $\beta$
be a self-dual $2$-form on $\hat{W}$ obtained by multiplying 
$d(\exp(2t)\alpha)$ with a cut-off function which equals $1$ on $t\geq 1$. 
Then $\ast_g d\beta\in \F_g$, and hence there exists a $u_g\in \E_g$
such that $\ast_g d u_g=-\ast_g d\beta$. We set $\omega_g\equiv u_g+\beta$. 
It is clear that $\omega_g$ is self-dual $g$-harmonic. To show that
$\omega_g$ has the property in (1) of the proposition, we note that 
on $[T,+\infty)\times (\s^3/G)$, $\ast_g d u_g=-\ast_g d\beta=0$, which
implies that there is a $\alpha_t$, with $\omega_g-d(\exp(2t)\alpha)=u_g
=d\alpha_t$, such that $\alpha_t$ and $\frac{d}{dt}\alpha_t$ converge 
to zero exponentially fast as $t\rightarrow +\infty$ (see the proof of 
Lemma 1 below). Finally, observe that such an $\omega_g$ is unique. 

To prove (2) and (3) of the proposition, we consider the evaluation maps 
$$
ev_q:(\omega_g,g)\mapsto\omega_g(q)\in\bigwedge_g^{2,+}(T^\ast\hat{W}_q), 
\mbox{ where }
q\in\hat{W}\setminus [T,+\infty)\times (\s^3/G).
$$
Suppose $\omega_g(q)=0$. Then as in \cite{Ho}, the differential of $ev_q$ 
is given by $(ev_q)_\ast(h)=v(q)\in\bigwedge_{g}^{2,+}(T^\ast\hat{W}_q)$,
where $(v,h)$ obeys 
$$
\Delta_{g}v +\frac{d}{dt}(\Delta_{g+th})|_{t=0}\omega_{g}=0.
$$
We shall prove first that for any $q\neq p_0$, $(ev_q)_\ast$ is surjective,
which gives (2) of the proposition by a standard argument. 

To this end, for any $\tau>T+2$, we set 
$W_\tau\equiv \hat{W}\setminus (\tau,+\infty)\times (\s^3/G)$, and let 
$DW_\tau$ be the double of $W_\tau$, which is given with the natural metric 
and orientation. Denote by $\Delta_\tau=d^\ast d+d^\ast d$ the Laplacian 
on $DW_\tau$, and let $\gamma_\tau$ be the first eigenvalue of $\Delta_\tau$ 
on the space of $L^2$ self-dual $2$-forms on $DW_\tau$. We will need the
following lemma. 

\vspace{3mm}

\noindent{\bf Lemma 2}{\em\hspace{2mm}
There exist a $\tau_0>T+2$ and a constant $c>0$ such that $\gamma_\tau\geq c$
for all $\tau\geq \tau_0$.
}

\vspace{3mm}

Assuming the validity of Lemma 2, we fix a $\tau\geq \tau_0+10$, and 
decompose $v=v_1+v_2$ with $v_1\equiv (1-\rho_\tau)v$ and 
$v_2\equiv \rho_\tau v$, where $\rho_\tau$ is a cut-off function
which equals $1$ on $t\geq \tau$. 
Then we have $\Delta_g v_1=-\frac{d}{dt}(\Delta_{g+th})|_{t=0}\omega_{g}
-\Delta_g v_2$. Note that $h$ is supported in $\hat{W}\setminus
(T,+\infty)\times (\s^3/G)$, so that the above equation may be regarded 
as an equation on $DW_\tau$ because $\Delta_g v_2=\Delta_g v=0$ on
$t\geq \tau$. Moreover, $|\Delta_g v_2|\leq c\cdot \exp(-\delta\tau)$ 
for some constants $c>0$ and $\delta>0$, where $c$ is a multiple of the
$C^0$-norm of $v$ on $\hat{W}\setminus (T,+\infty)\times (\s^3/G)$,
hence is bounded by a multiple of the norm of $h$ via the standard 
elliptic estimates. Therefore 
$|\Delta_g v_2|\leq C\cdot \exp(-\delta\tau)\cdot ||h||$. 
Let $G_\tau$ be the Green's function for the Laplacian $\Delta_\tau$
on $DW_\tau$ (the existence of the Green's function $G_\tau$ on the 
orbifold $DW_\tau$ is a straightforward generalization of that in the 
compact manifold case, cf. e.g. \cite{PR}), then for any 
$q\in\hat{W}\setminus [T,+\infty)\times (\s^3/G)$,
$$
v(q)=v_1(q)=\pm\int_{DW_\tau}\langle dd^\ast G_\tau(q,y),\ast(D_h\ast)
\omega_{g}(y)\rangle_g -(\Delta_\tau)^{-1}(\Delta_g v_2),
$$
where $D_h\ast=\frac{d}{dt}(\ast_{g+th})|_{t=0}$. By Lemma 2 and the
standard elliptic estimates, the last term $(\Delta_\tau)^{-1}(\Delta_g v_2)$
in the above equation is bounded by a multiple of $\exp(-\delta \tau)\cdot
||h||$, hence can be neglected by taking $\tau$ sufficiently large. 
The surjectivity of $(ev_q)_\ast:h\mapsto v(q)$ for $q\neq p_0$ follows 
from the surjectivity of 
$$
\Psi_q:h\mapsto \pm\int_{DW_\tau}\langle dd^\ast G_\tau(q,y),\ast(D_h\ast)
\omega_{g}(y)\rangle_g 
$$
as in \cite{Ho}, which we have recalled at the beginning. 

It remains to prove (3) of the proposition, i.e., for a generic metric 
$g$ which is sufficiently close to $g_0$, $\omega_g$ does not
vanish at the singular point $p_0$. To this end, we consider $g_0$ first,
and assume $\omega_{g_0}(p_0)=0$ (otherwise the claim is trivially true).
Identify a local uniformizing system at $p_0$ with $(\B^4,G)$, which 
is given with a flat metric. Then the bundle of self-dual $2$-forms has 
a local basis $\omega_0,\omega_1,\omega_2$, where $\omega_0=\sqrt{-1}\sum_i
dz_i\wedge d\bar{z}_i$, $\omega_1=Re(dz_1\wedge dz_2)$ and
$\omega_2=Im(dz_1\wedge dz_2)$. Note that $\omega_0$ is invariant
under the action of $G$. With this understood, we claim that for $g_0$,
$(ev_{p_0})_\ast$ is transverse to the subspace spanned by 
$\omega_1,\omega_2$. To see this, we pick a $y$ sufficiently close to 
$p_0$ such that $\omega_{g_0}(y)\neq 0$, and denote by $y_0,y_1,\cdots,y_N$ 
the inverse images of $y$ in $\B^4$. Then according to \cite{Ho} as we
recalled earlier, there exists a $h_0\in Sym^2\R^4$ such that 
$\bigwedge^2 R_{y_0}\circ i_{\omega_{g_0}(y_0)}(h_0)=\omega_0$. 
Let $h_i$, $i=0,1,\cdots,N$, be the orbit of $h_0$ under the action of $G$, 
then because $\omega_0$ is invariant under the action of $G$,
we have $\bigwedge^2 R_{y_i}\circ i_{\omega_{g_0}(y_i)}(h_i)=\omega_0$ for
$i=0,1,\cdots,N$. Now observe that $(ev_{p_0})_\ast: h\mapsto v(p_0)$ is 
given by the equation (with $\tau \gg 0$)
$$
(ev_{p_0})_\ast(h)=\pm\int_{DW_\tau}
\langle dd^\ast G_\tau(p_0,y),\ast(D_h\ast)
\omega_{g_0}(y)\rangle_{g_0}-(\Delta_\tau)^{-1}(\Delta_{g_0} v_2).
$$
It is clear that we can use $\{h_i\}$ to define
a $G$-equivariant section $h\in C^\infty(Sym^2 T^\ast \B^4)$, which is
supported in a small neighborhood of $\{y_i\}$, such that
the projection of $(ev_{p_0})_\ast(h)$ to the $\omega_0$ factor is
nonzero. Hence for $g_0$, $(ev_{p_0})_\ast$ is transverse to the subspace 
spanned by $\omega_1,\omega_2$, so is it for any $g$ sufficiently
close to $g_0$. As a corollary, let $Z$ be the subbundle spanned by 
$\omega_1,\omega_2$ over a sufficiently small, $G$-invariant neighborhood 
of $0\in\B^4$, then for any generic metric $g$ sufficiently close to $g_0$,
$\omega_g^{-1}(Z)$ is a $3$-dimensional manifold in $\B^4$ which is invariant 
under the action of $G$. If $0\in\omega_g^{-1}(Z)$, then the tangent space of 
$\omega_g^{-1}(Z)$ at $0\in\B^4$ is invariant under the action of $G$, 
which is possible only when $G=\{1,-1\}$. Hence when $G\neq\{1,-1\}$,
$\omega_g(p_0)$ is not contained in $Z$, therefore $\omega_g(p_0)\neq 0$. 

When $G=\{1,-1\}$, note that all $\omega_0,\omega_1,\omega_2$ are invariant
under the action of $G$. The above argument then shows that $(ev_{p_0})_\ast$
is surjective for $g_0$. It follows easily that for a generic
metric $g$ which is sufficiently close to $g_0$, if $\omega_g(p_0)=0$,
then one of the components of $\omega_g^{-1}(0)$ in $\hat{W}$ is a
compact $1$-dimensional manifold with boundary, whose boundary consists
of the singular point $p_0$. This is a contradiction. Hence the proposition.

\hfill $\Box$

\vspace{3mm}

The rest of the appendix is occupied by the proofs of Lemma 1 and Lemma 2.

\vspace{3mm}

\noindent{\bf Proof of Lemma 1}

\vspace{3mm}

We first prove that the kernel of $\ast_g d:\E_g\rightarrow\F_g$ is trivial.
By elliptic regularity, it suffices to show that if $\omega\in\E_g$ is
smooth and satisfies $\ast_g d\omega=0$, then $\omega=0$.

First of all, note that $\omega=dt\wedge\alpha_t+\ast_3\alpha_t$
on the cylindrical end with $d^\ast_3\alpha_t=0$ for $t\geq T$.
Moreover, the equation $\ast_g d\omega=0$ is equivalent to
$\frac{d}{dt}\alpha_t-\ast_3 d_3\alpha_t=0$. If we write $\alpha_t=
\sum_if_i(t)\alpha_i$, where $\{\alpha_i\}$ is a complete set of eigenforms
for the self-adjoint operator $\ast_3 d_3$ on the space of $L^2$ co-closed
$1$-forms on $\s^3/G$, with $\ast_3 d_3\alpha_i=\lambda_i\alpha_i$,
then the functions $f_i(t)$ satisfy $f_i^\prime(t)-\lambda_i f_i(t)=0$.
It follows easily, since $\omega$ has a bounded $L^2$-norm, that each
$f_i(t)=c_i\exp(\lambda_i t)$ for some constant $c_i$ with $\lambda_i<0$.
Set $\delta\equiv\min_i|\lambda_i|>0$. Then $|\alpha_t|\leq c\exp(-\delta t)$
for a constant $c>0$.

Since $H^2_{dR}(\hat{W})=0$, $\ast_g d\omega=0$ implies that there exists
a $1$-form $\gamma$ on $\hat{W}$ such that $d\gamma=\omega$. We write
$\gamma=f_t dt+g_t$ on the cylindrical end, then $f_t,g_t$ satisfy
$$
\frac{d}{dt}g_t=d_3f_t+\alpha_t \mbox{ and } d_3g_t=\ast_3\alpha_t.
$$
Set $\bar{f}_t\equiv \int_{t_0}^t f_sds$ and $\bar{\alpha}_t\equiv
\int_{t_0}^t\alpha_sds$. Then
$g_t=d_3\bar{f}_t+\bar{\alpha}_t+\mbox{constant}$.
With this understood, we have
$$
0=\int_{\hat{W}}d^\ast d\gamma\wedge\ast\gamma=\int_{\hat{W}}d\gamma
\wedge \ast d\gamma\pm\lim_{t\rightarrow +\infty}\int_{\{t\}\times
(\s^3/G)} d\gamma\wedge\gamma,
$$
where $d\gamma\wedge\gamma=\ast_3\alpha_t\wedge g_t=\ast_3\alpha_t\wedge
(d_3\bar{f}_t+\bar{\alpha}_t+\mbox{constant})$. Since $\alpha_t
\rightarrow 0$ exponentially fast along the cylindrical end, and
$\ast_3\alpha_t=d_3 g_t$, it follows easily that
$\lim_{t\rightarrow +\infty}\int_{\{t\}\times (\s^3/G)}
d\gamma\wedge\gamma=0$, which implies $\omega=d\gamma=0$.

Next we show that the cokernel of $\ast_g d:\E_g\rightarrow\F_g$ is
trivial. To this end, note that the formal adjoint of $\ast_g d$ is
$d^{+}$. By elliptic regularity, it suffices to prove that for any
smooth $1$-form $\theta\in\F_g$, $d^{+}\theta=0$ implies $\theta=0$.
Note that on the cylindrical end, $\theta=f_tdt+\beta_t$ where $f_t=0$
and $d^\ast_3\beta_t=0$ when $t\geq T$. Thus $d^{+}\theta=0$ implies
that $\frac{d}{dt}\beta_t+\ast_3 d_3\beta_t=0$. Similarly, there exists
a $\delta>0$,  such that $|\theta|=|\beta_t|\leq c\exp(-\delta t)$ for 
some constant $c>0$. As a consequence, since $d^{+}\theta=0$, we have
$$
\int_{\hat{W}} d\theta\wedge\ast d\theta=-\int_{\hat{W}}
d\theta\wedge d\theta=-\lim_{t\rightarrow +\infty}\int_{\{t\}\times
(\s^3/G)}\theta\wedge d\theta=0,
$$
and hence $d\theta=0$. Now $H^1_{dR}(\hat{W})=0$ implies that $\theta=
df$ for some smooth function $f$ on $\hat{W}$. On the cylindrical end,
$\theta=df=\frac{\partial f}{\partial t}\cdot dt+d_3f$, so that when
$t\geq T$, $\frac{\partial f}{\partial t}=0$. This implies that $f$
is bounded on $\hat{W}$, and therefore
$$
\int_{\hat{W}} |df|^2=\int_{\hat{W}}\langle d^\ast df,f\rangle
=\int_{\hat{W}}\langle d^\ast\theta,f\rangle=0,
$$
which implies $\theta=df=0$. This proves the lemma.

\hfill $\Box$

\vspace{1mm}

\noindent{\bf Proof of Lemma 2}

\vspace{3mm}

Suppose to the contrary that there exists a sequence $\tau_n\rightarrow
+\infty$ for which $\gamma_{\tau_n}\rightarrow 0$. Let $\omega_n$ be
a corresponding sequence of self-dual $2$-forms on $DW_{\tau_n}$ such that
$\Delta_{\tau_n}\omega_n=\gamma_{\tau_n}\omega_n$.

To fix the notation, we identify the cylindrical neck of $DW_{\tau_n}$ with
$[T,2\tau_n-T]\times (\s^3/G)$. We denote by $W_1,W_2$ the two
components of $DW_{\tau_n}\setminus (T+1,2\tau_n-T-1)\times (\s^3/G)$.
With this understood, we rescale each $\omega_n$ and pass to a subsequence
if necessary, so that the following conditions hold:
$$
\int_{W_2}|\omega_n|^2\leq\int_{W_1}|\omega_n|^2=1
$$
By the interior elliptic estimates and the fact that $\Delta_{\tau_n}
\omega_n=\gamma_{\tau_n}\omega_n$ with $\gamma_{\tau_n}$ bounded,
there exists a constant $M_0>0$,  such that $|\omega_n|\leq M_0$ holds on
$W_1\setminus (T+\frac{1}{2},T]\times (\s^3/G)$ and $W_2
\setminus [2\tau_n-T-1,2\tau_n-T-\frac{1}{2})\times (\s^3/G)$.

On the other hand, if we write $\omega_n=dt\wedge\alpha_{n,t}+\ast_3
\alpha_{n,t}$ and set 
$$
f_n(t)\equiv\int_{\{t\}\times (\s^3/G)}|\omega_n|^2
=\int_{\s^3/G}|\alpha_{n,t}|^2
$$ 
on the cylindrical neck $[T,2\tau_n-T]\times (\s^3/G)$, then 
$\Delta_{\tau_n}\omega_n=\gamma_{\tau_n}\omega_n$ implies that 
$$
-\frac{d^2}{dt^2}\alpha_{n,t}+
\Delta_3\alpha_{n,t}=\gamma_{\tau_n}\alpha_{n,t}
$$
(here $\Delta_3=d_3^\ast d_3+d_3d^\ast_3$ is the Laplacian on $\s^3/G$),
and consequently, we have
\begin{eqnarray*}
\frac{d^2f_n}{dt^2} & = & \int_{\s^3/G}\frac{d^2}{dt^2}|\alpha_{n,t}|^2
=2(\int_{\s^3/G}\langle\frac{d^2}{dt^2}\alpha_{n,t},\alpha_{n,t}\rangle+
\int_{\s^3/G}|\frac{d}{dt}\alpha_{n,t}|^2)\\
& = & 2(\int_{\s^3/G}\langle(\Delta_3-\gamma_{\tau_n})\alpha_{n,t},
\alpha_{n,t}\rangle+\int_{\s^3/G}|\frac{d}{dt}\alpha_{n,t}|^2)>0\\
\end{eqnarray*}
for sufficiently large $n>0$. By the maximum principle, $f_n$ will attain
its maximum at the end point $t=T$ or $2\tau_n-T$. This implies that
for any $a\in [T,2\tau_n-T-1]$, the integral $\int_{[a,a+1]\times
(\s^3/G)}|\omega_n|^2$ is uniformly bounded, hence by the interior elliptic
estimates, $|\omega_n|\leq M_1$ for some constant $M_1>0$ on the
cylindrical neck of $DW_{\tau_n}$. Set $M\equiv\max (M_0,M_1)$, then
$|\omega_n|\leq M$ on $DW_{\tau_n}$.

By the standard elliptic estimates, there exists a subsequence of
$\omega_n$ (still denoted by $\omega_n$ for simplicity), and a self-dual
$2$-form $\omega$ on $\hat{W}$, such that
$\omega_n\rightarrow\omega$ in $C^\infty$ on any given compact subset
of $\hat{W}$. In particular, the $2$-form $\omega$ obeys (1)
$\int_{\hat{W}\setminus [T+1,+\infty)\times (\s^3/G)} |\omega|^2=1$,
and (2) $\Delta_g\omega=0$ and $|\omega|\leq M$ on $\hat{W}$. 
The lemma is proved by observing that (2) above
implies that $\omega=0$, which contradicts (1) above.

\hfill $\Box$

{\Small {\it Current Address}: Department of Mathematics and 
Statistics, University of Massachusetts, Amherst, MA 01003.
{\it E-mail:} wchen@@math.umass.edu}

\end{document}